\documentclass[a4paper]{elsarticle}
\usepackage[utf8]{inputenc}
\usepackage[T1]{fontenc}
\usepackage[english]{babel}
\usepackage{amsmath}
\usepackage{amssymb}
\usepackage{amsthm}
\usepackage{color}
 
\usepackage{times}
\usepackage{a4wide}
\usepackage{xspace}
\usepackage{enumerate}
\usepackage[colorlinks=true,citecolor=red,linkcolor=blue,%
urlcolor=RubineRed,pdfpagetransition=Blinds,pdftoolbar=false,pdfmenubar=false]{hyperref}

\usepackage{amsfonts}
\usepackage{eucal}
\usepackage{xspace}
\usepackage{color}
\usepackage{enumitem}
\usepackage{version}

\newcommand{\fin}{\hfill $\square$}

\newcommand{\field}[1]{\ensuremath{\mathbb{#1}}}
\newcommand{\C}{\field{C}\xspace}
\newcommand{\Z}{\field{Z}\xspace}

\newcommand{\Q}{\field{Q}\xspace}
\newcommand{\R}{\field{R}\xspace}
\newcommand{\N}{\field{N}\xspace}

\newcommand{\calD}{\mathcal{D}}
\newcommand{\calS}{\mathcal{S}}
\newcommand{\calH}{\mathcal{H}}
\newcommand{\calE}{\mathcal{E}}
\newcommand{\calN}{\mathcal{N}}
\newcommand{\calB}{\mathcal{B}}

\newcommand{\calC}{\mathcal{C}}
\newcommand{\calK}{\mathcal{K}}
\newcommand{\calA}{\mathcal{A}}
\newcommand{\calP}{\mathcal{P}}

\newcommand{\ens}[1]{ \left\{#1\right\} }
\newcommand{\paren}[1]{ \left( #1 \right) }
\newcommand{\vabsolut}[1]{ \left| #1 \right| }
\newcommand{\norma}[1]{ \left\| #1 \right\| }

\newtheorem{theorem}{Theorem}[section]
\newtheorem{corollary}[theorem]{Corollary}
\newtheorem{lemma}[theorem]{Lemma}

\newtheorem{proposition}[theorem]{Proposition}

\theoremstyle{definition}
\newtheorem{definition}{Definition}[section]
\newtheorem{remark}[definition]{Remark}

\numberwithin{equation}{section}
\usepackage[pagewise]{lineno}

\title{Sharp estimates for biorthogonal families to exponential functions associated to complex sequences  without gap conditions} 

\date{}

\author[EDAN]{Manuel \textsc{Gonz\'alez-Burgos}\fnref{fn1}\corref{cor1}}
\author[AMU]{Lydia \textsc{Ouaili}}
\ead{morgan.morancey@univ-amu.fr}

\fntext[fn1]{Supported by grant MTM2016-76990-P, Ministry of Economy and Competitiveness (Spain).}

\cortext[cor1]{Corresponding author}

\address[EDAN]{Dpto.~Ecuaciones Diferenciales y An\'alisis Num\'erico and Instituto de Matem\'aticas de la Universidad de Sevilla (IMUS), Facultad de Matem\'aticas, Universidad de Sevilla, C/ Tarfia S/N, 41012 Sevilla, Spain.  Email address: \texttt{manoloburgos@us.es}}

\address[AMU]{Aix Marseille Universit\'e, CNRS, Centrale Marseille, I2M, UMR 7373, 39 rue F. Joliot Curie 13453 Marseille, France. Email address: \texttt{ouaili.lydia@gmail.com}}

\begin{document}

\begin{abstract}
The general goal of this work is to obtain upper and lower bounds for the $L^2$-norm of biorthogonal families to complex exponential functions associated to sequences $\ens{ \Lambda_k}_{k \ge 1} \subset \C$ which satisfy appropriate assumptions but without imposing a gap condition on the elements of the sequence. As a consequence, we also present new results on the cost of the boundary null controllability of parabolic systems at time $T > 0$. In this case, the eigenvalues of the generator of the $C_0$-semigroup associated to this parabolic system accumulate, do not satisfy a gap condition and can develop a positive minimal time for the null controllability.
\end{abstract}

\maketitle


\section{Introduction and main results}\label{s1}

In the last years, an increasing number of authors have addressed the problem of the null controllability of coupled parabolic systems with less controls than equations (see~\cite{AKBDGB}, \cite{FCGBdT}, \cite{AKBGBdTJMPA}, \cite{LZ},...). One of the most important problems in this framework is obtaining necessary and sufficient conditions that allow the system to be controlled with a reduced number of distributed or boundary controls. 

Another important problem is the study of the dependence of the so-called control cost with respect to the final observation time $T > 0$, when $T$ is small enough and the corresponding null controllability result holds at any time $T > 0$. Regarding this latter problem, we highlight the works~\cite{FR}, \cite{FR2}, \cite{Seidman}, \cite{Gui}, \cite{Hansen}, \cite{FCZ}, \cite{Miller}, \cite{tenenbaum}, \cite{BBGBO}, \cite{cindea}, \cite{L}, \cite{CMV1}, \cite{CMV2}, etc., where the authors study an estimate of the control cost $\mathcal{K}(T)$ (for the definition, see~\eqref{controocost}) in the case of scalar parabolic problems (problems that, under general assumptions, are null controllable for any $T > 0$). 
Most of the previous works uses the moment method in order to obtain an estimate of the control cost. 

In order to solve both problems, a classical tool in Control Theory is the use of biorthogonal families to appropriate sequences of exponentials in $L^2 (0, T; \C)$ and, to be precise, sharp estimates on the $L^2$-norm of the elements of the biorthogonal family. We will provide more details in what follows. 

Given $\ens{\Lambda_k}_{k \ge 1} \subset \C$, a complex sequence of pairwise distinct elements, we will use the following notation:
	\begin{equation}\label{f0}
e_k (t) = e^{- \Lambda_k t}, \quad \forall t \in (0, T),
	\end{equation}
where $T >0$ is fixed. With this notation, we define


\begin{definition} \label{biorth}
Let $\Lambda=\{\Lambda_k\}_{k\geq1} \subset \C$ be a complex sequence and $T>0.$ We say that the family of functions $\{q_k\}_{k\geq1} \subset L^2(0,T; \mathbb{C})$ is a biorthogonal family to the sequence of complex exponentials $\ens{e_k}_{k \ge 1} $ in $L^2(0,T; \mathbb{C} )$, if for every $k,n \in \mathbb{N}$, one has
	$$
\int^T_0 e_k (t) \, \overline{q}_n (t) \,dt = \delta_{kn} ,
	$$
where the function $e_k$ is given in~\eqref{f0}. 
\end{definition}

Given $T > 0$, the general objective of this work is 
\begin{enumerate}
\item to analyze the existence of families $\ens{q_k}_{ k \ge 1 }$ biorthogonal to $\ens{e_k }_{k \ge 1} $ in $L^2 (0, T ; \C)$ for general sequences $\Lambda = \{ \Lambda_k \}_{k \geq 1} \subset \C$ of complex numbers pairwise distinct ($e_k$ is the exponential function defined in~\eqref{f0}) and
\item to obtain sharp and explicit estimates of $\norma{q_k}_{L^2 (0,T; C)}$ with respect to $T$, $\Lambda_k $ and some appropriate parameters associated to the sequence $\Lambda$.  
\end{enumerate}

Biorthogonal families play a crucial role in the moment method. This method was developed by Fattorini and Russell (see~\cite{FR} and~\cite{FR2}) to study the boundary null controllability of the heat equation and uses in a key way the existence and estimates of biorthogonal families to $\ens{e_k}_{k \ge 1}$. 

In~\cite{FR} the authors provide an approach that allows to construct biorthogonal families $\ens{q_k}_{ k \ge 1 }$ to the sequence $\ens{e_k}_{k \ge 1}$ in $L^2 (0, T)$ ($ T > 0$) with explicit bounds of the $L^2$-norm of $q_k$ with respect to the final time $T > 0$. To be precise, for increasing sequences $\Lambda = \{ \Lambda_k \}_{k \geq 1} \subset \R$ satisfying
	\begin{equation}\label{f0_0}
\Lambda_k \in (0, \infty), \quad \Lambda_k = A (k + \omega)^2 + o(k), \quad \forall k \ge 1,
	\end{equation}
with $ A >0 $ and $\omega \in \R $, there exist $\calC_0, \tau_0 \in (0, \infty)$ and a family $\ens{q_k}_{k \ge 1} $ biorthogonal to $\ens{e_k}_{k \ge 1}$ in $L^2 (0, T)$ ($ T > 0$) such that
	\begin{equation}\label{f0_1}
\norma{q_k}_{L^2(0, T)} \le \calC_0 e^{ \calC_0 \paren{\sqrt{ \Lambda_k} + \frac{1}{T}} } , \quad \forall T \in (0, \tau_0), \quad \forall k \ge 1,
	\end{equation}
(see for instance~\cite{FR} and~\cite{Miller}). 

As a consequence of inequality~\eqref{f0_1}, in~\cite{FR}, the authors prove that the one-dimensional heat equation 
	\begin{equation} \label{problem1}
	\left\{
	\begin{array}{ll}
\partial_ty - \partial_{xx} y = 0   & \hbox{in } (0,T)\times (0, L), \\
	\noalign {\smallskip}
y(\cdot,0) = v , \quad y(\cdot, L )=0   &\hbox{on } (0,T), \\ 
	\noalign{\smallskip}
y(0 ,\cdot)=y_0   &\hbox{in }  (0, L ),
	\end{array}
	\right.
	\end{equation}
($ L > 0$) is null controllable in $H^{-1} (0 , L)$ at any time $T > 0$ with controls $v \in L^2(0,T)$. In fact, they prove the existence of a constant $C_0$ (only depending on $L$) such that for any $y_0 \in H^{-1} (0 , L) $ there exists a control $v \in L^2(0,T)$ satisfying
	\begin{equation*}
\norma{v}_{ L^2(0,T) } \le C_0 e^{ \frac {C_0}{T}} \norma{y_0 }_{H^{-1 } (0, L) },
	\end{equation*}
and such that the corresponding solution to~\eqref{problem1}  satisfies $y(T, \cdot ) = 0$ in $ (0, L)$. Thus, the set 
	$$
\calC_{T}\paren{y_0} :=\left \{ v \in L^2(0,T) : y(T, \cdot )=0 \hbox{ in } (0, L), \ y \hbox{ solution of } \eqref{problem1} \right \},
	$$
is non empty and we can define the so-called \textit{control cost} of system~\eqref{problem1} at time $T >0$:
	\begin{equation} \label{controocost}
\calK(T) := \sup_{\|y_0\|_{H^{-1}}=1}  \inf_{v \in \calC_{T} \paren{y_0} } \|v\|_{L^2(0,T)}. 
	\end{equation}
Therefore, for system~\eqref{problem1}, one has
	\begin{equation}\label{f0_2}
\calK(T) \le C_0 e^{ \frac{C_0}{T}} , \quad \forall T >0,
	\end{equation}
for a positive constant $C_0$ only depending on $L$.

In the framework of $N$-dimensional scalar parabolic problems,~\cite{Seidman} and~\cite{FCZ} give an estimate of the cost $\calK ( T )$ similar to~\eqref{f0_2} using different approaches: In~\cite{Seidman} the authors use the exact controllability of the wave equation to prove inequality~\eqref{f0_2} for the null-controllability of the heat equation. In~\cite{FCZ}, inequality~\eqref{f0_2}  is deduced from appropriate global Carleman inequalities for general parabolic operators. 

Estimate~\eqref{f0_2} is known to be optimal thanks to the work~\cite{Gui}: under assumption~\eqref{f0_0}, there exists a positive constant $\calC_1$ such that for any sequence $ \{q_k\} _{k\geq1} \subset L^2(0,T) $  biorthogonal to $\{ e_k \}_{k\geq1}$ in $L^2(0,T)$, one has
	\begin{equation}\label{f0_3}
\norma{q_k}_{L^2(0, T)} \ge \frac{M(k)}{\sqrt{T}}  e^{ \frac{\calC_1}{T}}  , \quad \forall T > 0, \quad \forall k \ge 1,
	\end{equation}
where $M (k) $ is a positive constant only depending on $k$ and $L$. In particular, inequality~\eqref{f0_3} implies the existence of new positive constants $C_1$ and $\tau_1$ (only depending on $L$) such that the control cost for system~\eqref{problem1} satisfies
	\begin{equation}\label{f0_4}
\calK(T) \ge C_1 e^{ \frac{ C_1}{T}} , \quad \forall T \in (0, \tau_1). 
	\end{equation}

The existence of biorthogonal families $\ens{q_k}_{k \ge 1}$ to sequences of exponentials $\ens{e_k}_{k \ge 1}$ ($e_k$ is the function given in~\eqref{f0}) and estimates~\eqref{f0_1},~\eqref{f0_2},~\eqref{f0_3} and~\eqref{f0_4} strongly depend on the properties of the sequence $\Lambda = \{ \Lambda_k \}_{k \geq 1}$. Our next objective is to provide some general properties for real or complex sequences $\Lambda $ appearing in the literature which assure the existence of sequences $\ens{ q _k }_{k \ge 1} $ biorthogonal to $\ens{e_k}_{k \ge 1}$ in $L^2 (0, T; \C)$ ($ T > 0 $) satisfying~\eqref{f0_1} or~\eqref{f0_3}.

As said before, the first results on existence and estimates of families $\ens{q_k}_{k \ge 1}$ biorthogonal to sequences of exponentials $\ens{e_k}_{k \ge 1}$ was proved in~\cite{FR},~\cite{FR2} and~\cite{Gui} (see also~\cite{Korevaar},~\cite{Seidman00}, \cite{Miller}, \cite{TeTu}, \cite{L0} and~\cite{L1}) for increasing real sequences that satisfy~\eqref{f0_0}.

The results on existence of biorthogonal families to $\ens{ e_k }_{k \ge 1}$ has been extended to the complex case in~\cite{Hansen},~\cite{AKBGBdTJMPA},~\cite{AKBGBdTJFA} and~\cite{BBGBO}. In the three first works, the authors prove the existence of biorthogonal sequences $\ens{q_k}_{k \ge 1}$ under general assumptions on the sequence $\Lambda$ and prove appropriate estimates of $\| q_k \|_{L^2 (0, T; \C )}$ (in the case of~\cite{AKBGBdTJFA}, the authors prove the results without imposing gap conditions on the sequence $\Lambda$). Nevertheless, in these works the authors use a technique that does not allow them to obtain an explicit dependence of this estimate with respect to the final time $T > 0$. Therefore, inequality~\eqref{f0_2} cannot be deduced from these works (for the details, see~\cite{Hansen} and~\cite{AKBGBdTJFA}).

Let us describe the result on existence and estimates of biorthogonal families to complex exponentials proved in~\cite{BBGBO}. One has:


\begin{theorem}[\cite{BBGBO}] \label{Olive}
Let $\Lambda = \ens{\Lambda_k}_{k \ge 1} \subset \C$ be a sequence satisfying assumptions~\ref{item1}--\ref{item5}, in Definition~\ref{d1}, the gap condition
	\begin{equation}\label{f29}
\inf_{k , n \ge 1 : k \not= n} \vabsolut{  \Lambda_{k}  -  \Lambda_n  }   > 0,
	\end{equation}
and
	$$
\vabsolut{p \sqrt{r} - \calN (r) } \le \alpha, \quad \forall r > 0,
	$$
($\mathcal{N}$ is the counting function associated with the sequence $\Lambda $, defined in~\eqref{counting}), for some parameters $\beta \in [0, \infty)$, $ \rho , p, \alpha \in (0, \infty) $ and $q \in \N$. Then, there exists $T_0 > 0$ such that  for every  $T \in (0 , T_0)$, there exists a sequence of $\mathbb{C}$-valued functions 
	$$
\{q_k\}_{k\geq1} \subset L^2(0 ,T ; \C )
	$$
biorthogonal to the exponentials $\{ e_k \}_{k\geq1}$ in $L^2(0 ,T; \C )$, $e_k $ given in~\eqref{f0}, which, in addition, satisfies~\eqref{f0_1} for a positive constant $\calC_0$ independent of $k $ and $T$. 
\end{theorem}

The previous result can be applied to a large range of scalar and coupled parabolic problems. It assures that the  system under consideration is null controllable at any time $T > 0$. In addition, Theorem~\ref{Olive} provides the inequality~\eqref{f0_2} for the control cost $\mathcal{K}(T)$ as in the scalar case ($C_0$ is a positive constant). In order to prove Theorem~\ref{Olive}, the authors use the Fourier transform with the help of the Paley-Wiener theorem (see~\cite{BBGBO} for the details).

The existence of biorthogonal families to real exponentials that satisfy~\eqref{f0_1} and~\eqref{f0_4} has been also treated by some authors with assumptions on the sequence $\Lambda $ different from~\eqref{f0_0}. In~\cite{cindea}, the authors consider a real increasing sequence $\Lambda = \ens{ \Lambda_k }_{ k \ge 1}$ that is given as 
	\begin{equation*}
\Lambda = \ens{ \lambda_k^{(1)}}_{ k \ge 1} \cup \ens{ \lambda_k^{(2)}}_{ k \ge 1}, 
	\end{equation*}
with $\ens{ \lambda_k^{(1)}}_{ k \ge 1} $ and $ \ens{ \lambda_k^{(2)}}_{ k \ge 1} $ two increasing sequences of positive real numbers satisfying 
	\begin{equation}\label{f28}
	\left\{
	\begin{array}{l}
\displaystyle \vabsolut{  \lambda_k^{(1)} - \frac{1}{\pi^2_1 }  k^2 } \le c_1 k, \quad \vabsolut{  \lambda_k^{(2)} - \frac{1}{ \pi^2_2} k^2 } \le c_1 k, \quad \forall k \ge 1, \\
	\noalign{\smallskip}
\displaystyle \inf_{n \ge 1} \vabsolut{ \sqrt{ \lambda_{k}^{(2)}} - \sqrt{ \lambda_n^{(1)}} }  \ge \frac rk , \quad \forall k \ge 1,
	\end{array}
	\right.
	\end{equation}
and the strong gap condition 
	\begin{equation}\label{f28b}
\displaystyle \sqrt{ \lambda_{k+1}^{(1)}} - \sqrt{ \lambda_{k}^{(1)}}  \ge c_2 , \quad \sqrt{ \lambda_{k+1}^{(2)}} - \sqrt{ \lambda_{k}^{(2)}}  \ge c_2, \quad \forall k \ge 1, 
	\end{equation}
for some positive constants $\pi_1$, $\pi_2$, $c_1, c_2$ and $r$. For this class of sequences, the authors prove the existence of a sequence $\ens{q_k}_{k \ge 1 } \subset L^2 (0,T)$ ($T > 0 $ is given) biorthogonal to $\ens{e_k}_{k \ge 1 } $ ($e_k$ given in~\eqref{f0}) in $L^2(0,T)$ which satisfies~\eqref{f0_1} for a positive constant $\calC_0$  independent of $k$ and $T$ and uniform for the class of sequences $ \Lambda$ satisfying the previous assumptions.

In~\cite{CMV1} and~\cite{CMV2}, the authors again consider increasing real positive sequences $\Lambda = \{ \Lambda_k \}_{k \geq 1} \subset \R$ satisfying a ``global gap condition'':
	\begin{equation}\label{f27}
\gamma_{0} \le \sqrt{\Lambda_{k +1}} - \sqrt{\Lambda_k} \le \gamma_{1}, \quad \forall k \ge 1 ,
	\end{equation}
and an ``asymptotic gap condition'':
	\begin{equation*}
\gamma_{0}^\star \le \sqrt{\Lambda_{k +1}} - \sqrt{\Lambda_k} \le \gamma_{1}^\star , \quad \forall k \ge N ,
	\end{equation*}
where $ N  $ is a positive integer and $\gamma_{0}, \gamma_{1}, \gamma_{0}^\star, \gamma_{1}^\star \in (0, \infty)$ are such that $0 < \gamma_{1}^\star - \gamma_{0}^\star < \gamma_{1} - \gamma_{0}$. Under these assumptions on $\Lambda$ the authors obtain general and precise upper and lower bounds for biorthogonal families as~\eqref{f0_1} or~\eqref{f0_3}, paying attention to the dependence of the constant $\calC_0$ and $\calC_1$ with respect to the parameters $\gamma_{0}$, $\gamma_{1}$, $\gamma_{0}^\star$ and $ \gamma_{1}^\star$.

It is interesting to observe that in all the previous works the authors impose conditions on the sequence $\Lambda = \ens{\Lambda_k}_{k \ge 1} \subset \C $, which, in particular, imply that it satisfies the gap condition~\eqref{f29}. This is easy to check for increasing real sequences fulfilling condition~\eqref{f27} and can be checked for the class of sequences $\Lambda = \ens{ \lambda_k^{(1)}}_{ k \ge 1} \cup \ens{ \lambda_k^{(2)}}_{ k \ge 1} $ under the assumptions of~\cite{cindea}. In fact, if the sequence $\Lambda  $ satisfies~\eqref{f28} and~\eqref{f28b}, then $\Lambda $ also satisfies the assumptions in Theorem~\ref{Olive} (see Proposition~\ref{p3} and Remark~\ref{r6}). To our knowledge, assumptions~\ref{item1}--\ref{item5} and \eqref{f29} are the most general hypotheses on the sequence $\Lambda$ that guarantee the existence of a family $ \{q_k\}_{k\geq1} \subset L^2(0 ,T ; \C )$ biorthogonal to the exponentials $\{ e_k \}_{k\geq1}$ in $L^2(0 ,T; \C )$, $e_k $ given in~\eqref{f0}, that satisfies~\eqref{f0_1} for a positive constant $\calC_0$ independent of $k $ and $T$.

The work~\cite{L} is of special relevance because in it, the author studies the cost of the controllability of the one-dimensional heat equation with a pointwise control at point $x_0$ and, in this framework, there might exist a positive minimal time of null-controllability $T_0 \in [0, \infty]$ (which depends on $x_0$ and could take any arbitrary value in $[0, \infty]$, see~\cite{Dolecki}). In this work the eigenvalues satisfy~\eqref{f29} and the minimal time comes from the action of the control. In particular, the author proves that, if $T_0 > 0$, the cost of the controllability at time $T > T_0$ when $T$ is close to $T_0$, may explode in any arbitrary way.

As said before, the analysis of the control cost in the framework of the controllability of coupled parabolic systems has been addressed in~\cite{BBGBO}. As in the previous works, the authors impose appropriate assumptions on the sequence $\Lambda$ which include a gap condition on the terms of the sequence (see Theorem~\ref{Olive} and~\eqref{f29}). In particular, conditions in Theorem~\ref{Olive} assure that the system under study is null controllable at any positive time $T $ and the control cost $\calK(T)$ satisfies~\eqref{f0_2} for a positive constant $C_0$. 

In the framework of the controllability of non-scalar parabolic problems, new phenomena associated with the vectorial nature of the problem arise (hyperbolic phenomena): minimal time of null controllability and dependence of the controllability result on the position of the control domain (see~\cite{AKBGBdTJFA},~\cite{AKBGBdT2}, \cite{Samb},~\cite{O},...). This minimal time may come from the control action itself (as in~\cite{Dolecki} and~\cite{L}) or from the condensation index of the sequence of eigenvalues of the generator of the semigroup associated to the system (see~\cite{AKBGBdTJFA}). In this latter case, the sequence $\Lambda$, in general, does not satisfy the gap condition~\eqref{f29}. Let us provide more details in the case of systems with a minimal time which comes from the condensation index of the sequence.

Assume that the sequence $\Lambda = \ens{\Lambda_k}_{k \ge 1} \subset \C$ satisfies 
	\begin{equation}\label{f18}
	\left\{
	\begin{array}{l}
\Lambda_i \not= \Lambda_k , \quad \forall i,k \in \N \hbox{ with } i \not= k, \\
	\noalign{\smallskip}
\displaystyle \Re \left( \Lambda _{k}\right) \geq \delta \left\vert \Lambda _{k}\right\vert >0, \quad \forall k\geq 1, \quad \hbox{and} \quad \sum_{k\geq 1}\frac{1}{\left\vert \Lambda _{k}\right\vert }<\infty ,
	\end{array}
	\right.
	\end{equation}
for a positive constant $\delta$. Observe that, in general, a sequence $\Lambda$ satisfying~\eqref{f18} does not fulfill the gap condition~\eqref{f29}.

Properties~\eqref{f18}  for the sequence $\Lambda$ imply that the family of exponentials $\ens{e_k}_{k \ge 1}$ is minimal\footnote{A sequence $\{ x_ k\}_{k \ge 1}$ in a Hilbert space $H$ is said to be minimal if it satisfies $x_n \not\in \overline{\hbox{span}} \, \{ x_k : k \not= n\} $ for any $ n \ge 1$.}
in $L^2(0,T; \C)$ for any $T > 0$ and, therefore, there exists a biorthogonal family $\ens{\widetilde q_k}_{k \ge 1}$ to $\ens{e_k}_{k \ge 1}$ in $L^2(0,T; \C)$ (see for instance~\cite{S}, \cite{Re}, \cite{AKBGBdTJMPA}, Theorem~4.1 in~\cite{AKBGBdTJFA},...). In addition, in~\cite{AKBGBdTJFA}, the authors prove that there exist two positive constant $C_1$ and $C_2$ (only depending on $\Lambda$ and $T$) such that
	\begin{equation}\label{f17}
C_1 \frac{\left|1 + \Lambda_k \right|^2 }{\left| \Lambda_k \right|} W_k  \le \| \widetilde q_k \|_{L^2 (0,T; \C)} \le C_2 \frac{\left|1 + \Lambda_k \right|^2 }{\left| \Lambda_k \right|} W_k,
	\end{equation}
where $C_1$ and $C_2$ are positive constant depending on $T$ and $W_k$ is the infinite Blaschke product given by
	\begin{equation*}
W_k = \frac{1 }{2  \Re ( \lambda_k) }\prod_{\substack{ n \ge 1  \\ n \not= k }} \vabsolut{\frac{\Lambda_n + \Lambda_k }{\Lambda_n - \Lambda_k }} \, . 
	\end{equation*}
Nevertheless, the authors do not provide an explicit dependence of the constants $C_1$ and $C_2$ in~\eqref{f17} with respect to the final time $T > 0$. This is due to the method used by the authors to prove~\eqref{f17}: these inequalities are first obtained in $L^2(0, \infty; \C)$ ($T = \infty$) and, then, proved in $L^2(0,T; \C)$ ($T \in (0, \infty)$) after a contradiction argument (see~\cite{AKBGBdTJFA} for the details).

From inequality~\eqref{f17}, among other properties, in~\cite{AKBGBdTJFA}, the authors prove a general result of null controllability for abstract parabolic problems that develop a minimal time $T_0 \in [0, \infty]$ of controllability: the system is null-controllable at any time $T>T_0$ and not null-controllable for $T<T_0$. This minimal time is related to the Bernstein's condensation index of the sequence of eigenvalues $\Lambda = \ens{\Lambda_k}_{k \ge 1}$ of the generator of the semigroup (see~\cite{AKBGBdTJFA} and~\cite{BBM} for further details). 

Let us now revisit some one-dimensional non-scalar parabolic systems with a generator whose sequence of eigenvalues $\Lambda = \ens{\Lambda_k}_{k \ge 1} \subset \C$ satisfies~\eqref{f18} and not inequality~\eqref{f29}. To this end, we consider a boundary controllability problem for the generic $2 \times 2$ system
	\begin{equation} \label{generic}
	\left\{ 
	\begin{array}{ll}
\partial_{t} y + L y = 0 & \text{in }   (0,T) \times (0, \pi)  , \\ 
 	\noalign{\smallskip}
y( \cdot, 0 ) = B v, \quad y( \cdot, \pi ) = 0 &\hbox{on  } (0,T),\\
  	\noalign{\smallskip}
y(0, \cdot) = y_0 & \hbox{in  } (0,\pi),\\
	\end{array}
	\right. 
	\end{equation}
where $L$ is a second order elliptic operator, with domain $D (L) = H^2(0, \pi; \R^2 ) \cap H_0^1 (0, \pi; \R^2 )$, $y_0 \in H^{-1} \paren{ 0,1,\R^2 }$ is the initial datum, $B\in \R^2$ is the control vector and $v \in L^2 (0,T )$ is a scalar control. 

The null controllability properties of the first example has been analyzed in~\cite{AKBGBdTJFA}. We consider system~\eqref{generic} with $L = L_1 = - \left(D_1 \partial_{xx} + A _1 \right)$, with domain $D (L_1) = H^2(0, \pi; \R^2 ) \cap H_0^1 (0, \pi; \R^2 )$, and 
	$$
D_1 := \hbox{diag} \, (1, d), \ d>0, \ d \not= 1, \quad \hbox{and} \quad A_1 := \left( 
	\begin{array}{cc}
0 & 1 \\ 
0 & 0
	\end{array}
\right),
	$$
(see system~\eqref{f32}). Observe that the sequence of eigenvalues associated to the operator $L_1 $ is $\Lambda^{(1)} = \ens{k^2}_{k \ge 1} \cup \ens{d k^2}_{k \ge 1} $. If $\sqrt{d } \not\in \Q$ (and this condition is necessary for the approximate controllability at time $T>0$ of the system~\eqref{generic} with the previous data, i.e., system~\eqref{f32}), the sequence $\Lambda^{(1)}$ can be rearranged as an increasing sequence $\Lambda^{(1)} = \ens{\Lambda^{(1)}_k}_{k \ge 1} \subset \R$ that fulfills property~\eqref{f18}. It is clear that $\Lambda^{(1)} $ does not satisfy, in general, the gap condition~\eqref{f29}. As a consequence, system~\eqref{f32} has a minimal time $ T_0 = T_0 (d) \in [0, \infty ] $ which, for some $d$, with $ \sqrt{d } \not\in \Q$, is positive. Therefore, the system is not null controllable at time $T$ when $T < T_0$ (see~\cite{AKBGBdTJFA} for the details).

The controllability properties of our second example has been analyzed in~\cite{O}. Let us consider system~\eqref{generic} with

	\begin{equation}\label{L2}
L = L_2 := \left(
	\begin{array}{cc}
- \partial_{xx}  &0\\0& - \partial_{xx} + Q
	\end{array}
	\right) ,
\quad D(L_2 ) = H^2(0, \pi;\R^2) \cap H_0^1(0,\pi; \R^2) ,
	\end{equation} 
with $ Q \in L^2(0,\pi )$. In this case, the sequence of eigenvalues of the vectorial operator $L_2 $ is given by $\Lambda^{(2)} = \ens{  k^2  }_{k \ge 1} \cup \ens{ \lambda_k^{(2)} }_{k \ge 1} \subset \R $, where $\ens{ \lambda_k^{(2)} }_{k \ge 1}$ is the sequence of eigenvalues of the operator $ - \partial_{xx} + Q$ with domain $H^2(0,1) \cap H_0^1(0,\pi)$. When $Q \in L^2(0, \pi)$ satisfies
	\begin{equation*}
\int_0^\pi Q(x) \, dx = 0,
	\end{equation*}
then 
	\begin{equation*}
\lambda_k^{(2)} =  k^2 + \varepsilon_k , \quad \forall k \ge 1 ,
	\end{equation*}
with $\ens{ \varepsilon_k }_{k \ge 1} \in \ell^2$. In particular, $\lim \varepsilon_k = 0$ and $ \Lambda^{(2)} $ does not fulfill the gap condition~\eqref{f29}. Assume that $ \lambda_k^{(2)} \not=  n^2 $ for any $k,n \ge 1$ (that, in fact, is a necessary condition for the approximate controllability of system~\eqref{generic} with $L = L_2 $, see~\cite{O} and Section~\ref{s5}). In this case, the sequence $\Lambda^{(2)} $ satisfies property~\eqref{f18}. Again, system~\eqref{generic} has a minimal time $ T_0 = T_0 (Q) \in [0, \infty ] $ and there exists coefficients $Q \in L^2(0,\pi)$ such that $T_0 (Q) > 0$. Thus, the system is not null controllable at time $T$ when $T < T_0$ (see~\cite{O} and Section~\ref{s5} for the details).

Let us consider a third example of non-scalar parabolic system. In~\cite{GB-SN} the authors study the boundary null controllability of a phase field system of Caginalp type which is a model describing the transition between the solid and liquid phases in solidification/melting processes of a material occupying the  interval $(0, \pi)$. For that purpose, they consider the nonlinear system
	\begin{equation}\label{PFSy}
	\left\{
	\begin{array}{ll}
\displaystyle {\theta}_t - \xi{\theta}_{xx} + \dfrac{1}{2}\rho\xi{\phi}_{xx} + \dfrac{\rho}{\tau}{\theta} = f({\phi})	& \mbox{in } (0,T) \times (0, \pi), 
	\\
	\noalign{\smallskip}
\displaystyle {\phi}_t - \xi{\phi}_{xx} - \dfrac{2}{\tau}{\theta} = -\frac 2 \rho f({\phi}) & \mbox{in }  (0,T) \times (0, \pi) , 
	\\
	\noalign{\smallskip}
\displaystyle {\theta}(\cdot, 0) = v,\ {\phi} (\cdot, 0) = c,\ {\theta}(\cdot, \pi )=0 , \ {\phi}(\cdot, \pi ) = c & \mbox{on }  (0,T),	 
	\\
	\noalign{\smallskip}
\displaystyle {\theta}(0, \cdot ) = {\theta}_0, \  {\phi}(0, \cdot ) ={\phi}_0 & \mbox{in }  (0,\pi),	\end{array}
	\right. 
	\end{equation}
where: ${{\theta} = {\theta}(t, x)}$ is the temperature of the material; ${{\phi} = {\phi} (t, x)}$ is the phase-field function used to identify the solidification level of the material; ${c \in \{ - 1,0,1\}}$; ${f}$ is the nonlinear term which comes from the derivative of the classical regular double-well potential $W$:
	$$ 
\displaystyle f({\phi}) = -\frac{\rho}{4\tau}\left( {\phi}-{\phi}^3 \right) . 
	$$
On the other hand, 
$\rho > 0$, $ \tau > 0$ and $ \xi > 0$ are, resp., the latent heat, a relaxation time and the thermal diffusivity. Finally, ${v \in L^2(0,T)}$ is the control function, and ${ \theta_0,  \phi_0}$ are the initial data.

The null controllability property of the nonlinear system~\eqref{PFSy} depends on the coefficients $\rho$, $\tau$ and $ \xi $. This property is obtained from the corresponding one of the linear version of~\eqref{PFSy} around the constant trajectory $(0, c)$ (see~\cite{GB-SN} for more details). This linear system is as system~\eqref{generic} with $y = (\theta , \phi)$ and $L = L_3$ given by
	\begin{equation}\label{f78}
	\left\{
	\begin{array}{l}
\displaystyle L = L_3 := - D_2  \partial_{xx} + A_ 2, \quad \hbox{with } \\
	\noalign{\smallskip}
\displaystyle D = D_2 := \left(
	\begin{array}{cc}
\xi & -\dfrac{1}{2} \rho \xi \\
	\noalign{\smallskip}
0 &  \xi
	\end{array}
	\right),
\quad
A = A_2 := \left(
	\begin{array}{cc}
\dfrac{\rho}{\tau} & -\dfrac{\rho}{2\tau} \\
	\noalign{\smallskip}
-\dfrac{2}{\tau} & \dfrac{1}{\tau}
	\end{array}
	\right),
\quad
B =
	\left(
	\begin{array}{cc}
1 \\ 0
	\end{array}
	\right). 
	\end{array}
	\right.
	\end{equation} 
In this case the sequence of eigenvalues of the operator $L_3 $, with domain $D(L_3 ) = H^2(0,1;\R^2) \cap H_0^1(0,1; \R^2)$, is given by $ \Lambda^{(3)} = \ens{\lambda_k^{(3,1)}, \lambda_k^{(3,2)} }_{k\geq1}$ with
	\begin{equation}\label{f78'}
\lambda_k^{(3,1)} = \xi k^2 + \dfrac{\rho + 1}{2\tau} - r_k ,\quad \lambda_k^{(3,2)} = \xi k^2 + \dfrac{\rho + 1}{2\tau} + r_k, \quad \forall k \ge 1,			
	\end{equation}
where $r_k$ is given by
	\begin{equation}\label{f78''}
r_k:= \sqrt{ \dfrac{\xi\rho}{\tau} k^2 + \left(\dfrac{\rho+1}{2\tau}\right)^2 } , \quad \forall k \ge 1.
	\end{equation}

If $\lambda_k^{(3,1)} \not= \lambda_n^{(3,2)} $ for any $k,n \ge 1$ (which in fact is a condition equivalent to the approximate controllability of the linear system~\eqref{generic} with $L = L_3 $), the sequence $\Lambda^{(3)}$ can be rearranged in such a way that $\Lambda^{(3)} = \ens{\Lambda^{(3)}_k}_{k \ge 1} $ is an increasing sequence that satisfies~\eqref{f18} for appropriate $\delta >0 $. However, if for some integer $j \ge 1$ one has
	\begin{equation}\label{f77}
\xi = \dfrac{1}{j^2}\dfrac{\rho}{\tau},
	\end{equation}
then, the eigenvalues of $L_3$ concentrate and one has
	$$
\inf_{k \ge 1} \paren{\Lambda^{(3)}_{k + 1} - \Lambda^{(3)}_k} = 0,
	$$
and condition~\eqref{f29} does not hold (see~\cite{GB-SN} and Section~\ref{s5} for the details). Therefore, we have another system where the associated sequence of eigenvalues does not satisfy the gap condition~\eqref{f29}.

\begin{remark}\label{r0}
It is interesting to observe that the objective of the work~\cite{GB-SN} is to study the exact boundary controllability to constant trajectories at time $T$, $T > 0$ arbitrary, of the nonlinear system~\eqref{PFSy}. To this end, the authors follow a technique developed in~\cite{LTT}. This methodology consists of obtaining a null controllability result at time $T$ for system~\eqref{generic}, with $L = L_3$, and an estimate of the cost of fast controls like~\eqref{f0_2}. In order to obtain inequality~\eqref{f0_2} for the linear version of system~\eqref{PFSy}, the authors assume the condition
	\begin{equation*}
\xi \not= \dfrac{1}{j^2}\dfrac{\rho}{\tau}, \quad \forall j \in \N . 
	\end{equation*}
This condition is crucial in~\cite{GB-SN} because it assures that the sequence $\Lambda^{(3)}$ satisfies~\eqref{f29} and the conditions in Theorem~\ref{Olive}. Thus, system~\eqref{generic}, with $L = L_3$, is null controllable at time $T$ for any $T>0$ and the control cost $\calK(T) $ satisfies~\eqref{f0_2} for a positive constant $C_0$ only depending on $\rho$, $\tau$ and $\xi$. \fin
\end{remark}

We have seen three examples of sequences of eigenvalues satisfying~\eqref{f18} and for which the gap condition~\eqref{f29} fails. The corresponding parabolic systems could have a positive minimal time of null controllability $T_0$ and the system is not null controllable at time $T$ when $T \in (0, T_0)$. Even if $T_0 = 0$, it is not clear that the control cost of the associated system fulfills inequality~\eqref{f0_2} or inequality~\eqref{f0_4} and this is an open problem.

In order to obtain sharp estimates of the control cost $\calK(T)$ associated to non-scalar parabolic systems, it is very important to prove sharp estimates for biorthogonal families to the exponentials associated to the corresponding sequence of eigenvalues of the generator when this sequence does not satisfy a gap condition. This is the objective of this work: Given a complex sequence $ \Lambda = \ens{\Lambda_k}_{k \ge 1}$ satisfying appropriate assumptions and such that inequality~\eqref{f29} does not hold, is there a biorthogonal family $\ens{q_k}_{k \ge 0} $ to $\ens{e_k }_{k \ge 1} $ in $L^2 (0,T ; \C)$ ($e_k$ is given in~\eqref{f0}) satisfying an appropriate estimate for $\norma{q_k}_{L^2 (0,T ; \C)}$ which, in particular, provides an estimate of the control cost $\calK (T)$? Understanding the behavior of the control cost $\calK (T)$ for general systems as~\eqref{generic} would allow us to extend the null controllability result in the one-dimensional case to some parabolic systems in any dimension (see for instance \cite{BBGBO,allonsius:hal-01827044}) and to some nonlinear parabolic equations using the method of Liu, Takahashi and Tucsnak introduced in \cite{LTT} (see for instance,~\cite{GB-SN} and~\cite{tesisLydia}). 

%


Let us  now present the main results of this work. To this aim, let us first introduce the class of complex sequences we will work with throughout this work:


\begin{definition}\label{d1}

Let $\Lambda =\left \{ \Lambda_k \right \}_{k \geq 1}$ be a complex sequence and let us fix constants $\beta \in [0, \infty)$, 
	$$ 
\rho , p_0, p_1, p_2, \alpha \in (0, \infty) 
	$$ 
and $q \in \N$. We say that the sequence $\Lambda$ is in the class $ \mathcal{L}(\beta,\rho,q,p_0, p_1, p_2,  \alpha)$, if the following properties hold:
\begin{enumerate}[label={(H\arabic*)},align=left]
\item \label{item1} $\Lambda_k\neq \Lambda_n$ for all $n,k\in \mathbb{N}^* $ with $n\neq k$; 
\item \label{item2}$\Re(\Lambda_n)>0$ for every $n\geq1$;
\item \label{item3} $\vert \Im(\Lambda_n) \vert \leq \beta \sqrt{\Re(\Lambda_n)}$, for any $ n\geq1$;
\item \label{item4} $\{ \Lambda_n\}_{n\geq1}$ is nondecreasing in modulus, i.e., $\vert \Lambda_n\vert \leq \vert \Lambda_{n+1} \vert$, for any $ n\geq1$;
\item \label{item5} $ \rho \left| k^2 -n^2 \right| \leq \left\vert \Lambda_{k}-\Lambda_{n} \right\vert $ for any $ n,k\geq 1: \vert k-n \vert\geq q; 	$
\item \label{item7} $p_1,  p_2 \ge p_0 $ and one has
	\begin{equation*}
 -\alpha+p_1\sqrt{r}  \leq\mathcal{N}(r)  \leq \alpha+ p_2 \sqrt{r}, \quad \forall r>0,
	\end{equation*}
where $\mathcal{N}$ is the counting function associated with the sequence $\Lambda $, defined by 
	\begin{equation}\label{counting}
\mathcal{N}(r)= \# \ens{ k :\, \vert \Lambda_{k}\vert \leq r}, \quad \forall r>0.
	\end{equation}
\end{enumerate}
\end{definition}


\begin{remark}\label{r1}
Observe that from the definition of the counting function $\mathcal{N}$ (see~\eqref{counting}) associated with the sequence $\Lambda = \left \{ \Lambda_k \right \}_{k\geq} \in \mathcal{L}(\beta,\rho,q,p_0, p_1, p_2,  \alpha)$ ($\beta \in [0, \infty)$, $ \rho , p_0, p_1, p_2, \alpha \in (0, \infty) $ and $q \in \N$ are given), we deduce the following properties: 
\begin{enumerate}
\item For any $r > 0$, one has
	\begin{equation*}
\mathcal{N}(r) = k \iff | \Lambda_k | \le r \quad \hbox{and} \quad | \Lambda_{k + 1} | > r. 
	\end{equation*}
\item If for some $ k_1,  k_2 \ge 1$ and $ r_1,  r_2 > 0$ one has $ \vabsolut{\Lambda_{ k_1}} \le  r_1 $ and $\vabsolut{\Lambda_{ k_2} } > r_2 $, then
	$$
k_1 \le \calN ( r_1 ) \quad \hbox{and} \quad k_2 \ge \calN ( r_2 ) +1.
	$$
We will use these properties throughout this work. \fin
\end{enumerate}
\end{remark}



\begin{remark}\label{r9}
The parameter $q \in \N $ in Definition~\ref{d1} plays an important role in this paper. Observe that in this work we are dealing with sequences $\Lambda$ that, in general, do not satisfy condition~\eqref{f29} and whose terms could condense. With condition~\ref{item5} and the parameter $q $ we mesure the maximal cardinal of the condensation groupings of the sequence $\Lambda$, that is to say, the maximal number of elements in $\Lambda$ around the term $\Lambda_k$ that do not satisfy~\ref{item5} and could condense. At the end of Section~\ref{s2} we will see, by means of an example, that the parameters $p_1$ and $p_2$ increasingly depend on $q$, even for real sequences $\Lambda$ that satisfies the gap condition~\eqref{f29}.  \hfill $\Box$
\end{remark}


We will see in Section~\ref{s2} that the class $ \mathcal{L}(\beta,\rho,q,p_0, p_1, p_2,  \alpha)$ includes sequences $\Lambda = \ens{ \Lambda_k }_{ k \ge 1}$ satisfying~\eqref{f0_0} or condition~\eqref{f27}, and sequences $\Lambda = \ens{ \lambda_k^{(1)}}_{ k \ge 1} \cup \ens{ \lambda_k^{(2)}}_{ k \ge 1}$ under assumptions~\eqref{f28} and~\eqref{f28b}. Also, it includes sequences that do not satisfy the gap condition~\eqref{f29} as $\Lambda = \ens{k^2}_{k \ge 1} \cup  \ens{d k^2}_{k \ge 1} $ ($\sqrt{d} \not\in \Q$), or $\Lambda = \ens{  k^2  } \cup \ens{  k^2 + \varepsilon_k } $, with $\ens{ \varepsilon_k }_{k \ge 1} \in \ell^2$, or the sequence considered in~\cite{GB-SN} (see~Remark~\ref{r0}).

We are now in a position to establish the first main result of this work. It reads as follows: 

\begin{theorem}\label{principal theorem} 
Let $\Lambda =\left \{ \Lambda_k \right \}_{k\geq}\subset\mathbb{C}$ be a sequence satisfying $\Lambda \in \mathcal{L}(\beta,\rho,q,p_0, p_1, p_2,  \alpha) $ with $\beta \in [0, \infty)$, $ \rho , p_0, p_1, p_2, \alpha \in (0, \infty) $ and $q \in \N$. Then, given $T>0$, there exists a family of complex functions 
	$$
\{q_k\} _{k\geq1} \subset L^2(0,T; \C),
	$$
biorthogonal to $\{e_k \}_{k\geq1}$ in $L^2(0,T; \C)$ ($e_k$ is given in~\eqref{f0}) which, in addition, satisfies
	\begin{equation} \label{bounds}
\displaystyle \|q_k \|_{L^2(0,T; \C)} \leq \calH_1(\rho, q, p_1,p_2)  \exp \left[C \paren{
1 +\calH_2 (\rho, q, p_1,p_2, T) \sqrt{\vabsolut{\Lambda_k } } + \frac{\paren{1+p_2 }^2}{T}} \right] \mathcal P_k ,
	\end{equation}
for every $k\geq1$. In~\eqref{bounds}, $C$ is a positive constant only depending on $\vabsolut{\Lambda_1}$, $\beta$, $p_0$ and $\alpha$ (increasing with respect to $\alpha$), and $\mathcal P_k$, $\calH_1 (\rho, q, p_1,p_2)$ and $\calH_2 (\rho, q, p_1,p_2, T)$ are respectively given by
	\begin{equation}\label{Pk}
\mathcal P_k := \frac{1}{\displaystyle \prod_{ \{ n \ge 1:  \ 1\leq \left | k-n \right | <q \} } \vabsolut{ \Lambda_k - \Lambda_n } }, \quad \forall k \ge 1, \quad \hbox{ if } q \ge 2, 
	\end{equation}
$\calP_k := 1$, for every $k \ge 1$, if $q = 1$, 
	\begin{equation}\label{f13'}
	\left\{
	\begin{array}{l}
\displaystyle \calH_1 (\rho, q, p_1,p_2) = \paren{\frac{ 1 + \rho p_2^2 + q^2 }{\rho^2 p_1^4} }^{2q - 2} , \\
	\noalign{\smallskip}
\displaystyle \calH_1(\rho, q, p_1,p_2) = \paren{\frac{ 1 + \rho p_2^2 }{\rho^2 p_1^4} }^{2q - 2} , \quad \hbox{when $ \Lambda $ is real}. 
	\end{array}
	\right.
	\end{equation}
and
	\begin{equation}\label{f13}
	\left\{
	\begin{array}{l}
\displaystyle \calH_2 (\rho, q, p_1,p_2, T) = 1 + q + \sqrt{T} + \frac{1 + q}{\rho^2 p_1^2} + p_2, \\
	\noalign{\smallskip}
\displaystyle \calH_2 (\rho, q, p_1,p_2, T) = 1 + q + \sqrt{T} + \frac{1}{\rho^2 p_1^2} + p_2, \quad \hbox{when $ \Lambda $ is real}. 
	\end{array}
	\right.
	\end{equation}
\end{theorem}

%
%
%
%


\begin{remark}\label{r1.5}
It is clear that if $\Lambda =\left \{ \Lambda_k \right \}_{k \geq 1}$ is a sequence satisfying the assumptions in Theorem~\ref{Olive} for some parameters $\beta \in [0, \infty)$, $ \rho , p, \alpha \in (0, \infty) $ and $q \in \N$, then $\Lambda$ belongs to $ \mathcal{L}(\beta,\rho,q,p, p, p,  \alpha) $, and satisfies 
	\begin{equation*}\
\vabsolut{  \Lambda_{k}  -  \Lambda_n  } \ge \gamma   > 0, \quad \forall k , n \ge 1 : k \not= n,
	\end{equation*}
for a positive constant $\gamma$. As a consequence, we can apply Theorem~\ref{principal theorem} and deduce the existence of $\{q_k\} _{k\geq1} \subset L^2(0,T; \C)$, a biorthogonal family to $\{e_k \}_{k\geq1}$ in $L^2(0,T; \C)$, satisfying~\eqref{bounds}. Thanks to the previous gap condition, we get $\calP_k = 1$, if $q=1$, or
	$$
\calP_k \le \gamma^{ 2 - 2q}, \quad \forall k \ge 1, \quad \hbox{ if } q \ge 2. 
	$$
Combining this inequality and~\eqref{bounds} we deduce~\eqref{f0_1} for a positive constant $\calC_0$ independent of $k $ and $T$. Therefore, Theorem~\ref{principal theorem} is a generalization of Theorem~\ref{Olive} to the case of complex sequences that do not satisfy the gap condition~\eqref{f29}. 

On the other hand, we will see in Section~\ref{s2} that if sequence $\Lambda = \ens{\Lambda_k}_{k \ge 1} \subset \C$ satisfies~\eqref{f0_0}, or~\eqref{f28}--\eqref{f28b}, or~\eqref{f27}, then $\Lambda$ belongs to $ \mathcal{L}(0,\rho,q ,p_0, p_1, p_2,  \alpha)$, for appropriate $ \rho , p_0, p_1, p_2, \alpha \in (0, \infty) $ and $q \in \N$, and satisfies the gap condition~\eqref{f29}. Therefore, Theorem~\ref{principal theorem} generalizes the results on bounds of biorthogonal families to exponentials proved in~\cite{FR},~\cite{Miller},~\cite{cindea},~\cite{CMV1} and~\cite{CMV2}.        \fin
\end{remark}

The quantity $\calP_k$ in Theorem~\ref{principal theorem} provides a mesure of the condensation of the sequence $\Lambda$. When condition~\eqref{f29} holds, then, there exists a constant $\calC >0$ such that $\vabsolut{\calP_k } \le \calC$ for any positive integer $k$. But in general, $\calP_k$ could have any explosive behavior with respect to $k$ (see for instance Remark~\ref{r7}). 

In the next result we will prove that inequality~\eqref{bounds} is optimal with respect to $\calP_k$. This is our second main result:


\begin{theorem}\label{estimbas famillbio}
Let $\Lambda =\left \{ \Lambda_k \right \}_{k\geq}\subset\C$ be a complex sequence satisfying
	\begin{equation} \label{item6}
\vabsolut{ \Lambda_{k}-\Lambda_{n} } \leq \nu \vabsolut{ k^2 -n^2 }, \quad \forall k, n \geq1,
	\end{equation}
for $\nu >0$, and $\Lambda \in \mathcal{L} (\beta,\rho,q,p_0, p_1, p_2,  \alpha) $ with $\beta \in [0, \infty)$, $ \rho , p_0, p_1, p_2, \alpha \in (0, \infty) $ and $q \in \N$. Then, for any sequence $ \{q_k\} _{k\geq1} \subset L^2(0,T; \C) $  biorthogonal to $\{e_k \}_{k\geq1}$ in $L^2(0,T; \C)$ ($e_k$ is given in~\eqref{f0}), one has
	\begin{equation}\label{lowerbound}
\norma{q_{k} }_{L^2(0,T; \C)} \geq \max \ens { \frac{6}{\pi^2}  \calB_k \, e^{\frac{1}{T\nu }} ,\calE_k } \, \mathcal P_k ,  \quad \forall k \geq 3 ,
	\end{equation}
where $\mathcal P_k$ is given in~\eqref{Pk}, 
	\begin{equation}\label{f49}
\calB_k =
	\left\{
	\begin{array}{ll} 
\displaystyle \nu^{ k + q - 2 } \frac{ \paren{q - 1}! }{(q + 3)!} \paren{k + q }! \frac{ \paren{\nu T}^{k + 1} }{ \paren{1 + \nu T }^{2 k + q + 1} }  \frac{ \paren{2 k + q - 1} ! }{\paren{ 2 k + q + 1} !}  \sqrt{\delta \vabsolut{\Lambda_1 } + \frac{1}{2T} } , & \hbox{if }  1 \le k < q, \\
	\noalign{\smallskip}
\displaystyle \nu^{ 2 (q-1)} \frac{\left[ \paren{q - 1}! \right]^2}{(q + 3)!}  \frac{  \paren{k + q }! k  }{ (2 k - q)! }\frac{ \paren{\nu T}^{k + 1} }{ \paren{1 + \nu T }^{2 k + q + 1} }  \frac{ \paren{2 k + q - 1} ! }{\paren{ 2 k + q + 1} !}  \sqrt{\delta \vabsolut{\Lambda_1 } + \frac{1}{2T} } , & \hbox{if }   k \ge q,
	\end{array}
	\right.
	\end{equation}
	\begin{equation}\label{f49'}
\calE_k =
	\left\{
	\begin{array}{ll}
\displaystyle \frac{(k + q - 2) !}{T^{k + q - 2 }}  \paren{ \frac{2(k + q) - 3 }{2 T} + \delta \vabsolut{ \Lambda_1 } }^{1/2}  , & \hbox{if } 1 \le k < q, \\
	\noalign{\smallskip}
\displaystyle \frac{(2 q - 2) !}{T^{2 (q - 1) }}  \paren{ \frac{4q - 3 }{2 T} + \delta \vabsolut{ \Lambda_{ k + 1 - q } } }^{1/2} , & \hbox{if } k \ge q,
	\end{array}
	\right.
	\end{equation}
and $\delta $ is a positive constant only depending on $\beta$ ($ \delta = 1$ when $\beta = 0$). 
\end{theorem}


\begin{remark}
Theorem~\ref{estimbas famillbio} generalizes the results proved in~\cite{Gui},~\cite{CMV1} and~\cite{CMV2} to general complex sequences that might not satisfy the gap condition~\eqref{f29}. \fin
\end{remark}

As an application of Theorems~\ref{principal theorem} and~\ref{estimbas famillbio}, we will study the cost of fast controls $\calK(T)$ for system~\eqref{generic} in two situations in which condition~\eqref{f29} does not hold: 
\begin{enumerate}
\item First, we will analyze system~\eqref{generic} when the operator $L = L_2$ is given by~\eqref{L2} with $ Q \in L^2(0,\pi )$ a function such that  the sequence of eigenvalues of $L_2$ is given by
	$$
\Lambda^{(2)} = \ens{k^2 , k^2 + e^{- k ^{2 \gamma}}}_{k \ge 1}
	$$
and $\gamma \in (0, 1)$. In this example the minimal time associated to system~\eqref{generic} with $L = L_2$ is $T_0 (Q ) = 0$. Observe that the sequence $\Lambda^{(2)}$ does not satisfy~\eqref{f29} and, therefore, Theorem~\ref{Olive} cannot be applied. We will see that the sequence $\Lambda^{(2)}$ fulfills the assumptions in Theorems~\ref{principal theorem} and~\ref{estimbas famillbio} and, as a consequence, we will obtain new estimates (even with $T_0 (Q) = 0$) from above and from below for the control cost $\calK (T)$ associated to system~\eqref{generic} for $L = L_2$ (see Theorems~\ref{t8} and~\ref{t9}). These estimates show that the fast controls for system~\eqref{generic} with $L = L_2$ are more violent than those of the heat equation. This violent behavior comes from the condensation of the eigenvalues of the elliptic operator $L_2$.
\item We will also study system~\eqref{generic} with $L = L_3 $ (see~\eqref{f78}), and $\rho$, $\tau$ and $\xi$ positive constant satisfying~\eqref{f77} for an integer $ j \ge 1$. In this case we will check that system~\eqref{generic} is null controllable for any $T > 0$ and the corresponding control cost $\calK (T) $ satisfies~\eqref{f0_2} for a constant $C_0 = C_0 (\rho, \tau, \xi) > 0 $. With this example we generalize the null controllability result obtained in~\cite{GB-SN} for the linear version of~\eqref{PFSy}. 
\end{enumerate}

In a forthcoming paper (see~\cite{BGB}) we will carry out a more in-depth analysis of the cost of fast controls $\calK (T )$ of parabolic systems with a positive minimal time $T_0$ which comes from the condensation index associated to the sequence of eigenvalues of the generator of the corresponding $C_0$-semigroup.

\smallskip

The plan of the paper is the following: In Section~\ref{s2}, we will study some general properties of the sequences $\Lambda $ in the class $ \mathcal{L}(\beta,\rho,q,p_0, p_1, p_2,  \alpha)$, with $\beta \in [0, \infty)$, 
	$$ 
\rho , p_0, p_1, p_2, \alpha \in (0, \infty) 
	$$ 
and $q \in \N$. We will also provide in this section some examples of sequences $\Lambda $ in the literature that satisfy the conditions in Definition~\ref{d1}. Sections~\ref{s3} and~\ref{s4} will be respectively devoted to the proofs of the main results of this work, namely, Theorem~\ref{principal theorem} and Theorem~\ref{estimbas famillbio}. Finally, in Section~\ref{s5} we will apply the results on general bounds of biorthogonal families to complex sequences that do not satisfy the gap condition~\eqref{f29} to system~\eqref{generic} when $L = L_{2}$ (see~\eqref{L2}) and
	$$
\sigma (L_2) = \ens{k^2 , k^2 + e^{- k ^{2 \gamma}}}_{k \ge 1}
	$$
with $\gamma \in (0,1)$, and when $L = L_3 $ (see~\eqref{f78}) is such that $\rho, \tau,\xi \in (0, \infty)$ satisfy~\eqref{f77} for an integer $j \ge 1$. Some results presented in this fifth section have been announced in~\cite{GB-O}.

%
%
%


\section{Some general properties of sequences under the assumptions of Definition~\ref{d1}. Some examples}\label{s2}
We will devote this section to prove some general properties of sequences $\Lambda$ in the class of Definition~\ref{d1}, $ \mathcal{L}(\beta,\rho,q,p_0, p_1, p_2,  \alpha)$, with $\beta \in [0, \infty)$, $ \rho , p_0, p_1, p_2, \alpha \in (0, \infty)$ and $q \in \N$. These properties will be used in the proof of Theorems~\ref{principal theorem} and~\ref{estimbas famillbio}. We also complete this section with some examples of sequences $\Lambda$ that fulfill assumptions in Definition~\ref{d1}.

Let us first analyze the conditions which appear in Definition~\ref{d1} and condition~\eqref{item6} because in some particular cases they are redundant. To be precise, let us first check that the properties~\ref{item1}--\ref{item5} and~\eqref{item6} imply property~\ref{item7} for some $p_0$, $p_1$, $p_2$ and $\alpha$. One has:


\begin{proposition}\label{p0}
Let $\Lambda = \ens{\Lambda_k}_{k \ge 1} \subset \C$ be a sequence satisfying~\ref{item1},~\ref{item4}, \ref{item5} and~\eqref{item6} for some $\rho, \nu > 0$ and $q \ge1$. Then,~\ref{item7} holds, with
	\begin{equation*}
p_0 = \frac{1}{\sqrt{\nu}} , \quad p_1 =  \frac{1}{\sqrt{\nu}} , \quad p_2 = \frac{1}{\sqrt{\rho}}  \quad \hbox{and} \quad  \alpha = \max \ens{q - \sqrt{ \frac{ \vabsolut{\Lambda_1}}{\rho }} , \sqrt{ \frac{ \vabsolut{\Lambda_1}}{\rho } + 1} , \sqrt{\frac{ \vabsolut{\Lambda_1} }{\nu}} + 1}.
	\end{equation*}
\end{proposition}

\begin{proof}
Let us take $\Lambda = \ens{\Lambda_k}_{k \ge 1} \subset \C$, a sequence under the assumptions of the proposition, and let us prove that~\ref{item7} holds for appropriate parameters $p_0$, $p_1$, $p_2$ and $\alpha$. 

From~\ref{item5} and~\eqref{item6}, we have
	$$
\rho \paren{k^2 -n^2} \le \vabsolut{\Lambda_k} + \vabsolut{\Lambda_n}, \quad \forall k,n: k \ge n +q,  \quad \hbox{and} \quad \vabsolut{\Lambda_k} - \vabsolut{\Lambda_n} \le \nu \paren{k^2 -n^2}, \quad \forall k,n: k \ge n. 
	$$
In particular,
	\begin{equation}\label{f41}
	\left\{ 
	\begin{array}{l}
\vabsolut{\Lambda_k} \ge \rho \paren{k^2 - 1 } - \vabsolut{\Lambda_1}, \quad \forall k \ge q+1, \\
	\noalign{\smallskip}
\vabsolut{\Lambda_k} \le \nu \paren{k^2 - 1 } + \vabsolut{\Lambda_1}, \quad \forall k \ge 1. 
	\end{array}
	\right.
	\end{equation}

Let us consider $r \ge \vabsolut{\Lambda_{q + 1}}$. Taking into account the first item in Remark~\ref{r1}, if $\calN (r) = k$, then $k \ge q + 1$, $\vabsolut{\Lambda_k} \le r $ and $\vabsolut{\Lambda_{k +1}}  > r $. The first inequality in~\eqref{f41} gives $r \ge \rho \paren{k^2 - 1 } - \vabsolut{\Lambda_1}$, i.e.,
	$$
\calN (r) = k \le \sqrt{\frac 1\rho r + \frac{\vabsolut{\Lambda_1}}{\rho } + 1} \le \frac{1}{\sqrt{\rho}} \sqrt{r} + \sqrt{ \frac{ \vabsolut{\Lambda_1}}{\rho } + 1}, \quad \forall r \ge \vabsolut{\Lambda_{q + 1}} . 
	$$
On the other hand, the second inequality in~\eqref{f41} also provides $r < \nu \left[ \paren{k + 1}^2 -1 \right] + \vabsolut{ \Lambda_1} $ and
	$$
\calN (r) = k > -1 + \sqrt{\frac 1\nu r - \frac{\vabsolut{\Lambda_1} }{\nu} + 1} > \sqrt{\frac 1\nu r - \frac{\vabsolut{\Lambda_1} }{\nu}} -1 \ge  \frac{1}{\sqrt{\nu}} \sqrt{r} - \sqrt{\frac{\vabsolut{\Lambda_1} }{\nu}} -1. 
	$$
Observe that this inequality is also valid when $0 < r < \vabsolut{\Lambda_{q + 1}} $. In the previous reasoning we have used the inequalities
	\begin{equation}\label{f71}
	\left\{
	\begin{array}{l}
\displaystyle \sqrt{a + b } \le \sqrt{a} + \sqrt{ b }, \quad \forall a, b \in [0, \infty),\\
	\noalign{\smallskip}
\displaystyle \sqrt{a - b } \ge \sqrt{a} - \sqrt{ b }, \quad \forall a, b \in [0, \infty), \quad a \ge b. 
	\end{array}
	\right.
	\end{equation}

Let us now take $r$ such that $\vabsolut{\Lambda_{1}} \le r < \vabsolut{\Lambda_{q + 1}}$. In this case, 
	$$
\calN (r ) \le q \le \frac{1}{\sqrt{\rho}} \sqrt{r} + q - \frac{1}{\sqrt{\rho}} \sqrt{\vabsolut{\Lambda_{1}} } . 
	$$
Finally, when $r$ is such that $0 < r < \vabsolut{\Lambda_{1}} $, $\calN(r) = 0 \le  \sqrt{r} / \sqrt{\rho} $. We deduce then that $\Lambda $ satisfies~\ref{item7} with $p_0$, $p_1$, $p_2$ and $\alpha $ given in the statement. This proves the result.
\end{proof}


\begin{remark}\label{r2}
Property~\ref{item5} does not imply, in general,~\ref{item7}, even for increasing positive real sequences. Indeed, $\Lambda = \ens{k^3}_{k \ge 1}$ is an increasing positive real sequence that satisfies~\ref{item5}, with $\rho = 1$ and $q = 1$, and does not satisfy~\ref{item7}. 

Something similar can be said for property~\eqref{item6}: $\Lambda = \ens{k}_{k \ge 1}$ is an increasing positive real sequence that satisfies~\eqref{item6} with $\nu = 1$ and does not satisfy~\ref{item7}. 

On the other hand, sequences $\Lambda$ that satisfy~\ref{item1}--\ref{item5} and~\eqref{item6} for some $\beta \ge 0 $, $\rho, \nu > 0$ and $q \ge1$, also satisfy condition~\ref{item7} with parameters $p_0$, $p_1$, $p_2$ and $\alpha $ given in the statement of Proposition~\ref{p0}. In conclusion, $\Lambda \in \mathcal{L}(\beta,\rho,q,p_0, p_1, p_2,  \alpha)$. \fin
\end{remark}

As a consequence of the previous result, we also have a relation between the different parameters that appear in~\eqref{item6} and in Definition~\ref{d1}. To be precise, one has:


\begin{corollary}\label{c1}
Let $\Lambda = \ens{\Lambda_k}_{k \ge 1} \subset \C$ be a sequence satisfying~\ref{item1}, \ref{item4} and~\ref{item7} for some positive constants $ p_0 , p_1, p_2$. Then, 
\begin{enumerate}
\item If $\Lambda$ satisfies property~\ref{item5}, then 
	\begin{equation}\label{f23}
p_1 \le \frac{1}{\sqrt{\rho}}. 
	\end{equation}
\item If \eqref{item6} holds, then
	\begin{equation*}
\frac{1}{\sqrt{\nu}} \le p_2 . 
	\end{equation*}
\end{enumerate}
\end{corollary}


\begin{proof}
Let us consider a sequence $\Lambda = \ens{\Lambda_k}_{k \ge 1} \subset \C$ satisfying~\ref{item1}, \ref{item4} and~\ref{item7} for some parameters $ p_0 , p_1, p_2 , \alpha \in (0, \infty) $. In addition, let us assume that property~\ref{item5} holds. In particular, from~\ref{item7}, we get $p_1 \sqrt r - \alpha \le \calN(r)$, for any $r >0$. On the other hand, thanks to property~\ref{item5}, one also obtain (see the proof of Proposition~\ref{p0}), 
	$$
\calN (r) \le \sqrt{\frac 1\rho r + \frac{\vabsolut{\Lambda_1}}{\rho } + 1} \le \frac{1}{\sqrt{\rho}} \sqrt{r} + \sqrt{ \frac{ \vabsolut{\Lambda_1}}{\rho } + 1} , \quad \forall r \ge \vabsolut{\Lambda_{q + 1}}, 
	$$
that is to say, 
	$$
 p_1 \sqrt{r} - \alpha \le \frac{1}{\sqrt{\rho}} \sqrt{r} + \sqrt{ \frac{ \vabsolut{\Lambda_1}}{\rho } + 1} , \quad \forall r \ge \vabsolut{\Lambda_{q + 1}} . 
	$$
From this inequality, it is not difficult to deduce~\eqref{f23}. This proves item 1. 

The proof of the second item of the result follows the same ideas as above. We leave it to the reader. This ends the proof. 
\end{proof}


\begin{remark}\label{r3}
Observe that, if $\Lambda$ is a sequence under the conditions of Corollary~\ref{c1}, from inequality~\eqref{f23} we also deduce
	\begin{equation}\label{f23p}
0 < \rho \le \frac{1}{p_0^2} \, ,   \quad \rho p_1^2 \le 1  , \quad  \hbox{and} \quad  \rho p_1 \le \sqrt{\rho} \le \frac 1{p_0} .  
	\end{equation}
These estimates will be used later. \fin
\end{remark}

Let us now analyze the case of increasing positive real sequences $\Lambda = \ens{\Lambda_k}_{k \ge 1} \subset (0, \infty)$. This case is specially interesting because some assumptions in Definition~\ref{d1} are direct. For instance, $\Lambda$ satisfies~\ref{item1}--\ref{item4} for $\beta = 0$. In addition, if~\ref{item7} holds for some parameters $p_0$, $p_1$,  $p_2$ and $\alpha $ satisfying $p_1 = p_2 = p \ge p_0 >0$, some assumptions in Definition~\ref{d1} are redundant. To be precise, in this particular case,~\ref{item7} implies~\ref{item5} and the additional property~\eqref{item6}. One has:


\begin{proposition}\label{p1}
Let $\Lambda = \ens{\Lambda_k}_{k \ge 1}$ be a positive real sequence satisfying~\ref{item1}, \ref{item4} and~\ref{item7} for some $ p_0 , \alpha \in (0, \infty) $, with $p_1 = p_2 = p \ge p_0$. Then, $\Lambda \in \mathcal{L}( 0 ,\rho,q,p_0, p_1, p_2,  \alpha)$ and~\eqref{item6} holds, with
	\begin{equation}\label{f30}
q = {3 \alpha} , \quad \rho = \frac{1}{3 p^2} \quad \hbox{and} \quad \nu = \frac 13 \paren{\frac{2 + \alpha }{p}}^2.
	\end{equation}

\end{proposition}

\begin{proof}
Let us take $\Lambda = \ens{\Lambda_k}_{k \ge 1}$, a positive real sequence satisfying~\ref{item1}, \ref{item4} and~\ref{item7} for some $ p_0 , \alpha \in (0, \infty) $, with $p_1 = p_2 = p \ge p_0$. It is clear that $\Lambda$ satisfies~\ref{item2} and~\ref{item3} for $\beta = 0$. 

Let us see that $\Lambda$ also satisfies~\ref{item5} for appropriate positive constants $\rho$ and $q$. Indeed, using~\ref{item1} and \ref{item4} we infer that $\Lambda $ is an increasing positive real sequence. Thus, $\calN (\Lambda_k) = k$, for any $k \ge 1$, (see~\eqref{counting}) and, from~\ref{item7} ($p_1 = p_2 = p$), we deduce 
	\begin{equation}\label{f25}
k - \alpha \le p \sqrt{\Lambda_k} \le k + \alpha , \quad \forall k \ge 1. 
	\end{equation}

If $k , n \in \N$ are such that $k - n \ge 3 \alpha$, then, $k \ge \alpha$ and inequality~\eqref{f25} provides
	$$
\frac{p^2  \paren{\Lambda_k - \Lambda_n}}{k^2 - n^2 } \ge \frac{\paren{k - \alpha}^2  - \paren{n + \alpha}^2 }{k^2 - n^2 } = \frac{k - n - 2 \alpha}{k - n } = 1 - \frac{2 \alpha }{k - n } \ge \frac 13. 
	$$
Therefore,  sequence $\Lambda $ satisfies assumption~\ref{item5} for $q$ and $\rho$ as in the statement of the proposition. 

Let us now check property~\eqref{item6}. To this end, we will again use~\eqref{f25}. Without loss of generality, we can assume that $\alpha \ge 1$. Thus, if $\alpha < n < k$, one has
	$$
\frac{p^2  \paren{\Lambda_k - \Lambda_n}}{k^2 - n^2 } \le \frac{\paren{k + \alpha}^2  - \paren{n - \alpha}^2 }{k^2 - n^2 } = \frac{k - n + 2 \alpha}{k - n } = 1 + \frac{2 \alpha }{k - n } \le 1 + 2 \alpha  \le \frac 13 \paren{ 2 + \alpha}^2 . 
	$$
On the other hand, if $n \le \alpha < k$, i.e., if $n \le \lfloor \alpha \rfloor < \lfloor \alpha \rfloor +1 \le k$ ($ \lfloor \cdot \rfloor $ is the floor function: given $ x \in \R$, $ \lfloor x \rfloor $ is the greatest integer less than or equal to $x$), we also deduce 
	$$
\frac{p^2  \paren{\Lambda_k - \Lambda_n}}{k^2 - n^2 } \le \frac{\paren{k + \alpha}^2 }{k^2 -\lfloor \alpha \rfloor^2 } \le \frac{ \paren{ \lfloor \alpha \rfloor + \alpha + 1 }^2}{2 \lfloor \alpha \rfloor + 1 }  \le \frac 13 \paren{ 2 + \alpha}^2 . 
	$$
In the previous inequality we have used that $\alpha \ge 1$.

Finally, let us assume that $\alpha \ge 2$ and take $n < k \le \alpha$. We can write
	$$
\frac{p^2  \paren{\Lambda_k - \Lambda_n}}{k^2 - n^2 } \le \frac{ \paren{k + \alpha}^2}{2 k -1}  \le \frac 13 \paren{ 2 + \alpha}^2  .
	$$
Summarizing, property~\ref{item5} holds for $\nu $ given in~\eqref{f30}. This ends the proof.
\end{proof}




\begin{remark}\label{r4}
Let us consider $\Lambda = \ens{\Lambda_k}_{k \ge 1}$, an increasing positive sequence, satisfying property~\ref{item7} with $p_1 = p_2 = p > 0 $. In this case, this condition can be written under the equivalent form
	\begin{equation}\label{f42}
\Lambda_k = \frac{1}{p^2 } k^2 + O(k), \quad \forall k \ge 1 . 
	\end{equation}
Indeed, from~\ref{item7} with $p_1 = p_2 = p$, we infer~\eqref{f25} and
	$$
\vabsolut{p \sqrt{ \Lambda_k } -k} \le \alpha, \quad \forall k \ge 1,
	$$
i.e., $p \sqrt{ \Lambda_k } = k + O(1)$ for any $k \ge 1$. So,~\eqref{f42} holds. 

On the other hand, from~\eqref{f42} we deduce 
	\begin{equation}\label{f43}
\frac{1}{p^2} k^2 - \alpha_1 k \le \Lambda_k \le \frac{1}{p^2} k^2 + \alpha_1 k, \quad \forall k \ge 1,
	\end{equation}
with $\alpha_1 \ge 0$. Thus, given $ r >0$, if $\calN (r) = k $, then, (see Remark~\ref{r1}) we also have $\Lambda_k \le r$ and $\Lambda_{k + 1 } > r$. Using the previous inequalities, we obtain
	$$
\frac{1}{p^2} k^2 - \alpha_1 k - r \le 0 \quad \hbox{and} \quad \frac{1}{p^2} \paren{k + 1}^2 + \alpha_1 \paren{k + 1} - r > 0. 
	$$
In particular,  
	$$
	\left\{
	\begin{array}{l}
\calN (r) = \displaystyle k \le \frac 12 \paren{  p^2 \alpha_1 + p \sqrt{p^2 \alpha_1^2 + 4r }} \le p \sqrt r + p^2 \alpha_1, \\
	\noalign{\smallskip}
\displaystyle \calN (r) + 1 = k + 1 > \frac 12 \paren{  - p^2 \alpha_1 + p \sqrt{p^2 \alpha_1^2 + 4r }} > p \sqrt r - \frac 12 p^2 \alpha_1. 
	\end{array}
	\right.
	$$
Therefore,~\ref{item7} holds with $p_0 = p_1 = p_2 = p$ and
	$$
\alpha = \max \ens{p^2 \alpha_1 , \frac 12 p^2 \alpha_1 + 1}.
	$$

Observe that, in particular, if $\Lambda = \ens{\Lambda_k}_{k \ge 1}$ is an increasing real sequence such that~\eqref{f0_0} holds, then $\Lambda$ also satisfies~\eqref{f42} with $L = 1/p^2 $. As a consequence of Proposition~\ref{p1}, we can conclude that if $\Lambda = \ens{\Lambda_k}_{k \ge 1}$ is an increasing real sequence satisfying~\eqref{f0_0}, then  $\Lambda \in \mathcal{L}(\beta,\rho,q,p_0, p_1, p_2,  \alpha)$ and~\eqref{item6} holds for $\beta = 0$, $p_0 = p_1 = p_2 = p = 1/ \sqrt{C}$, $\alpha \in (0, \infty)$ and $q$, $\rho$ and $\nu $ as in~\eqref{f30}. Therefore, Theorems~\ref{principal theorem} and~\ref{estimbas famillbio} generalize the results on estimates of biorthogional families established in~\cite{FR} and~\cite{Gui}. \fin
\end{remark}


Let us continue showing some properties for sequences $\Lambda$ in the class $\Lambda \in \mathcal{L}(\beta,\rho,q,p_0, p_1, p_2,  \alpha)$. One has:
 

\begin{lemma} \label{lemmma1}
Let $\Lambda =\left \{ \Lambda_k \right \}_{k\geq}\subset\mathbb{C}$ be a sequence satisfying $\Lambda \in \mathcal{L}(\beta,\rho,q,p_0, p_1, p_2,  \alpha)$ with $\beta \in [0, \infty)$, $ \rho , p_0, p_1, p_2, \alpha \in (0, \infty) $ and $q \in \N$. Then, 
	\begin{equation} \label{convergence serie}
\sum_{k\geq1} \frac{1}{\vert \Lambda_k\vert} < \infty\quad \hbox{and} \quad  \vert \Lambda_k\vert \leq  \Re(\Lambda_k)+\beta \sqrt{\Re(\Lambda_k)} , \quad \forall k \ge 1. 
	\end{equation}
On the other hand, there exists a positive constant $C$, only depending on $\left | \Lambda_1 \right |$, $\beta$, $p_0$ and $\alpha$ (increasing with respect to $\alpha$), such that 
	\begin{equation} \label{condi*} 
\frac{1}{p_2} ( k - \alpha ) \leq \sqrt{\vert  \Lambda_k\vert} \leq \frac{1}{p_1} k + \frac{C ( 1 + q )}{\rho p_1^2} ,  \quad  \forall k\geq1.  
	\end{equation}
\end{lemma}


\begin{proof}
Let us take a sequence $\Lambda = \ens{\Lambda_k}_{k \ge 1}$ under assumptions of the proposition. From items~\ref{item4} and~\ref{item7} of Definition~\ref{d1}, we have that:
	\begin{equation*}
\sum_{k\geq1} \frac{1}{\vert \Lambda_k\vert}= \int_{\vert \Lambda_1 \vert}^{\infty} \frac{1}{r} \, d\mathcal{N}(r) = \int_{\vert \Lambda_1 \vert}^{\infty} \frac{1}{r^2} \mathcal{N}(r) \, dr \leq  \int_{\vert \Lambda_1 \vert}^{\infty} \frac{\alpha+p_2\sqrt{r}}{r^2} \, dr = \frac{\alpha}{\vert \Lambda_1 \vert}+ \frac{2p_2}{\sqrt{\vert\Lambda_1 \vert}}< \infty.
	\end{equation*}

On the other hand, using assumption \ref{item3}, we deduce that 
	$$
\vert \Lambda_k \vert^2 = \Re(\Lambda_k)^2+ \Im(\Lambda_k)^2\leq \Re(\Lambda_k)^2+\beta^2 \Re(\Lambda_k)\leq \left( \Re(\Lambda_k)+\beta \sqrt{\Re(\Lambda_k)} \right)^2.
	$$
Therefore, we have the proof of~\eqref{convergence serie}.

Let us now prove property~\eqref{condi*}.  Let us first assume that $\Lambda \in \mathcal{L}(\beta,\rho,q,p_0, p_1, p_2,  \alpha)$ is a positive real sequence ($\beta = 0$). We  have that $\mathcal N ( \Lambda_k ) = k$, for any $ k \ge 1$. In particular, taking $r = \Lambda_k $ in assumption~\ref{item7}, we deduce 
	$$
\frac{k - \alpha}{p_2} \le  \sqrt{\Lambda_k} \le \frac{ k + \alpha}{p_1} = \frac {k}{p_1} + \frac{\alpha}{ p_1 } \le 
 \frac {k}{p_1} + \frac{\alpha}{ p_0 }   \frac{1}{ \rho p_1^2 }  \rho p_1^2  \le   \frac {k}{p_1} + \frac{\alpha}{p_0 } \frac{1}{ \rho p_1^2 }, \quad \forall k \ge 1.
	$$
In the previous inequality we have used~\eqref{f23p}. This shows inequality~\eqref{condi*} in the real case.

Let us now assume that the sequence $\Lambda \in \mathcal{L}(\beta,\rho,q,p_0, p_1, p_2,  \alpha)$ is complex, i.e., $ \beta > 0$. As before, we are going to work with property~\ref{item7} with $r = \left| \Lambda_k \right|$ ($k \ge 1$). From Remark~\ref{r1},~\ref{item4} and~\ref{item7} (see Definition~\ref{d1}), we can write that, if $n = \calN \paren{\left| \Lambda_k \right|}$, then $k \le n $, $\left| \Lambda_k \right| = \left| \Lambda_n \right|$ and
	\begin{equation}\label{f1}
\frac{- \alpha + n}{p_2} = \frac{ - \alpha + \calN \paren{\vabsolut{\Lambda_k }}}{p_2} \le  \sqrt{\vabsolut{\Lambda_k }} \le \frac{\alpha + \calN \paren{\vabsolut{\Lambda_k }}}{p_1} =\frac{ \alpha + n}{p_1}, \quad \forall k \ge 1.
	\end{equation}
In particular, $k \le n$ and
	$$
- \alpha + k \le - \alpha + n \le p_2 \sqrt{\vabsolut{\Lambda_k }} , \quad \forall k \ge 1.
	$$
This proves the first inequality in~\eqref{condi*} in the complex case.

In order to show the second inequality in~\eqref{condi*}, let us estimate $n = \calN \paren{\left| \Lambda_k \right|}$. As $\left| \Lambda_k \right| = \left| \Lambda_n \right|$, using property~\ref{item3}, we infer
	$$
\left| \Re(\Lambda_n )^2-\Re(\Lambda_k)^2 \right| = \left| \Im(\Lambda_n)^2 - \Im(\Lambda_k)^2\right| \leq \beta^2(\Re(\Lambda_k )+\Re(\Lambda_n)),
	$$
that is to say,
	$$
\left |  \Re(\Lambda_k )-\Re(\Lambda_n)  \right | \leq \beta^2.
	$$
Again, assumption~\ref{item3} also provides the inequality
	$$
\left| \Lambda_k-\Lambda_n \right| \leq \left| \Re(\Lambda_k )-\Re(\Lambda_n)  \right| + \left| \Im(\Lambda_k ) -\Im(\Lambda_n)\right| \leq \beta^2 + 2\beta \sqrt{\left| \Lambda_k \right |}.
	$$
If $\vabsolut{k - n} \ge q$, combining the previous inequality and assumption~\ref{item5} we obtain
	$$
\rho \left| k -n \right| \paren{k + n} = \rho \left| k^2 -n^2 \right| \le \left| \Lambda_k-\Lambda_n \right| \leq \beta^2 + 2\beta \sqrt{\left| \Lambda_k \right |}.
	$$
Thus,
	\begin{equation*}
n - k = \left| k-n \right| \leq \max\left \{ q, \frac{\beta^2+2 \beta \sqrt{\left | \Lambda_k  \right |}}{\rho \paren{k + n} } \right \},
	\end{equation*}
i.e.,
	\begin{equation*}
n \leq k + \max\left \{ q, \frac{\beta^2+2 \beta \sqrt{\left | \Lambda_k  \right |}}{\rho \paren{k + \calN \paren{\left| \Lambda_k \right|}} } \right \},
	\end{equation*}
and, from~\eqref{f1}
	\begin{equation}\label{f2}
p_1 \sqrt{\vabsolut{\Lambda_k }} \le \alpha + k + \max\left \{ q, \frac{\beta^2+2 \beta \sqrt{\left | \Lambda_k  \right |}}{\rho \paren{k + \calN \paren{\left| \Lambda_k \right|}} } \right \}. 
	\end{equation}

If the maximum in~\eqref{f2} is $q$, in particular,
	$$
p_1 \sqrt{\vabsolut{\Lambda_k }} \le  k  + \alpha + q .
	$$
Taking into account inequalities~\eqref{f23} and~\eqref{f23p}, we also deduce
	$$
p_1 \sqrt{\vabsolut{\Lambda_k }} \le  k  + \frac{\alpha + q}{\rho p_1} \rho p_1 \le k + \frac{\paren{\alpha + q}/{p_0 }}{\rho p_1} .
	$$
Thus, we get the second inequality in~\eqref{condi*} for a positive constant $C$ only depending on $\alpha$ and $p_0$ and increasing with respect to $\alpha$. 

Let us now assume that the maximum in~\eqref{f2} is given by the second term. Using again~\ref{item7} and~\eqref{f23p}, for $k \ge \alpha$, we can write
	$$
p_1 \sqrt{\vabsolut{\Lambda_k }} \le k + \alpha + \frac{\beta^2 + 2 \beta  \sqrt{\left | \Lambda_k  \right |}}{\rho \paren{k + \calN \paren{\left| \Lambda_k \right|}} } \le  k + \alpha + \frac{\beta^2 + 2 \beta  \sqrt{\left | \Lambda_k  \right |}}{\rho \paren{k - \alpha + p_1 \sqrt{\left| \Lambda_k \right|}} } \le k + \frac{\alpha}{\rho p_1} \frac 1{p_0} + \frac{\beta^2 + 2 \beta  \sqrt{\left | \Lambda_k  \right |}}{\rho p_1  \sqrt{\left | \Lambda_k  \right |}}.
	$$
This inequality provides the second inequality in~\eqref{condi*} when $k \ge \alpha $ for a positive constant $C$ only depending on $\left | \Lambda_1 \right |$, $\beta$, $p_0$ and $\alpha$ (of course, increasing with respect to $\alpha$). 

Finally, let us consider the case $k < \alpha $. Thus, there exists a positive constant $C$ (only depending on $\alpha$ and increasing with respect to $\alpha$) such that
	$$
\sqrt{\vabsolut{\Lambda_k}} \le C \le \frac{k}{p_1} + \frac{C}{\rho p_1^2 } .
	$$
In the previous inequality we have used~\eqref{f23p}. 

Finally, it is not difficult to see that the constant $C$ appearing in the second inequality of~\eqref{condi*} is increasing with respect to the parameter $\alpha$. This ends the proof.
\end{proof}



\begin{remark}\label{r2.4'}
Analyzing the proof of Lemma~\ref{lemmma1} we deduce that, in fact, if the sequence $\Lambda$ is real and satisfies the assumptions of the result, then the second inequality in~\eqref{condi*} can be written as follows: there exists a positive constant $C$, only depending on $p_0$ and $\alpha$ (increasing with respect to $\alpha$) such that 
	\begin{equation} \label{condi**} 
\frac{1}{p_2} ( k - \alpha ) \leq \sqrt{\vert  \Lambda_k\vert} \leq \frac{1}{p_1} k + \frac{C }{\rho p_1^2} ,  \quad  \forall k\geq1.  
	\end{equation}
In particular the previous inequalities are independent of $q$. We will use this property for real sequences $\Lambda$ throughout the paper. \fin
\end{remark}



\begin{remark}\label{r2.4}
From the previous result we deduce that, if the sequence $\Lambda = \{ {\Lambda}_k\}_{k\geq1} $ is in the class $ \mathcal{L} (\beta,\rho,q,p_0, p_1, p_2,  \alpha)$, with $\beta \in [0, \infty)$, $ \rho , p_0, p_1, p_2, \alpha \in (0, \infty) $ and $q \in \N$, then one also has~\eqref{f18} for some $ \delta > 0$, only depending on $\beta$ ($ \delta = 1$ when $\beta = 0$).  

As said before, property~\eqref{f18} implies that the family of exponentials $\ens{e_k}_{k \ge 1}$, $e_k$ is given in~\eqref{f0},  is minimal in $L^2(0,T; \C)$ for any $T > 0$. Thus, there exists a biorthogonal family $\ens{\widetilde q_k}_{k \ge 1}$ to $\ens{e_k }_{k \ge 1}$ in $L^2(0,T; \C)$ (see for instance~\cite{S}, \cite{Re}, \cite{AKBGBdTJMPA}, Theorem~4.1 in~\cite{AKBGBdTJFA},...). \fin
\end{remark}


Let us complete this section providing some examples of sequences $\Lambda = \{ \Lambda_k \}_{k \ge 1 }$ such that $\Lambda \in \mathcal{L}(\beta,\rho,q,p_0, p_1, p_2,  \alpha) $ for some $\beta \in [0, \infty)$, $ \rho , p_0, p_1, p_2, \alpha \in (0, \infty) $ and $q \in \N$. In order to have a clearer exposition, we will present the results and we will include the corresponding proofs in an appendix, at the end of this paper.

Firstly, we will analyze the case of real sequences that fulfill the assumptions in~\cite{CMV1} and~\cite{CMV2}. In particular, this class of sequences also satisfies a gap condition and, therefore, the general assumptions in~\cite{BBGBO}. One has:


\begin{proposition}\label{p6}
Let $\Lambda = \ens{ \Lambda_k }_{ k \ge 1} \subset (0, \infty)$ be a real sequence satisfying~\eqref{f27} for two positive constants $\gamma_0$ and $\gamma_1 $. Then, $\Lambda = \ens{ \Lambda_k }_{ k \ge 1} \in \mathcal{L} (\beta ,\rho,q,p_0, p_1, p_2,  \alpha)  $ and~\eqref{item6} holds with $\beta = 0$,
	\begin{equation*}
	\left\{
	\begin{array}{l}
\displaystyle p_0 = p_1 = \frac{1}{\gamma_1} , \quad p_2 = \frac{1}{\gamma_0}, \quad \alpha = \max \ens{1 - \frac{\sqrt{\Lambda_1 }}{\gamma_0} , \frac{\sqrt{\Lambda_1 }}{\gamma_1}}, \\
	\noalign{\smallskip}
 \displaystyle q = 1,  \quad \rho = \min \ens{ \gamma_0^2, \frac 13 \gamma_0^2 + \frac 23 \gamma_0 \sqrt{\Lambda_1 } }  \quad \hbox{and} \quad \nu = \max \ens{ \gamma_1^2, \frac 13 \gamma_1^2 + \frac 23 \gamma_1 \sqrt{\Lambda_1 } }. 
	\end{array}
	\right.
	\end{equation*} 
In particular, the gap condition~\eqref{f29} holds. 
\end{proposition}

For the proof, see~\ref{a1}.


\begin{remark}
As said before, sequences $\Lambda = \ens{ \Lambda_k }_{ k \ge 1} \subset (0, \infty)$ under the assumptions of Proposition~\ref{p6} satisfy the general assumptions that assure the existence of a family $\{q_k\} _{k\geq1} \subset L^2(0,T)$ biorthogonal to $\{ e_k \}_{k\geq1}$ in $L^2(0,T)$ ($e_k$ is given in~\eqref{f0}) satisfying Theorem~\ref{principal theorem} and Theorem~\ref{estimbas famillbio}  with parameters given in the statement of the proposition. Observe, in particular, that $q=1$ and $\calP_k = 1$. Therefore, Theorems~\ref{principal theorem} and~\ref{estimbas famillbio} cover the results in~\cite{CMV1} and~\cite{CMV2}. \fin
\end{remark}

We continue our analysis of real sequences that fulfill general assumptions previously discussed in the literature. More specifically, we will analyze real sequences that fulfill the assumptions in~\cite{cindea}. One has:


\begin{proposition}\label{p3}
Let us consider two increasing sequences of positive real numbers 
	\begin{equation*}
\ens{ \lambda_k^{(1)}}_{ k \ge 1} \quad \hbox{and} \quad \ens{ \lambda_k^{(2)}}_{ k \ge 1}
	\end{equation*}
satisfying~\eqref{f28} and 
	\begin{equation}\label{f31}
\lambda^{(1)}_{k +1} - \lambda^{(1)}_k  \ge c_0 \quad \hbox{and} \quad  \lambda^{(2)}_{k + 1} - \lambda^{(2)}_k \ge c_0 , \quad \forall k \ge 1,
	\end{equation}
for some positive constants $\pi_1$, $\pi_2$, $c_0$, $c_1$ and $r$. Then, the sequence 
	$$
\Lambda  = \ens{ \lambda_k^{(1)}}_{ k \ge 1} \cup \ens{ \lambda_k^{(2)}}_{ k \ge 1}
	$$
can be rearranged as an increasing sequence $\Lambda  = \ens{\Lambda_k}_{k \ge 1} $ satisfying $\Lambda \in  \mathcal{L} (\beta,\rho,q,p_0, p_1, p_2,  \alpha) $,~\eqref{item6} and the gap condition~\eqref{f29}, with $\beta = 0$, $p_0 = \min \{ \pi_1, \pi_2 \}$,  $ p_1 = p_2 = p = \pi_1 + \pi_2$, 
	$$
\alpha = \max \ens{2 + \dfrac 12 c_1\paren{ \pi_1^2 + \pi_2^2}  , c_1\paren{ \pi_1^2 + \pi_2^2}  }
	$$
and $q$, $\rho$ and $\nu$ given in~\eqref{f30}. 
\end{proposition}

The proof of this result can be seen in~\ref{a2}. 


\begin{remark}\label{r6}
In~\cite{cindea} the authors consider families of positive real numbers $\ens{ \lambda_k^{(1)}}_{ k \ge 1}$ and $\ens{ \lambda_k^{(2)}}_{ k \ge 1}$ satisfying~\eqref{f28}, for some positive constants $\pi_1$, $\pi_2$, $c_1$ and $r$, and the strong gap condition~\eqref{f28b}, with $c_2$ a positive constant. In particular, these sequences fulfill assumptions~\eqref{f28} and~\eqref{f31} in Proposition~\ref{p3} and, therefore, the general hypotheses imposed to general complex sequences $\ens{\Lambda_k}_{k \ge 1}$ in~\cite{BBGBO} (see assumptions in Theorem~\ref{Olive}). Thus, the results on existence and sharp estimates of biorthogonal families established  in~\cite{cindea} can be deduced from the corresponding results proved in~\cite{BBGBO}. Of course, Theorem~\ref{principal theorem} generalizes the results in~\cite{cindea} and in~\cite{BBGBO} to complex sequences that do not satisfy the gap condition~\eqref{f29}. \fin
\end{remark}


As said before, in~\cite{AKBGBdTJFA} the authors prove the existence of a minimal time of controllability for some parabolic problems. This minimal time is related to the condensation index of the sequence of eigenvalues of the corresponding operator. In order to illustrate the existence of this minimal time, the authors consider the system
	\begin{equation} \label{f32}
	\left\{ 
	\begin{array}{ll}
\partial_{t} y - (D_1 \partial_{xx} + A_1) y = 0 & \text{in }   (0,T) \times (0, \pi)  , \\ 
 	\noalign{\smallskip}
y( \cdot, 0 ) = B v, \quad y( \cdot, \pi ) = 0 &\hbox{on  } (0,T),\\
  	\noalign{\smallskip}
y(0, \cdot) = y_0 & \hbox{in  } (0,\pi),\\
	\end{array}
	\right. 
	\end{equation}
where $B \in \R^2 $, $v \in L^2 (0, T)$ is the control,
	$$
D_1 := \hbox{diag} \, (1, d), \ d>0, \ d \not= 1, \quad \hbox{and} \quad A_1 := \left( 
	\begin{array}{cc}
0 & 1 \\ 
0 & 0
	\end{array}
\right).
	$$
The sequence of eigenvalues associated to the operator $L_1 = - (D_1 \partial_{xx} + A_1  )$, with domain $D (L_1) = H^2(0, \pi; \R^2 ) \cap H_0^1 (0, \pi; \R^2 )$, is given by $\Lambda = \ens{k^2}_{k \ge 1} \cup \ens{d k^2}_{k \ge 1} $. Remember that the condition $\sqrt{d } \not\in \Q$ is necessary for the approximate controllability of the system~\eqref{f32} at time $T>0$. On the other hand, under this assumption, there exists a minimal time $ T_0 = T_0 (d) \in [0, \infty ] $ such that the system is not null controllable at time $T$ when $T < T_0$ (see~\cite{AKBGBdTJFA} for the details). In our second example we will consider the sequence of eigenvalues associated to this system:

\begin{proposition}\label{p4}
Let us consider $d \in (0, \infty)$ such that $\sqrt{d} \not\in \Q$. Then, the sequence
	\begin{equation*}
\Lambda = \ens{k^2}_{k \ge 1} \cup  \ens{d k^2}_{k \ge 1} 
	\end{equation*}
can be rearranged as an increasing sequence $\Lambda = \ens{\Lambda_k}_{k \ge 1}$  satisfying $ \Lambda \in  \mathcal{L} (\beta,\rho,q,p_0, p_1, p_2,  \alpha) $ and condition~\eqref{item6} with $\beta = 0$, $p_0 = 1$,
	\begin{equation}\label{f33}
p_1 = p_2 = p = 1 + \frac{1 }{\sqrt{d}  }, \quad \alpha = 2, \quad q = 2,  \quad \rho = \frac 58 \frac{1}{ p^2} \quad \hbox{and} \quad \nu = \frac 83 \frac{1}{ p^2} .
	\end{equation}
\end{proposition}

The proof of Proposition~\ref{p4} can be found in~\ref{a3}. 

\smallskip


Let us now analyze a fourth example of sequence $\Lambda $ which satisfy~\eqref{item6} and the general conditions appearing in Definition~\ref{d1}. With this example we cover the kind of sequences associated to some parabolic problems studied in~\cite{O}:

\begin{proposition}\label{p5}
Let us consider two real positive sequences $\Lambda_1 = \ens{ \lambda_k^{(1)}}_{ k \ge 1}$ and $\Lambda_2 = \ens{ \lambda_k^{(2)}}_{ k \ge 1}$. Assume that $\Lambda_1$ satisfies $\Lambda_1 \in  \mathcal{L} ( 0 ,\rho_ 1, 1 , \pi_0, \pi_1, \pi_2,  \alpha_1 ) $, for $ \rho_1 , \pi_0, \pi_1, \pi_2 , \alpha_1 \in (0, \infty) $, and~\eqref{item6}, for $ \nu = \nu_1 \in (0, \infty)$. On the other hand, assume
	\begin{equation*}
\displaystyle \lambda_k^{(2)} = \lambda_k^{(1)} + \varepsilon_k , \quad \forall k \ge 1, \quad  \lambda^{(2)}_k\neq \lambda^{(2)}_n, \quad \forall k,n \ge 1, \hbox{ with } k \not= n , \quad \hbox{and} \quad \lambda^{(1)}_k\neq \lambda^{(2)}_n, \quad \forall k, n \ge 1,
	\end{equation*}
where $\ens{\varepsilon_k}_{k \ge 1}$ is a real bounded sequence. Let us take $\varepsilon_0 = \sup_{k \ge 1 } | \varepsilon_k |$. Then, the sequence 
	$$
\ens{ \lambda_k^{(1)}}_{ k \ge 1} \cup \ens{ \lambda_k^{(2)}}_{ k \ge 1}
	$$
can be rearranged as a positive increasing sequence $\Lambda  = \ens{\Lambda_k}_{k \ge 1} $ satisfying $\Lambda \in  \mathcal{L} (0 ,\rho, q, \pi_0, p_1, p_2,  \alpha) $ and~\eqref{item6}, with $\beta = 0$, $ p_1 = 2 \pi_1$,  $p_2 = 2 \pi_2$, $\alpha = \pi_2 \sqrt{\varepsilon_0 } + 2\alpha_1$, $q = 2$ and $\rho$ and $\nu$ positive constants only depending, resp., on $\rho_1$ and $\varepsilon_0$ and on $\rho_1$, $\nu_1$ and $\varepsilon_0$. 
\end{proposition}

For the proof, see~\ref{a4}.



\begin{remark}\label{r8}
Proposition~\ref{p5} covers the sequence of eigenvalues of operator $L $ in system~\eqref{generic} when $L = L_2 $ (see~\eqref{L2}). We will use this proposition in Section~\ref{s5}. \fin
\end{remark}



\begin{remark}\label{r5}
Under assumptions of Proposition~\ref{p5}, observe that the sequence $\Lambda$, in general, does not satisfy the gap condition~\eqref{f29}. In fact, it is easy to see that condition~\eqref{f29} holds if and only if 
	$$
\liminf \vabsolut{\varepsilon_k} > 0. 
	$$

On the other hand, analyzing the proof of Proposition~\ref{p5}, it is possible to provide some additional information about parameters $\rho$ and $\alpha $ in Proposition~\ref{p5} when the sequence $\ens{\varepsilon_k}_{k \ge 1}$ satisfies appropriate properties. Indeed, when the bounded sequence $\ens{\varepsilon_k}_{k \ge 1}$ is such that $\varepsilon_0 = \sup_{k \ge 1 } | \varepsilon_k |$ satisfies 
	$$
\vabsolut{\varepsilon_k } \le \varepsilon_0 \le \frac{ \rho_1 }{4} , \quad \forall k \ge 1,
	$$
then $k_0 = 1$ and the sequence $\Lambda $ can be explicitly defined by~\eqref{f40} for any $ k \ge 1 $ (see~\ref{a4}), that is to say,
	\begin{equation}\label{f40'}
\Lambda_k = \left\{
	\begin{array}{ll}
\displaystyle \min \ens{\lambda_{ \ell }^{(1)}, \lambda_{ \ell }^{(2)}}, & \hbox{if } k = 2 \ell - 1, \\
	\noalign{\medskip}
\displaystyle \max \ens{\lambda_{ \ell }^{(1)}, \lambda_{ \ell }^{(2)}}, & \hbox{if } k = 2 \ell,
	\end{array}
	\right.
	\end{equation}
for any $k \ge 1$. In addition, from the proof of Proposition~\ref{p5}, we can deduce
	$$
	\left\{
	\begin{array}{l}
\displaystyle \Lambda_{ k} - \Lambda_{ n} \ge \frac{ \rho_1 }{16}  \paren{k^2 -n^2}  , \quad \forall k,n \in \N: k \ge n + 2, \\
	\noalign{\smallskip}
\displaystyle \Lambda_{ k} - \Lambda_{ n} \le \frac{ \nu_1 + \varepsilon_0 }{2}  \paren{k^2 -n^2}  , \quad \forall k,n \in \N, 
	\end{array}
	\right.
	$$
i.e., we can take $\rho = \rho_1/16$ and $\nu = \paren{\nu_1 + \varepsilon_0}/2$ in Proposition~\ref{p5}. \fin
\end{remark}


As said in Remark~\ref{r9}, let us finalize this section with an academic example of a positive sequence $\Lambda$ in the class $\mathcal{L} (0 ,\rho, q, p_0, p_1, p_2,  \alpha)$ with a parameter $q$ which can be chosen as large as we want. With this example will see that the parameters $p_1$ and $p_2$ are increasing with respect to $q$. To this end, let us fix a positive integer $ m \ge 2 $. With this integer, we define   
	\begin{equation}\label{f72}
\Lambda = \left\{ k^2 + \frac{\ell - 1}m : k \ge 1, \quad 1 \le \ell \le m \right\}. 
	\end{equation}
It is clear that the set $\Lambda $ can be written as an increasing sequence $\Lambda = \ens{\Lambda_k}_{k \ge 1}$ that satisfies the gap condition~\eqref{f29}. Let us see that it also satisfies $\Lambda \in \mathcal{L} (0 ,\rho, q, p_0, p_1, p_2,  \alpha)$, for appropriate parameters $q \in \N$ and $\rho, p_0, p_1, p_2,  \alpha \in (0, \infty)$, and condition~\eqref{item6}, for $\nu > 0$. One has:


\begin{proposition}\label{p11}
Let us take a positive integer $ m \ge 2 $ and consider the sequence $\Lambda $ defined in~\eqref{f72}. Then, 
\begin{enumerate}
\item $\Lambda \in \mathcal{L} (0 ,\rho, q, p_0, p_1, p_2,  \alpha)$, with $q = m$, $p_0 = 2$, $ p_1 = p_2 = m $, $\alpha = m $ and 
	$$
\rho = \frac{2}{(2m-1)(2m+1)}. 
	$$
In fact, property~\ref{item5} does not hold if $q \le m-1$. 
\item The sequence $\Lambda $ satisfies~\eqref{item6} with 
	$$
\nu = \frac{4 m - 1}{m \left( 2m +1 \right) } . 
	$$
\end{enumerate}
\end{proposition}


\begin{proof}
If $m \ge 2$, it is clear that the sequence $\Lambda $, given in~\eqref{f72}, is an increasing sequence that satisfies items~\ref{item1}--\ref{item4}, with $\beta =0$. Let us check the other items in Definition~\ref{d1} and condition~\eqref{item6}:

\smallskip

\textbf{1.} Let us prove item~\ref{item7} for the sequence $\Lambda$. To be precise, let us see
	\begin{equation}\label{f73}
-m + m \sqrt r < \calN (r) \le m \sqrt r , \quad \forall r >0, 
	\end{equation}
where $\calN (r) $ is defined in~\eqref{counting}. First, if $r \in (0, 1)$, $\calN ( r ) = 0 $ and it is clear that~\eqref{f73} holds. Therefore, we will prove~\eqref{f73} when $ r \ge 1$. In this case, the function $\calN (r )$ is given by
	$$
\calN (r) = \sum_{\ell = 1}^m \#\left\{ k :  k^2 + \frac{\ell - 1}m \le r  \right\} = \sum_{\ell = 1}^m \left\lfloor \sqrt{r - \frac{\ell - 1}m} \right\rfloor \le  \sum_{\ell = 1}^m \left\lfloor \sqrt{r } \right\rfloor = m \sqrt{r }, \quad \forall r \ge 1. 
	$$

On the other hand, we can explicitly calculate $\calN (r) $: Given $r \ge 1$, there exists an integer $k \ge 1$ such that $r \in \left[ k^2 , (k + 1)^2 \right)$. In this case,
	\begin{equation}\label{f74}
\calN (r) = 
	\left\{
	\begin{array}{ll}
m \lfloor \sqrt r  \rfloor - m + \widetilde \ell , & \hbox{if} \quad r \in \left[ \displaystyle k^2 + \frac{\widetilde \ell - 1}m, k^2 + \frac{\widetilde \ell}m \right), \hbox{ with } \widetilde \ell \in \N: 1 \le \widetilde \ell \le m,     \\
	\noalign{\smallskip}
m \lfloor \sqrt r  \rfloor , & \hbox{if} \quad r \in \left[ \displaystyle k^2 + 1 , (k+1)^2 \right). 
	\end{array}
\right.
	\end{equation}
Indeed, if $r \in \left[ \displaystyle k^2 + \frac{\widetilde \ell - 1}m, k^2 + \frac{\widetilde \ell}m \right)$,  with $ \widetilde \ell \in \N: 1 \le \widetilde \ell \le m$, then, for any $\ell : 1 \le \ell \le \widetilde \ell $,
	$$
{r - \frac{\ell - 1}m} \in \left[ \displaystyle k^2 + \frac{\widetilde \ell - \ell }m, k^2 + \frac{\widetilde \ell - \ell + 1}m \right) \subset \left[ \displaystyle k^2 , (k + 1) ^2 \right) , 
	$$
and $\left\lfloor \sqrt{r - \frac{\ell - 1}m} \right\rfloor = k = \lfloor \sqrt r  \rfloor $. Also, if $\ell : \widetilde \ell + 1 \le \ell \le  m$, one has
	$$
{r - \frac{\ell - 1}m} \in \left[ \displaystyle k^2 - \frac{ \ell - \widetilde \ell }m, k^2 - \frac{ \ell - \widetilde \ell - 1}m \right) \subset \left[ \displaystyle (k - 1) ^2, k^2  \right) , 
	$$
and $\left\lfloor \sqrt{r - \frac{\ell - 1}m} \right\rfloor = k - 1 = \lfloor \sqrt r  \rfloor - 1  $. We deduce in this case
	$$
\calN (r) = \sum_{\ell = 1}^m \left\lfloor \sqrt{r - \frac{\ell - 1}m} \right\rfloor = m k - m + \widetilde \ell = m \lfloor \sqrt r  \rfloor - m + \widetilde \ell, 
	$$
and the first equality in~\eqref{f74}. 

Now, if  $ r \in \left[ \displaystyle k^2 + 1 , (k+1)^2 \right)$, we can apply the same reasoning as before and deduce
	$$
\left\lfloor \sqrt{r - \frac{\ell - 1}m} \right\rfloor = k  = \lfloor \sqrt r  \rfloor, \quad \forall \ell : 1 \le \ell \le m, 
	$$
and the second equality of~\eqref{f74}. 

Let us now prove the first inequality in~\eqref{f73} for $r \ge 1$. As before, $r \in \left[ k^2 , (k + 1)^2 \right)$, with $k \ge 1$ an integer. Thus, if 
	$$
r \in \left[ \displaystyle k^2 + \frac{\widetilde \ell - 1}m, k^2 + \frac{\widetilde \ell}m \right),
	$$
with $\widetilde \ell \in \N: 1 \le \widetilde \ell \le m$, then
	$$
	\left\{
\begin{array}{l}
\displaystyle  \calN (r) + m -m \sqrt r = mk + \widetilde \ell -m \sqrt r > m k - m \sqrt{ k^2 + \frac{\widetilde \ell}m } + \widetilde \ell  \\
	\noalign{\smallskip}
\displaystyle \phantom{\calN (r) + m -m \sqrt (r)} = \frac{ \left( m k + \widetilde \ell \right)^2 - \left( m^2 k^2 + m \widetilde \ell \right) }{ m k + \widetilde \ell + m \sqrt{ k^2 + \frac{\widetilde \ell}m } } = \frac{ \widetilde \ell^2 + m \widetilde \ell \left( 2 k -1 \right) }{ m k + \widetilde \ell + m \sqrt{ k^2 + \frac{\widetilde \ell}m } } > 0. 
	\end{array}
	\right.
	$$

Finally, if $ r \in \left[ \displaystyle k^2 + 1 , (k+1)^2 \right)$, we can write
	$$
\calN (r) =  m \lfloor \sqrt r  \rfloor > m \left( \sqrt r  - 1 \right). 
	$$
This proves~\eqref{f73} and property~\ref{item7} for the sequence $\Lambda$ with $p_0 = 2$, $p_1 = p_2 = m$ and $\alpha = m$.

\smallskip

\textbf{2.} Let us now see that property~\ref{item5} holds for $q = m$ (and an appropriate parameter $\rho > 0$) and is not valid if $q < m$. To this end, let us first provide the expression of the terms of the sequence $ \Lambda$. It is not difficult to see that, given an integer $k \ge 1$, this can be written as $k = m \widetilde k + \ell $, with $\widetilde k \ge 0$ and $\ell \in \N$ with $1 \le \ell \le m$. Thus,
	$$
\Lambda_k = \left( \widetilde k + 1 \right)^2 + \frac{\ell - 1 }{m }. 
	$$

Negative part: Fix $q \in \N$, with $1 \le q \le m - 1$, and take $n = m \widetilde k + 1$ and $k = m \widetilde k + q + 1 $, with $\widetilde k \ge 0$, an arbitrary integer. It is clear that $k - n = q \le m - 1$ and
	$$
\frac{\Lambda_{k } - \Lambda_{n}}{k^2 - n^2} = \frac{\left( \widetilde k + 1 \right)^2 + \frac{q }{m }  - \left( \widetilde k + 1 \right)^2 }{ \left( m \widetilde k + q + 1 \right)^2 - \left( m \widetilde k + 1 \right)^2 } = \frac{1}{ m  \left( 2 m \widetilde k + 2 + q \right) } \to 0, \quad \hbox{when } \widetilde k \to \infty. 
	$$
We deduce that property~\ref{item5} is not valid when $q \le m-1$. 

Positive part: Let us take $q = m$ and $n , k \ge 1$ with $k - n \ge q$. In this case,
	$$
k = m \widetilde k + \ell_2, \quad n = m \widetilde n + \ell_1, \quad \hbox{with }\widetilde n,  \ell_1, \widetilde k,  \ell_2 \in \Z, \  1 \le  \ell_1,  \ell_2 \le m \hbox{ and } \widetilde k , \widetilde n \ge 0. 
	$$
Observe that, thanks to the inequality $k - n \ge q = m$, we can conclude $ \widetilde k - \widetilde n  \ge 1$. So,
	$$
	\begin{array}{l}
\displaystyle \frac{\Lambda_{k} - \Lambda_{n}}{k^2 - n^2} = \frac{\left( \widetilde k + 1 \right)^2 + \frac{\ell_2 - 1 }{m }  - \left( \widetilde n + 1 \right)^2 - \frac{\ell_1 - 1 }{m }}{ \left( m \widetilde k + \ell_2 \right)^2 - \left( m \widetilde n + \ell_1 \right)^2 } = \frac{ \left( \widetilde k - \widetilde n \right) \left( \widetilde k + \widetilde n +  1 \right) + \left( \widetilde k - \widetilde n \right) + \frac{\ell_2 - \ell_1 }{m }} { \left( m \widetilde k + \ell_2 \right)^2 - \left( m \widetilde n + \ell_1 \right)^2 } \\ 
	\noalign{\smallskip}
\displaystyle \phantom{\frac{\Lambda_{k} - \Lambda_{n}}{k^2 - n^2} } \ge \frac{ \left( \widetilde k - \widetilde n \right) \left( \widetilde k + \widetilde n + 1 \right) + \widetilde k - \widetilde n - 1 + \frac{ 1 }{m }} { \left( m \widetilde k + m \right)^2 - \left( m \widetilde n + 1 \right)^2 } > \frac{\widetilde k - \widetilde n }{ m \left( \widetilde k - \widetilde n \right) + m - 1 } \cdot \frac{\widetilde k + \widetilde n + 1}{ m \left( \widetilde k + \widetilde n + 1 \right) + 1 } \\ 
	\noalign{\smallskip}
\displaystyle \phantom{\frac{\Lambda_{k} - \Lambda_{n}}{k^2 - n^2} } \ge \frac{1}{2m -1 } \cdot \frac{2}{2m +1 }. 
	\end{array}
	$$
This shows property~\ref{item5} for the sequence $\Lambda $ with $q = m $ and $ \rho $ given in the statement.

\smallskip

\textbf{3.} In order to finish the proof of this result, let us show property~\eqref{item6}. Again, let us take $k, n \in \N$ with $ k > n $. As before, 
	$$
k = m \widetilde k + \ell_2, \quad n = m \widetilde n + \ell_1, \quad \hbox{with }\widetilde n,  \ell_1, \widetilde k,  \ell_2 \in \Z, \  1 \le  \ell_1,  \ell_2 \le m \hbox{ and } \widetilde k , \widetilde n \ge 0 \hbox{ with } \widetilde k \ge \widetilde n \ge 0. 
	$$

Let us first analyze the case $ \widetilde k = \widetilde n = \widehat k \ge 0 $ and, of course, $1 \le \ell_1 < \ell_2 \le m $. We deduce,
	$$
\frac{\Lambda_{k } - \Lambda_{n}}{k^2 - n^2} = \frac{\left( \widehat k + 1 \right)^2 + \frac{\ell_2 - 1 }{m }  - \left( \widehat k + 1 \right)^2 -  \frac{\ell_1 - 1 }{m }}{ \left( m \widehat k + \ell_2 \right)^2 - \left( m \widehat k + \ell_1 \right)^2 } = \frac 1m \, \frac{1}{2 m \widehat k + \ell_2 + \ell_1 }  \le \frac 1{3m}. 
	$$

Now, if $  \widetilde k > \widetilde n $ and $1 \le \ell_1, \ell_2 \le m$, one gets
	$$
	\begin{array}{l}
\displaystyle \frac{\Lambda_{k} - \Lambda_{n}}{k^2 - n^2} 
= 
\frac{ \left( \widetilde k - \widetilde n \right) \left( \widetilde k + \widetilde n +  2 \right) + \frac{\ell_2 - \ell_1 }{m }} { \left( m \widetilde k + \ell_2 \right)^2 - \left( m \widetilde n + \ell_1 \right)^2 } 
\le 
\frac{ \left( \widetilde k - \widetilde n \right) \left( \widetilde k + \widetilde n + 2 \right) + 1 - \frac{ 1 }{m }} { \left( m \widetilde k + 1 \right)^2 - \left( m \widetilde n + m \right)^2 } \\ 
	\noalign{\smallskip}
\displaystyle \phantom{\frac{\Lambda_{k} - \Lambda_{n}}{k^2 - n^2} }  
= 
\frac{\widetilde k - \widetilde n }{ m \left( \widetilde k - \widetilde n \right) + 1 - m  } \cdot \frac{\widetilde k + \widetilde n + 2}{ m \left( \widetilde k + \widetilde n + 1 \right) + 1 } + \frac{ 1 - \frac{ 1 }{m }} { \left( m \widetilde k + 1 \right)^2 - \left( m \widetilde n + m \right)^2 } \\ 
	\noalign{\smallskip}
\displaystyle \phantom{\frac{\Lambda_{k} - \Lambda_{n}}{k^2 - n^2} } 
\le 
\frac{3}{2m + 1} + \frac{m - 1}{m \left( 2m +1 \right) } 
=
\frac{4 m - 1}{m \left( 2m +1 \right) }. 
	\end{array}
	$$

Taking into account that $m \ge 2 $, we can infer that 
	$$
\frac 1{3m} \le \frac{4 m - 1}{m \left( 2m +1 \right) }
	$$
and, therefore, the sequence $\Lambda$ fulfills inequality~\eqref{item6} with $\nu$ given in the statement. This ends the proof of the proposition. 
\end{proof}



\begin{remark}
It is interesting to point out that, thanks to Proposition~\ref{p1}, once property~\ref{item7} is proved for the sequence $\Lambda $ with $ p_0 = 2$, $p_1 = p_2 = m$ and $\alpha = m$, we can conclude that $\Lambda \in \mathcal{L}(0 ,\widetilde \rho, \widetilde q, 2 , m , m ,   m )$ and~\eqref{item6} holds, with (see~\eqref{f30})
	\begin{equation*}
\widetilde q = {3 m} , \quad \widetilde \rho = \frac{1}{3 m ^2} \quad \hbox{and} \quad \widetilde \nu = \frac 13 \paren{\frac{2 + m }{m}}^2.
	\end{equation*}
The parameters provided by Proposition~\ref{p11} are better than the previous values. Indeed, taking into account that $m \ge 2$, it is clear that $q = m < \widetilde q = {3 m} $,
	$$
\rho = \frac{2}{(2m-1)(2m+1)} > \widetilde \rho = \frac{1}{3 m ^2} \quad \hbox{and} \quad  \nu = \frac{4 m - 1}{m \left( 2m +1 \right) }   < \widetilde \nu = \frac 13 \paren{\frac{2 + m }{m}}^2. \eqno \square
	$$
\end{remark}


\begin{remark} 
We can apply~Theorems~\ref{principal theorem} and~\ref{estimbas famillbio} to the sequence $\Lambda $ given by~\eqref{f72} and conclude the existence of a sequence $ \{q_k\} _{k\geq1} \subset L^2(0,T)$, biorthogonal to $\{e_k \}_{k\geq1}$ in $L^2(0,T; \C)$ ($e_k$ is given in~\eqref{f0}), which satisfies~\eqref{bounds} and~\eqref{lowerbound}. If we make use of Proposition~\ref{p11}, these two inequalities can be written under the form 
	\begin{equation}\label{f75}
\calA_k^{(1)} ( m ) \, \mathcal P_k \le \|q_k \|_{L^2(0,T)} \leq \calA_k^{(2)} ( m ) \, \mathcal P_k, \quad \forall k \ge 3,
	\end{equation}
where $\calA_k^{(1)} ( m ) := \calE_k $ (see~\eqref{f49'}) and 
	\begin{equation*}
\calA_k^{(2)} ( m ) := \calH_1(\rho, q, p_1,p_2)  \exp \left[C \paren{ 1 + \calH_2 (\rho, q, p_1,p_2, T) \sqrt{\vabsolut{\Lambda_k } } + \frac{\paren{1+p_2 }^2}{T}} \right], 
	\end{equation*}
(see~\eqref{f13'} and~\eqref{f13} in the real case) with $\rho$, $q$, $p_1$, $p_2$ and $\nu$ given in Proposition~\ref{p11} (recall that the parameter $m$ is the maximal cardinal of the condensation groupings of the sequence $\Lambda$, that is to say the maximal number of elements in $\Lambda$ that do not satisfy~\ref{item5} and could condense). 

Observe that, taking into account Remark~\ref{r4}, the elements of the sequence $\Lambda$ satisfy
	\begin{equation*}
\Lambda_k = \frac{1}{m^2 } k^2 + O(k), \quad \forall k \ge 1
	\end{equation*}
and, therefore, one has
	$$
\lim_{m \to \infty} \calS (m) = \infty \quad \hbox{where} \quad \calS (m ) = \sum_{k \ge 1} \frac{1}{ \Lambda_k }, \quad \forall m \ge 2. 
	$$
In some sense, the family of exponentials $\ens{e_k}_{k \ge 1}$ ($e_k$ given in~\eqref{f0}) ``loses'' its property of minimality in $L^2(0,T)$ when $m$ tends to infinity. Thus, it is natural that the constants $\calA_k^{(1)} ( m )$ and $\calA_k^{(2)} ( m ) $ in~\eqref{f75} satisfy
	\begin{equation}\label{f76}
\lim_{m \to \infty} \calA_k^{(1)} ( m ) = \lim_{m \to \infty} \calA_k^{(2)} ( m ) = \infty, \quad \forall k \ge 1. 
	\end{equation}

Let us see that~\eqref{f76} holds. To this end, we will analyze the asymptotic behavior of $\calA_k^{(1)} ( m )$ and $\calA_k^{(2)} ( m )$ when $m \to \infty$. In what follows, we will provide an explicit expression of these constants when $3 \le k \le m$. 

\begin{enumerate}
\item Let us first analyze $\calA_k^{(1)} ( m )$. From the expression of $\calE_k $ for $3 \le k \le m$ (see~\eqref{f49'}) and Proposition~\ref{p11}, we can write
	$$
\calA_k^{(1)} ( m ) := \calE_k = \frac{(m + k - 2) !}{T^{m + k - 2 }}  \paren{ \frac{2(m + k ) - 3 }{2 T} + 1 }^{1/2}. 
	$$
Observe that Stirling's formula implies the existence of a positive constant $c_0 > 0$ such that
	$$
n! \ge c_0 \sqrt{2 \pi n} \paren{\frac ne}^n, \quad \forall n \in \N. 
	$$
In particular, for a new positive constant $c$ (independent of $m$), we deduce
	$$
\calA_k^{(1)} ( m ) \ge c \, \sqrt{(m + k -2 ) \left[ \frac{2(m + k ) - 3 }{2 T} + 1 \right]} \paren{\frac{ m + k - 2 }{e T }}^{m + k - 2}, 
	$$
which is valid for any $ m \ge 2 $ and any $k : 3 \le k \le m$. One has the first equality in~\eqref{f76}.

\item We continue with the analysis of $\calA_k^{(2)} ( m )$. Let us start with $\calH_1(\rho,q, p_1, p_2) $ (see~\eqref{f13'} in the real case). From Proposition~\ref{p11}, this constant only depends on $m$ and has the expression:
	$$
\calH_1(\rho,q, p_1, p_2) \equiv \calH_1 (m) = \left[ \frac{ \left( 6 m^2 - 1 \right) \left( 4m^2 - 1 \right) }{4 m^4} \right]^{2(m-1)}, \quad \forall m \ge 2. 
	$$
It is not difficult to see that
	$$
\lim_{m \to \infty }\frac{ \calH_1 (m) }{6^{2( m- 1) }} = 1,
	$$
and, then
	$$
c_1 6^{2( m- 1) } \le \calH_1 (m) \le c_2 6^{2( m- 1) }, \quad \forall m \ge 2. 
	$$
for two positive constants $c_1$ and $c_2$, independent of $m$. 

On the other hand, from the expression of $\calH_2 (\rho,q, p_1, p_2) $ (see~\eqref{f13} in the real case), we can write
	$$
\calH_2 (\rho, q, p_1,p_2, T) \equiv \calH_2 (m , T) = 4 m^2 + 2m -1 + \frac{1}{4 m^2 } + \sqrt{T}. 
	$$

Observe that in our case $\alpha = m$. We can conclude that $\calA_k^{(2)} ( m )$ is given by
	$$
\calA_k^{(2)} ( m ) := \calH_1(m)  \exp \left[C(m) \paren{ 1 +\calH_2 (m, T) \sqrt{\Lambda_k } + \frac{\paren{ 1+ m }^2}{T}} \right], \quad  m \ge 2,
	$$
with $C(m)$ a positive constant only depending on $m$ and increasing with respect to $m$ (see Theorem~\ref{principal theorem} with $\alpha = m$). Clearly, $\calA_k^{(1)} ( m )$  has an exponential behavior with respect to $m$ and we can write
	$$
\calA_k^{(2)} ( m ) \ge \exp \left[C \paren{ 1 + m^2 \left(  \sqrt{\Lambda_k } + \frac 1T \right) + \sqrt{T \Lambda_k } } \right], \quad \forall m \ge 2, \quad \forall k : 3 \le k \le m. 
	$$

We can conclude that $\calA_k^{(2)} ( m )$ has an exponential behavior with respect to $m$ and saisfies the second equality in~\eqref{f76}. \fin
\end{enumerate}
\end{remark}



\section{Proof of the first main result}
\label{s3} 
This section is devoted to prove Theorem \ref{principal theorem}. The main idea we will use is the Fourier transform together with the Paley-Wiener Theorem.
We need to introduce the following definition and  recall the Paley-Wiener Theorem.

\begin{definition} An entire function $f$ is said to be of exponential type $A$ if the inequality 
	\begin{equation*}
\left | f(z) \right | \leq B e^{A \left | z \right |}
	\end{equation*}
holds for some positive constants $A$ and $B$ and all values of $z\in \mathbb{C}.$
\end{definition}

Let us now present the Paley-Wiener Theorem:

\begin{theorem} \label{paleywiner}
Let $f$ be an entire function of exponential type $A > 0$ such that 
	\begin{equation*}
\left \| f \right \|_{L^2(\mathbb{R})} := \paren{\int_{- \infty}^\infty | f (x) |^2 \, dx}^{1/2} < \infty.
	\end{equation*}
Then, there exists a function $\phi \in L^2 (-A,A; \C)$ such that 
	\begin{equation*}
f(z) = \frac{1}{\sqrt{2 \pi}} \int_{-A}^{A}\phi (t) e^{izt}\,dt.
	\end{equation*}
Moreover, the Plancherel theorem gives
	$$
\| \phi \|_{L^2(- A, A; \C)} = \left \| f \right \|_{L^2(\mathbb{R})}.
	$$
\end{theorem}

For the proof of Theorem~\ref{paleywiner} we refer to \cite[Theorem 18. p.~101]{young1981introduction}.

\begin{remark}
In what follows, $C$ will denote a positive constant independent of $T$, $k\in \mathbb{N}$, $\rho$, $q$, $p_1$ and $p_2$, which may change from one line to another ($C$ may depend on $\vabsolut{\Lambda_1}$, $\beta$, $p_0$ and $\alpha$, and is increasing with respect to $\alpha$). In this work, the dependence of the constants with respect to the parameters $\rho$, $q$, $p_1$ and $p_2$ (see assumptions~\ref{item5} and~\ref{item7}) will be explicitly given. \fin
\end{remark}


Let us begin with a result of existence of entire functions satisfying appropriate properties. Our first main result will be a consequence of this theorem. One has:

\begin{theorem} \label{lemme technique}
Let $\Lambda =\left \{ \Lambda_k \right \}_{k\geq}\subset\mathbb{C}$ be a sequence satisfying $\Lambda \in \mathcal{L}(\beta,\rho,q,p_0, p_1, p_2,  \alpha)$ with $\beta \in [0, \infty)$, $ \rho , p_0, p_1, p_2, \alpha \in (0, \infty) $ and $q \in \N$. Then, for all $T>0$, there exists a sequence of entire functions $\{G_k\}_{k\geq1}$, with the following properties:
\begin{enumerate}
\item For any $k \ge 1$ and $\varepsilon > 0$, there exists a positive constant $C'_{T,k, \varepsilon }$ such that 
	\begin{equation}\label{f3}
\left | e^{-iz\frac{T}{2}} G_k(z)\right | \leq C'_{T,k, \varepsilon } e^{\left( \frac{T}{2} + \varepsilon \right) \left | z \right | }, \quad \forall z \in \C;
	\end{equation}
\item $\displaystyle G_k(i \overline{\Lambda}_n) = \frac{1}{\sqrt{2 \pi}} \delta_{kn}$, for all $k,n\geq1$;
\item $G_k$   belongs to $L^2(\mathbb{R})$, for any $k\geq1$, and there exists a positive constants $C>0$, only depending on $\vabsolut{\Lambda_1}$, $\beta$, $p_0$ and $\alpha$ (increasing with respect to $\alpha$), such that
	\begin{equation} \label{estimatt} 
\left \| G_k \right \|_{L^2(\R)}  \leq \calH_1(\rho, q, p_1,p_2)  \exp \left[C \paren{
1 + \calH_2 (\rho, q, p_1,p_2, T) \sqrt{\vabsolut{\Lambda_k } } + \frac{\paren{1 + p_2 }^2}{T}} \right] \mathcal P_k, 
	\end{equation}
for any $ k \geq 1$, where $\mathcal P_k$, $\calH_1(\rho, q, p_1,p_2)$ and $\calH_2 (\rho, q, p_1,p_2, T)$ are respectively given in~\eqref{Pk}, \eqref{f13'} and~\eqref{f13}.
\end{enumerate}
\end{theorem}

Theorem~\ref{principal theorem} is a direct consequence of Theorem~\ref{lemme technique}. Therefore, before providing the proof of the technical result established in Theorem~\ref{lemme technique}, we will complete the proof of Theorem~\ref{principal theorem}. 


\begin{proof}[Proof of Theorem~\ref{principal theorem}]
Let us consider a sequence $\Lambda =\left \{ \Lambda_k \right \}_{k\geq}\subset\mathbb{C}$ such that $\Lambda \in \mathcal{L}(\beta,\rho,q,p_0, p_1, p_2,  \alpha)$ with $\beta \in [0, \infty)$, $ \rho , p_0, p_1, p_2, \alpha \in (0, \infty) $ and $q \in \N$. On the other hand, let us fix $T > 0$. With the previous data, let us consider the function
	$$
F_k(z):= G_k(z){e^{-iz \frac{T}{2}}}, \quad z \in \C, \quad k \ge 1,
	$$
where $\{ G_k \}_{k \ge 1}$ is the sequence provided by Theorem~\ref{lemme technique}. Let us see some properties of the function $F_k$. First, $F_k$ is, for any $k \ge 1$, an entire function over $\C$. In fact, $F_k \in L^2 (\R)$ with 
	$$ 
\left \| F_k \right \|_{L^2(\R)} = \left \| G_k \right \|_{L^2(\R)}, \quad \forall k \ge 1.
	$$

Secondly, for any $\varepsilon > 0$ and $k \ge 1$, $F_k$ is an entire function of exponential type $T/2 + \varepsilon$ (see~\eqref{f3}). So, we can apply Payley-Wiener Theorem (see Theorem~\ref{paleywiner}) and deduce that there exists 
	$$
\psi_k \in L^2(- {T}/{2} - \varepsilon, {T}/{2}+ \varepsilon; \C)
	$$ 
such that
	$$
F_k(z) ={e^{-iz \frac{T}{2}}} G_k(z) = \frac{1}{\sqrt{2 \pi}} \int_{-\infty }^{\infty}  \psi_k(t)\, e^{iz t}\,dt, \quad \forall z \in \C, \quad \forall k \ge 1.
	$$

Observe that the support of the function $\psi_k$ is contained in $[- {T}/{2} - \varepsilon, {T}/{2}+ \varepsilon]$, for any $k \ge 1$ and for any $\varepsilon > 0$. We conclude that, in fact,  $\psi_k \in L^2(- {T}/{2} , {T}/{2}; \C)$ and
	\begin{equation}\label{f15}
F_k(z) ={e^{-iz \frac{T}{2}}} G_k(z) = \frac{1}{\sqrt{2 \pi}} \int_{-\frac{T}{2}}^{\frac{T}{2}}  \psi_k(t)\, e^{iz t}\,dt, \quad \forall z \in \C, \quad \forall k \ge 1.
	\end{equation}

Let us now consider the function
	\begin{equation}\label{f16}
q_k (t) := \psi_k \paren{t - \frac T2 } , \quad t \in [0,T], \quad k \ge 1.
	\end{equation}

It is clear that $q_k$ is well defined and $q_k \in L^2(0,T; \C)$ for any $k \ge 1$. The objective now is to prove that the sequence $\{ q_k\}_{k \ge 1} \subset L^2(0,T; \C)$ satisfies Theorem~\ref{principal theorem}. Let us first see that $\{ q_k\}_{k \ge 1}$ is biorthogonal to $\{e^{-\Lambda_k t}\}_{k\geq1}$ in $L^2(0,T; \C)$. Indeed, for any $k,n \ge 1$ and thanks to~\eqref{f15} and item 2 in Theorem~\ref{lemme technique}, we can write,
	\begin{equation*}
	\begin{split}
\int_0^T q_k(t) e^{- \overline \Lambda_n t} \, dt & =  \int_0^T \psi_k \paren{t - \frac T2 } e^{- \overline \Lambda_n t} \, dt = e^{- \overline \Lambda_n \frac T2} \int_{-\frac T2}^{\frac T2} \psi_k \paren{ t } e^{ - \overline \Lambda_n t} \, dt \\
	\noalign{\smallskip}
& = e^{- \overline \Lambda_n \frac T2} \sqrt{2 \pi} e^{ \overline \Lambda_n \frac T2} G( i \overline \Lambda_n) = \delta_{kn}.
	\end{split}
	\end{equation*}

As said before, $q_k \in L^2(0,T; \C)$. Let us now estimate $\| q_k \|_{L^2(0,T)}$. To this aim, we will use Plancherel Theorem and estimate~\eqref{estimatt}. From the expression of $q_k$ (see~\eqref{f16}), one has 
	\begin{equation*}
\| q_k \|_{L^2(0,T; \C )} = \|\psi_k\|_{L^2(-\frac{T}{2},\frac{T}{2}; \C)} = \| F_k\|_{L^2(\R)} = \| G_k\|_{L^2(\R)}.
	\end{equation*}
Combining the previous inequality and inequality~\eqref{estimatt} we deduce~\eqref{bounds}.
This completes the proof of Theorem~\ref{principal theorem}.
\end{proof}

Once Theorem~\ref{principal theorem} is proved, our next objective will be to show Theorem~\ref{lemme technique}. The proof of this result is very technical. In order to make it clearer, we will divide it in two subsections: 

\begin{enumerate}
\item In the first subsection (see Subsection~\ref{ss31})  we will introduce an entire function $f_k (z)$ ($ k \ge 1 $) with simple zeros at $\Lambda_n$ with $n \ge 1$ and $ n \not= k$. To this end, we will use the natural infinite product that satisfies the condition $f_k (\Lambda_n) = 0$ for any $ n \not= k$. We will show some properties of this function that, in particular, will imply item 2 in Theorem~\ref{lemme technique}. 
\item In the second subsection (see Subsection~\ref{ss32}) we will introduce a ``mollifier'' function that we will use in the definition of the entire function $G_k$ ($ k \ge 1$) in Theorem~\ref{lemme technique}. We will prove some properties of this function (which, in particular, will provide the property of item~3 in Theorem~\ref{lemme technique}) and we will complete the proof of Theorem~\ref{lemme technique}.
\end{enumerate}

\begin{remark}
Let us remark that, if the sequence $\Lambda = \{ {\Lambda}_k\}_{k\geq1}$ satisfies $\Lambda \in \mathcal{L} (\beta,\rho,q,p_0, p_1, p_2,  \alpha)$, with $\beta \in [0, \infty)$, $ \rho , p_0 , p_1, p_2,  \alpha \in (0, \infty) $ and $q \in \N$ (see Definition~\ref{d1}), then, the sequence $\overline{\Lambda} := \{ \overline{\Lambda}_k\}_{k\geq1}$ also belongs to $\mathcal{L} (\beta,\rho,q,p_0, p_1, p_2,  \alpha)$. As a consequence, we will prove Theorem~\ref{lemme technique} for the sequence $ \overline{\Lambda} $ instead of  $\Lambda $. \fin
\end{remark}


\subsection{An infinite product} \label{ss31}
In this section we will consider again a sequence $\Lambda = \{ {\Lambda}_k\}_{k\geq1}$ satisfying $\Lambda \in \mathcal{L} (\beta,\rho,q,p_0, p_1, p_2,  \alpha)$, for $\beta \in [0, \infty)$, $ \rho , p_0, p_1, p_2, \alpha \in (0, \infty) $ and $q \in \N$. Thus, for each $k\geq1$ and $z\in \mathbb{C}$, we define
	\begin{equation} \label{f14}
f_k(z):= \underset{n\neq k}{\prod_{n\geq1}} \left( 1-\frac{z}{\Lambda_n} \right), \quad z \in \C. 
	\end{equation}
The objective of this section is to prove some interesting properties satisfied by the function~$f_k$. 

First, observe that, by property~\eqref{convergence serie}, the previous product is uniformly convergent on compact sets of $\mathbb{C}$. Therefore,  $f_k$ is, for  any $k\geq1$, an entire function over $\C$ (see for instance \cite[p.~457]{E. Hille}). Moreover, 
	$$
f_k(\Lambda_n)=0, \quad \forall n\neq k.
	$$
In fact, the zeros of $f_k$ are exactly the elements of the sequence $\{\Lambda_n\}_{n\geq1,\, n\neq k}$ and they are zeros of multiplicity 1.

We have the following property of function $f_k$:


\begin{lemma} [\cite{BBGBO}]\label{lemmma3}
Let $\Lambda =\left \{ \Lambda_k \right \}_{k\geq}\subset\mathbb{C}$ be a sequence satisfying $\Lambda \in \mathcal{L} (\beta,\rho,q,p_0, p_1, p_2,  \alpha)$ with $\beta \in [0, \infty)$, $ \rho , p_0, p_1, p_2, \alpha \in (0, \infty) $ and $q \in \N$. Then, for  every $z\in \mathbb{C}$ and $k\geq 1$, we have
	\begin{equation} \label{p2}
\log \vert f_k(z) \vert \leq (p_2 \pi + 1) \sqrt{\vert z\vert}+C, 
	\end{equation}
where $p_2$ is given in assumption~\ref{item7} and $C$ is a positive constant only depending on $\alpha$ and $| \Lambda_1 |$ and increasing with respect to $\alpha$.
\end{lemma}


\begin{proof}
The proof of this result can be found in~\cite{BBGBO}. For completeness, we provide the proof here.

Given $z\in \mathbb{C}$, one has 
	\begin{equation*}
\log\left| f_k (z) \right| \leq \underset{n\neq k}{\sum_{n\geq1}} \log\left ( 1+ \frac{\left | z \right |}{\left | \Lambda_k\right |} \right ) \leq \int_{\left | \Lambda_1 \right |}^{ \infty} \log\left ( 1+ \frac{\left | z \right |}{\left | t\right |} \right ) \, d\calN (t).
	\end{equation*}
By assumption \ref{item7}, we get 
	$$
\lim_{t\rightarrow \infty}  \frac{\mathcal{N}(t)}{t}=0,
	$$
and an integration by parts provides 
	$$
\int_{\left | \Lambda_1  \right |}^{ \infty} \log \left ( 1+\frac{\left | z \right |}{t} \right ) \, d\calN (t) = \int_{\left | \Lambda_1  \right |}^{\infty} \frac{\left | z \right |}{t (\left | z \right |+t )}\, \calN (t) \, dt .
	$$
The change of variables $t= \left | z \right | s$ leads to  
	$$
\int_{\left| \Lambda_1  \right |}^{\infty}  \frac{\left | z \right |}{t(\left | z \right |+t)} \, \calN (t) \, dt = \int_{\frac{\left | \Lambda_1  \right |}{\left | z \right |}}^{\infty} \frac{ \mathcal{N}({\left | z \right | s})}{ s(s+1)} \, ds. 
$$
Using again assumption~\ref{item7}, we can conclude
	\begin{equation*}
	\begin{array}{l}
\displaystyle \int_{\frac{\left | \Lambda_1  \right |}{\left | z \right |}}^{\infty} \frac{ \mathcal{N}({\left | z \right | s})}{ s(s+1)} \, ds \leq p_2 \sqrt{\left | z \right |}  \int_{\frac{\left | \Lambda_1  \right |}{\left | z \right |}}^{\infty} \frac{1}{\sqrt{s} \paren{s + 1}}  \,ds + \alpha \int_{\frac{\left | \Lambda_1  \right |}{\left | z \right |}}^{ \infty} \frac{1}{s(s+1)} \,ds \\
	\noalign{\smallskip}
\displaystyle \phantom{\int_{\frac{\left | \Lambda_1  \right |}{\left | z \right |}}^{\infty} \frac{ \mathcal{N}({\left | z \right | s})}{ s(s+1)} \, ds} \leq p_2 \pi \sqrt{\left | z \right |}+ \alpha \log\left ( 1+ \frac{\left | z \right |}{\vabsolut{\Lambda_1} } \right ).
	\end{array}
	\end{equation*}
Finally, it is easy to check that there exists a positive constant $C$ (only depending on $\alpha$ and $| \Lambda_1 |$ and increasing with respect to $\alpha$) such that
	$$
\alpha \log\left ( 1+ \frac{\left | z \right |}{\vabsolut{\Lambda_1}} \right ) - \sqrt{\left | z \right |} \le C, \quad \forall z \in \C. 
	$$
Thus, we can conclude that inequality~\eqref{p2} holds. This finishes the proof. 
\end{proof}


Remember that our objective is to construct a sequence $\ens{G_k}_{k \ge 1}$ of entire functions over $\C$ satisfying items~1--3 in Theorem~\ref{lemme technique}. This construction will use the function $f_k (z) $ and an estimate from below of the non-zero quantity $\vabsolut{f_k (\Lambda_k)}$. This is one of the key points of this work and is established in the next


\begin{lemma}\label{l1}
Let $\Lambda =\left \{ \Lambda_k \right \}_{k\geq}\subset\mathbb{C}$ be a sequence satisfying $\Lambda \in \mathcal{L} (\beta,\rho,q,p_0, p_1, p_2,  \alpha)$ with $\beta \in [0, \infty)$, $ \rho , p_0, p_1, p_2, \alpha \in (0, \infty) $ and $q \in \N$. Then,  
	\begin{equation} \label{f19}
\vert f_k(\Lambda_k)\vert \geq \calH_1 (\rho, q, p_1,p_2)^{-1} e^{- C \, \calH_3 (\rho, q, p_1,p_2) \sqrt{\vert\Lambda_k\vert}} \, \calP_k^{- 1}, \quad \forall k\geq1,
	\end{equation}
where $C  $ is a positive constant, only depending on $\vabsolut{\Lambda_1}$, $\beta$, $p_0$ and $\alpha$ (increasing with respect to $\alpha$), $\calH_1 (\rho, q, p_1,p_2)$, $f_k$ and $\calP_k$ are respectively given in~\eqref{f13'}, \eqref{f14} and~\eqref{Pk}, and $ \calH_3$ is defined by
	\begin{equation*}
	\left\{
	\begin{array}{l}
\displaystyle \calH_3 (\rho, q, p_1,p_2) = 1 + q + \frac{1 + q}{\rho^2 p_1^2} + p_2, \\
	\noalign{\smallskip}
\displaystyle \calH_3 (\rho, q, p_1,p_2) = 1 + q + \frac{1}{\rho^2 p_1^2} + p_2, \quad \hbox{when $ \Lambda $ is real}. 
	\end{array}
	\right.
	\end{equation*}

\end{lemma}


\begin{proof}
As said before, if $\Lambda =\left \{ \Lambda_k \right \}_{k\geq}\subset\mathbb{C}$ satisfies $\Lambda \in \mathcal{L} (\beta,\rho,q,p_0, p_1, p_2,  \alpha) $ for constants $\beta \in [0, \infty)$, $ \rho , p_0, p_1, p_2, \alpha \in (0, \infty) $ and $q \in \N$, then $f_k$ (see~\eqref{f14}) is an entire function over $\C$ with simple zeros at the points $\ens{\Lambda_n}_{n \ge 1, n\not=k}$. Moreover, from assumption~\ref{item1}, we have
	$$
\vert f_k(\Lambda_k)\vert  = \underset{n\neq k}{\prod_{n\geq1}}  \left| \frac{\Lambda_n-\Lambda_k}{\Lambda_n} \right| \not= 0.
	$$

In order to obtain lower estimates of $| f_k (\Lambda_k) | $ let us decompose the set $\{ n \ge 1 : n \not= k \}$ into the following sets: 
	\begin{equation*}
	\left\{
	\begin{array}{l}
S_1 (k):=  \ens{ n \ge 1: 1 \le | n - k | < q }, \\
	\noalign{\smallskip}
S_2 (k) := \ens{ n \ge 1: | n - k| \ge q, \ \vert \Lambda_n\vert \leq 2 \vert \Lambda_k \vert },\\
	\noalign{\smallskip}
S_3 (k):= \ens{ n \ge 1: | n - k| \ge q, \ \vert \Lambda_n \vert>2\vert \Lambda_k \vert }.
	\end{array}
	\right.
	\end{equation*}
Then, 
	\begin{equation}\label{f22}
\left | f_k(\Lambda_k) \right | = \prod_{n\in S_1(k)} \left | 1-\frac{\Lambda_k}{\Lambda_n} \right | \prod_{n\in S_2(k)} \left | 1-\frac{\Lambda_k}{\Lambda_n} \right | \prod_{n\in S_3(k)} \left | 1-\frac{\Lambda_k}{\Lambda_n} \right | := \prod_{i=1}^3 P^{(k)}_i .
	\end{equation}
%


Let us estimate each term in~\eqref{f22} and, to this aim, let us take $n \in S_1( k )$. In particular, $n < k + q$ and, from~\ref{item4} and~\eqref{condi*} (or~\eqref{condi**} in the real case), we deduce
	\begin{equation*}
	\begin{split}
\vabsolut{\Lambda_n} & \le \vabsolut{\Lambda_{k + q}} \le \frac{2}{p_1^2} \vabsolut{ k + q }^2 + \frac{2C(1 + q)^2}{ \rho^2 p_1^4} \le  \frac{2}{p_1^2} \vabsolut{p_2 \sqrt{\vabsolut{\Lambda_k}} + \alpha + q}^2 + \frac{2C(1 + q)^2}{ \rho^2 p_1^4} \\
	\noalign{\smallskip}
& \le 4 \frac{p_2^2}{p_1^2} \vabsolut{\Lambda_k}+ \frac{2C(1 + q)^2}{ p_1^2} + \frac{2C(1 + q)^2}{ \rho^2 p_1^4} := 4 \frac{p_2^2}{p_1^2}\vabsolut{\Lambda_k} + A, \quad \forall n \in S_1(k),
	\end{split}
	\end{equation*}
(or
	$$
\vabsolut{\Lambda_n}  \le 4 \frac{p_2^2}{p_1^2} \vabsolut{\Lambda_k}+ \frac{2C}{ p_1^2} + \frac{2C}{ \rho^2 p_1^4} := 4 \frac{p_2^2}{p_1^2}\vabsolut{\Lambda_k} + A, \quad \forall n \in S_1(k),
	$$
when $\Lambda $ is a real sequence). In the previous inequalities,  $C$ is a positive constant independent of $\rho$, $q$, $p_1$ and $p_2$. 

One has
	$$
\log \left( 4 \frac{p_2^2}{p_1^2} x + A \right) = \log ( x)+  \log \left(4 \frac{p_2^2}{p_1^2}+  \frac{A}{x} \right) \leq \sqrt{x} +   \log  \left( 4\frac{p_2^2}{p_1^2} +  \frac{A}{\vabsolut{\Lambda_1}} \right)  , \quad \forall x \ge \vabsolut{\Lambda_1}. 
	$$
On the other hand, thanks to~\eqref{f23} and~\eqref{f23p}, we also deduce
	$$
	\begin{array}{l}
\displaystyle 4 \frac{p_2^2}{p_1^2} + \frac{A}{\vabsolut{\Lambda_1}} 
=  
4 \frac{p_2^2}{p_1^2} + \frac{C(1 + q^2 )}{\vabsolut{\Lambda_1}} \paren{\frac{1}{ p_1^2} + \frac{1}{ \rho^2 p_1^4}} 
=  
\frac{4 \vabsolut{\Lambda_1} \rho^2 p_1^2 p_2^2 + C(1 + q^2 ) \left( \rho^2 p_1^2 + 1 \right)}{\vabsolut{\Lambda_1} \rho^2 p_1^4} \\
	\noalign{\smallskip}
\displaystyle \phantom{4 \frac{p_2^2}{p_1^2} + \frac{A}{\vabsolut{\Lambda_1}} }
\le  
\frac{C \left( 1 + \rho p_2^2 + q^2 \right) }{\rho^2 p_1^4},
	\end{array}
	$$
(or 
	$$
4 \frac{p_2^2}{p_1^2} + \frac{A}{\vabsolut{\Lambda_1}} \le   \frac{C \left( 1 + \rho p_2^2 \right) }{\rho^2 p_1^4}
	$$
when $\Lambda $ is real). 

Thus,
	\begin{equation}\label{f20}
	\left\{
	\begin{array}{l}
\displaystyle P^{(k)}_1 = \prod_{n\in S_1(k)} \frac{\left |\Lambda_n -\Lambda_k \right |}{\left | \Lambda_n  \right |} 
\geq 
\prod_{n\in S_1(k)} \frac{\left |\Lambda_n -\Lambda_k \right |}{  \frac{ 4p_2^2}{p_1^2} \left | \Lambda_k  \right |  + A } 
= 
\frac{1}{\paren{ \frac{ 4p_2^2}{p_1^2} \left | \Lambda_k  \right | + A }^{2q - 2}} \prod_{n\in S_1(k)} \left |\Lambda_n -\Lambda_k \right | \\
	\noalign{\smallskip}
\displaystyle \phantom{ P^{(k)}_1} \geq \frac{1}{\left( \frac{ 4p_2^2}{p_1^2}  + \frac{A}{\vabsolut{\Lambda_1}}  \right)^{2q - 2}} e^{- (2q -2) \sqrt{\left | \Lambda_k  \right |}} \, \calP_k^{- 1} 
\ge 
C \paren{\frac{\rho^2 p_1^4}{ 1 + \rho p_2^2 + q^2 }}^{2q - 2} e^{- (2q -2) \sqrt{\left | \Lambda_k  \right |}} \, \calP_k^{- 1}, 
	\end{array}
	\right.
	\end{equation}
where $\calP_k$ is given in~\eqref{Pk} and $C$ is a new a positive constant independent of $\rho$, $q$, $p_1$ and $p_2$. In the real case, we deduce the following inequality for $P^{(k)}_1$:
	\begin{equation}\label{f20'}
\displaystyle P^{(k)}_1 \ge C \paren{\frac{\rho^2 p_1^4}{ 1 + \rho p_2^2 }}^{2q - 2} e^{- (2q -2) \sqrt{\left | \Lambda_k  \right |}} \, \calP_k^{- 1}. 
	\end{equation}

Let us now estimate the product $P_2^{(k)}$. At this point we will use the gap condition assumed in hypothesis~\ref{item5} when $ \vabsolut{n - k} \ge q$. We will follow some ideas from~\cite{BBGBO}. Using again Lemma~\ref{lemmma1} (or inequality~\eqref{condi**} in the real case), we deduce 
	\begin{equation}\label{f24}
\frac{k + n}{\sqrt{\vabsolut{\Lambda_k}}} \ge \frac{k}{\frac {k}{p_1} + \frac{C(1 + q) }{\rho p_1^2}} = \frac{\rho p_1^2 k}{\rho p_1 k + C (1 + q) } \ge \frac{\rho p_1^2 }{\rho p_1  + C (1 + q) }:= B, \quad \forall n,k \geq 1.
	\end{equation}
(or
	\begin{equation}\label{f24'}
\frac{k + n}{\sqrt{\vabsolut{\Lambda_k}}} \ge \frac{k}{\frac {k}{p_1} + \frac{C }{\rho p_1^2}} = \frac{\rho p_1^2 k}{\rho p_1 k + C } \ge \frac{\rho p_1^2 }{\rho p_1  + C }:= B, \quad \forall n,k \geq 1.
	\end{equation}
when $\Lambda$ is real). 
Then, if $n \in S_2(k)$,
	\begin{equation*}
P_2^{(k)}= \prod_{n\in S_2 (k)}\left | \frac{\Lambda_n - \Lambda_k}{\Lambda_n} \right |\geq \prod_{n\in S_2 (k)} \frac{ \rho}{2} \frac{\vert k-n\vert \vert k+n\vert}{\vert \Lambda_k\vert} \ge \prod_{n \in S_2 (k)} \frac{B \rho}{2} \frac{\vert k-n\vert}{\sqrt{\vert \Lambda_k\vert}} ,
	\end{equation*}
where $B$ is given in~\eqref{f24} (or~\eqref{f24'} in the real case).

Let $r_k:=\#\ens{n \in S_2(k) : n \le k - q }$ and $s_k:=\# \ens{  n \in S_2(k): n \ge q + k }$. Then, from the previous estimate, one has  
	$$ 
P_2^{(k)}\geq r_k! \left(\frac{B \rho }{2 \sqrt{\vert \Lambda_k\vert}} \right)^{r_k} s_k! \left( \frac{B \rho }{2\,\sqrt{\vert \Lambda_k\vert}} \right)^{s_k} := P_{2,1}^{(k)} P_{2,2}^{(k)}.
	$$

Observe that Stirling's formula implies the existence of a positive constant $c_0 > 0$ such that
	$$
n! \ge c_0 \sqrt{2 \pi n} \paren{\frac ne}^n, \quad \forall n \in \N. 
	$$
On the other hand, for $c_1 = e^{-1}$ one has
	$$
x \log x \ge -c_1, \quad \forall x \in (0, \infty).
	$$
Thus,
	\begin{equation*}
	\left\{
	\begin{array}{l}
\displaystyle P_{2,1}^{(k)} = r_k! \left(\frac{B \rho }{2\sqrt{ \vert \Lambda_k\vert}} \right)^{r_k} \geq c_0 \left(\frac{ B \rho r_k}{ \,2\, e \sqrt{ \vert \Lambda_k\vert}} \right)^{r_k} =c_0 \,\exp\left( \frac{2 e\sqrt{\Lambda_k} } {B \rho } \frac{B \rho  r_k}{2 e\sqrt{\vert\Lambda_k}\vert } \log\left(\frac{ B \rho r_k }{2 e\sqrt{\vert\Lambda_k}\vert} \right) \right) \\
	\noalign{\smallskip}
\displaystyle \phantom{P_{2,1}^{(k)}} \geq c_0 \exp\left(\frac{-2 c_1 e}{ B \rho } \sqrt{ \vert \Lambda_k\vert} \right).
	\end{array}
	\right.
	\end{equation*}

Taking into account the expression of $B$ (see~\eqref{f24}, resp.,~\eqref{f24'} in the real case) and inequalities~\eqref{f23} and~\eqref{f23p}, we can conclude the existence of positive constant $C_1$ and $C_2$ (independent of $\rho$, $q$, $p_1$, $p_2$ and $T$) such that
	$$
-\frac{C_1(1 + q )}{\rho^2 p_1^2} \le  \frac{-1}{B \rho} \le -\frac{C_2(1 + q) }{\rho^2 p_1^2},  \quad \hbox{(resp.,} \quad -\frac{C_1}{\rho^2 p_1^2} \le  \frac{-1}{B \rho} \le -\frac{C_2 }{\rho^2 p_1^2} , \hbox{ in the real case)}. 
	$$
As a consequence,
	$$
P_{2,1}^{(k)} \ge C \exp\left(\frac{-C (1 + q) }{  \rho^2 p_1^2 } \sqrt{ \vert \Lambda_k\vert} \right),  \quad \hbox{(resp.,} \quad P_{2,1}^{(k)} \ge C \exp\left(\frac{-C }{  \rho^2 p_1^2 } \sqrt{ \vert \Lambda_k\vert} \right) , \hbox{ in the real case)}. 
	$$

A similar reasoning can be applied to $P_{2,2}^{(k)}$. Therefore, we have proved:  
	\begin{equation}\label{f21}
P_2^{(k)} \ge C \exp\left(\frac{-C (1 + q)}{  \rho^2 p_1^2 } \sqrt{ \vert \Lambda_k\vert} \right),  \quad \hbox{(resp.,} \quad P_{2}^{(k)} \ge C \exp\left(\frac{-C }{  \rho^2 p_1^2 } \sqrt{ \vert \Lambda_k\vert} \right) , \hbox{ in the real case)} ,
	\end{equation}
for any $ k\geq1$. Again, $ C $ is a positive constant independent of $\rho$, $q$, $p_1$, $p_2$ and $T$.

In order to finish, let us analyze the third product in~\eqref{f22}. To this aim, we will use the inequalities
	\begin{equation*}
	\left\{
	\begin{array}{l}
\displaystyle \log(1-x)\geq-2x,  \quad \forall x\in\left[0, \frac{1}{2}\right], \\
	\noalign{\smallskip}
\displaystyle \frac{\vert\Lambda_k \vert }{\vert \Lambda_n \vert}< \frac{1}{2}, \quad \forall n \in S_3(k).
	\end{array}
	\right.
	\end{equation*}
From these inequalities and~\ref{item7}, we can write $\calN (r) \le \alpha + p_2 \sqrt{r}$, for any $r > 0$, and
	\begin{equation*} 
	\begin{split}
\log P_3^{(k)} &\geq \underset{n\in S_3 (k)} {\sum} \log\left(1-\frac{\vert \Lambda_k \vert}{\vert \Lambda_n\vert}\right) \geq -2 \vert \Lambda_k\vert  \underset{n\in S_3 (k)} {\sum}  \frac{1}{\vert \Lambda_n \vert} \ge -2  \vert \Lambda_k\vert   \int_{2 \vert \Lambda_k \vert }^{\infty} \frac{1}{r}\, d\mathcal{N}(r) \\
	\noalign{\smallskip}
&= -2 \vert \Lambda_k\vert  \left(   -\frac{\mathcal{N}(2\vert \Lambda_k\vert )} {2\vert \Lambda_k\vert } +  \int_{2 \vert \Lambda_k\vert }^{\infty}  \frac{\mathcal{N}(r)}{r^2} \,dr \right) \geq -2 \vert \Lambda_k\vert  \int_{2\vert \Lambda_k\vert}^{\infty}\frac{\mathcal{N}(r)}{r^2} \,dr \\
	\noalign{\smallskip}
&\geq  -2 \vert \Lambda_k\vert   \int_{2\vert \Lambda_k\vert}^{\infty} \frac{\alpha+  p_2  \sqrt{r}}{r^2} \,dr = - 2 \vert \Lambda_k\vert \left( \frac{\alpha}{2\vert \Lambda_k\vert}+ \frac{2 p_2} { \sqrt{2\vert \Lambda_k\vert}} \right) \\
	\noalign{\smallskip}
&=-\alpha- 2 \sqrt{2}  p_2  \sqrt{\vert \Lambda_k\vert}. 
	\end{split}
	\end{equation*}

Coming back to~\eqref{f22} and putting together the previous inequality and inequalities~\eqref{f20} (or~\eqref{f20'} in the real case) and~\eqref{f21},  we conclude that inequality~\eqref{f19} holds. This ends the proof. 
\end{proof}


\subsection{Additional properties and proof of Theorem~\ref{lemme technique}} \label{ss32}

In this paragraph we will prove some additional properties that we will use in the proof of Theorem~\ref{lemme technique}. To this end, we will introduce a ``mollifier'' function and we will construct  the entire function $G_k$ ($ k \ge 1$) in Theorem~\ref{lemme technique} by means of this function and function $f_k$ (see~\eqref{f14}). In order to construct this ``mollifier'' function, we follow the strategy of~\cite{Seidman00,CMV2}

Let us take $T > 0$ and a sequence $\Lambda =\left \{ \Lambda_k \right \}_{k\geq 1} \in \mathcal{L} (\beta,\rho,q,p_0, p_1, p_2,  \alpha) $, with $\beta \in [0, \infty)$, $q \in \N$ and $ \rho, p_0, p_1, p_2, \alpha \in (0, \infty) $  (see Definition~\ref{d1}). With all these data, we fix an integer $N \ge 2$ and we define the sequence $\{ a_k \}_{k \ge 1} \subset (0, \infty)$ given by 
	\begin{equation}\label{f5}
a_k:= \frac{C_{N,T}}{k^2},  \quad \hbox{where } C_{N,T}:=\frac{T}{\displaystyle 2 \sum_{k \ge N} \frac{1}{k^2}},
	\end{equation}
in order to have 
	$$
\sum_{k \ge N}  a_k = \frac{T}{2}.
	$$
Observe that this choice implies
	$$
\frac{1}{N} = \int_{N}^{\infty}\frac{1}{y^2}\,dy \leq  \sum_{k \ge N} \frac{1}{k^2} \leq \int_{N-1}^{+\infty}\frac{1}{y^2}\,dy=\frac{1}{N-1},
	$$
and the estimate
	\begin{equation}\label{f6}
\left ( \frac{N-1}{2} \right )T \leq C_{N,T} \leq \frac{N}{2}T.
	\end{equation}

Consider now the function 
	\begin{equation*}
P_{N , T}(z):= e^{iz\frac{T}{2}}\,\prod_{k\geq N} \cos \paren{a_k z}, \quad z \in \C.
	\end{equation*}

With the previous data,  one has:

\begin{lemma}\label{suitablemollifier}
Under the previous conditions, the following properties hold:
\begin{enumerate}
\item The function $P_{N ,T}$ is entire over $\mathbb{C}$ and satisfies  
	\begin{equation*} 
	\left\{
	\begin{array}{l}
P_{N,T}(0)=1,   \\ 
	\noalign{\smallskip}
\left| P_{N,T}(z)  \right|  \leq 1, \quad \forall z \in\mathbb{C} \hbox{ such that } \Im (z) \geq0, \\ 
	\noalign{\smallskip}
\left| e^{-iz\frac{T}{2}} P_{N ,T}(z) \right| \leq e^{\left| z \right|  \frac{T}{2}}, \quad \forall z \in\mathbb{C}.
	\end{array}
	\right.
	\end{equation*}
\item There exist positive constants $\theta _0>0$ and $\theta_1>0$ (independent of $T$ and $N$) such that 
	\begin{equation} \label{f7}
	\left\{
	\begin{array}{l}
\displaystyle \left ( \frac{C_{N,T}}{\theta_0 } \right )^{\frac{1}{2}} \sqrt{\left | x \right |} +1 \geq N \Longrightarrow \log\left | P_{N,T}(x) \right | \leq \frac{-\theta_1 }{2^3} \left ( \frac{C_{N,T}}{\theta_0 } \right )^{\frac{1}{2}} \sqrt{\left | x \right |}, \\ 
	\noalign{\smallskip}
\displaystyle  \left ( \frac{C_{N,T}}{\theta_0 } \right )^{\frac{1}{2}} \sqrt{\left | x \right |} +1 \leq N \Longrightarrow \log\left | P_{N,T} (x)\right | \leq \frac{-\theta_1}{N^3} \left ( \frac{C_{N,T}}{\theta_0 } \right )^2 \left | x \right |^2.
	\end{array}
	\right.
	\end{equation}
\item There exists a positive constant $\theta_2>0$ (independent of $T$ and $N$)  such that 
	\begin{equation}\label{pnsurlesimaginair}
P_{N,T}(ix) \geq e^{-\theta_2 \sqrt{C_{N,T} x}}, \quad \forall x\geq0.
	\end{equation}
\end{enumerate}
\end{lemma}


For the proof of Lemma~\ref{suitablemollifier}, see~\cite{Seidman00,CMV1,CMV2}. 
%


We are ready to prove the fundamental result stated in Theorem~\ref{lemme technique}. We will follow some ideas of~\cite{BBGBO}.

\begin{proof}[Proof of Theorem~\ref{lemme technique}]
Remember that $T>0$ is given and $\Lambda =\left \{ \Lambda_k \right \}_{k\geq} \in \mathcal{L} (\beta,\rho,q,p_0, p_1, p_2,  \alpha) $, with $\beta \in [0, \infty)$, $ \rho , p_0, p_1, p_2, \alpha \in (0, \infty) $ and $q \in \N$, is a complex sequence. Let us define the function
	\begin{equation} \label{Gkk}
G_k(z):= \frac{1}{\sqrt{2\pi}}\frac{f_k(-i z)\, P_{N,T}(z+ \Im(\Lambda_k))}{f_k(\Lambda_k)\,P_{N,T}(i\, \Re ( \Lambda_k ))}; 
	\end{equation}
(the function $f_k$ is given in~\eqref{f14}). From the properties of the functions $P_{N,T}$ (see Lemma~\ref{suitablemollifier}) and $f_k$ we deduce that the function $G_k$ is well defined and is an entire function over $\C$. In addition,
	$$
G_k(i\Lambda_n ) = \frac{1}{\sqrt{2 \pi}} \delta_{kn}, \quad \forall k,n \ge 1.
	$$
Observe that the function $P_{N,T}$ only has real zeros ($ \ens{a_n}_{n \ge 1}$ is a real sequence) and, then, the sequence $\{\Lambda_n\}_{n\geq1,\, n\neq k}$ are zeros of $G_k$ of multiplicity $1$. This proves item~2 in Theorem~\ref{lemme technique}.

Let us now see that $  e^{-iz\frac{T}{2}}G_k$ satisfies inequality~\eqref{f3}.  From Lemmas~\ref{lemmma3} and~\ref{suitablemollifier}, one has 
	\begin{equation*}
\begin{split} 
\left | e^{-iz\frac{T}{2}} G_k(z) \right | & \leq \frac{e^{{\paren{p_2 \pi + 1} \sqrt{\left | z \right |}+C}}  e^{\vabsolut{z+ \Im( \Lambda_k )} \frac{T}{2}} } {\sqrt{2\pi}  \left | f_k(\Lambda_k ) \right |  \left | P_{N,T}(i \Re(\Lambda_k ))\right |}\\
	\noalign{\smallskip}
&\leq  \frac{e^{ \vabsolut{ \Im(\Lambda_k)} \frac{T}{2}+ \paren{p_2\pi + 1} \sqrt{\left | z \right |}+C }} {\sqrt{2\pi}  \left | f_k(\Lambda_k ) \right |  \left | P_{N,T}(i \Re(\Lambda_k ))  \right |} e^{\left | z \right | \frac{T}{2}}, \quad \forall k\geq1.
\end{split}
	\end{equation*}
If we combine the previous inequality with
	$$
\paren{p_2\pi + 1} \sqrt{\left | z \right |} \le \frac{1}{4 \varepsilon} \paren{p_2\pi + 1}^2 + \varepsilon | z|,
	$$
valid for any $\varepsilon > 0$, we conclude that there exists a positive constant $C'_{N,T,k, \varepsilon }$ such that one has~\eqref{f3}. This proves item~1 in Theorem~\ref{lemme technique}.

Let us prove that $G_k $ belongs to $ L^2(\mathbb{R})$ and satisfies estimate~\eqref{estimatt}. To this end, we will make the following choice of $N$:
	\begin{equation}\label{f4}
2+ \gamma \paren{p_2 \pi + 1}^2  \frac 1T \leq N \leq 4+ \gamma \paren{p_2 \pi + 1}^2  \frac 1T \quad \hbox{with} \quad \gamma = \frac{2^7 \theta_0}{\theta_1^2};
	\end{equation}
($p_2$ is given in assumption~\ref{item7}).

Using~\eqref{p2} and~\eqref{f7}, we have that for $\vert x \vert$ large enough one has
	\begin{equation*}
\left| G_k(x) \right| \leq  \frac{e^{ \paren{p_2 \pi + 1} \sqrt{\left | x \right |}+C -\frac{\theta_1}{2^3}\left ( \frac{C_{N,T}}{\theta_0} \right )^{\frac{1}{2}} \sqrt{\left| x + \Im(\Lambda_k) \right|}}} {\sqrt{2\pi}  \left| f_k(\Lambda_k ) \right|  \left| P_{N,T}(i \Re(\Lambda_k )) \right|}.
	\end{equation*}
Observe that if 
	\begin{equation*}
{p_2 \pi + 1} < \frac{\theta_1}{2^3}\left ( \frac{C_{N,T}}{\theta_0} \right )^{\frac{1}{2}},
	\end{equation*}
then $G_k\in L^2(\mathbb{R})$. In fact, thanks to assumption~\eqref{f4} the previous estimate is satisfied. Indeed, recall that $a_k $ and $ C_{N,T}$ are given in~\eqref{f5} and satisfies~\eqref{f6}. So, from~\eqref{f4}
	$$
\frac{\theta_1}{2^3}\left ( \frac{C_{N,T}}{\theta_0} \right )^{\frac{1}{2}} \ge \frac{\theta_1}{2^3} \left ( \frac{(N-1)T}{2 \theta_0} \right )^{\frac{1}{2}} > {p_2 \pi + 1}.
	$$
This proves $G_{k}\in L^2(\mathbb{R})$.

In what follows, we will estimate $\left\| G_k \right\|_{L^2(\R)}$. First, from the expression of $G_k$ (see~\eqref{Gkk}) and using~\eqref{p2},~\eqref{f6}, \eqref{pnsurlesimaginair} and~\eqref{f4},  one has 
	\begin{equation}\label{f9}
	\begin{split}
\int_{-\infty}^{\infty} \vert G_k(x) \vert^2 \,dx & \leq \frac{e^{2 \theta_2 \sqrt{C_{N,T} \Re(\Lambda_k)} + 2C }}{2\pi \left| f_k(\Lambda_k) \right|^2 }\int_{-\infty}^{\infty} e^{2 \paren{p_2 \pi + 1} \sqrt{\left | x \right |}} \vabsolut{P_{N,T}(x+\Im(\Lambda_k )) }^2 \, dx\\
	\noalign{\smallskip}
 & \leq \frac{e^{2 \theta_2 \sqrt{ \left[ 2T  + \gamma  \paren{p_2 \pi + 1}^2 /2 \right] \Re(\Lambda_k)} + 2C }}{2\pi \left| f_k(\Lambda_k) \right|^2 }\int_{-\infty}^{\infty} e^{2 \paren{p_2 \pi + 1} \sqrt{\left | x \right |}} \vabsolut{P_{N,T}(x+\Im(\Lambda_k )) }^2 \, dx \\
	\noalign{\smallskip}
& := 	\frac{e^{2 \theta_2 \sqrt{ \left[ 2T  + \gamma  \paren{p_2 \pi + 1}^2 /2 \right] \Re(\Lambda_k)} + 2C }}{2\pi \left| f_k(\Lambda_k) \right|^2 } \, I .
	\end{split}
	\end{equation}

Denote 
	\begin{equation*}
\calA_1 := \ens{x \in \R : \vabsolut{x + \Im \paren{\Lambda_k}} < X_{N,T}}, \quad \calA_2 := \ens{x \in \R : \vabsolut{x + \Im \paren{\Lambda_k}} \ge X_{N,T}},
	\end{equation*}
where
	\begin{equation*}
X_{N,T}:= \frac{\theta_0 (N-1)^2}{C_{N,T}}.
	\end{equation*}
Let us first observe that, thanks to inequalities~\eqref{f6} and~\eqref{f4}, it is not difficult to see the property
	\begin{equation}\label{f8}
\frac 12 \theta_0 \left( \frac 1T + \gamma \frac{ \paren{p_2 \pi + 1}^2}{T^2} \right) \le X_{N,T } \le 18 \theta_0 \left( \frac 1T + \gamma \frac{ \paren{p_2 \pi + 1}^2}{T^2} \right), 
	\end{equation}
with $\gamma$ given in~\eqref{f4}.

With the previous notations, we can write  
	\begin{equation*}
	\left\{
	\begin{array}{l}
\displaystyle I = \int_{\calA_1} e^{2 \paren{p_2 \pi + 1} \sqrt{\left | x \right |}} \vabsolut{P_{N,T}(x+\Im(\Lambda_k )) }^2 \, dx + \int_{\calA_2}  e^{2 \paren{p_2 \pi + 1} \sqrt{\left | x \right |}} \vabsolut{P_{N,T}(x+\Im(\Lambda_k )) }^2  \,dx \\
	\noalign{\smallskip}
\displaystyle \phantom{I} := I_1+I_2.
	\end{array}
	\right.
	\end{equation*}

The next objective is to provide an estimate of $I_1$ and $I_2$. To this end, we will use property~\eqref{f7} of Lemma~\ref{suitablemollifier}. Firstly, we estimate $I_1$. We have:
	\begin{equation}\label{f11}
	\left\{
	\begin{array} {l}
\displaystyle I_1  \leq e^{2 \paren{p_2 \pi + 1}\sqrt{\left | \Im(\Lambda_k ) \right |}}  \int_{\calA_1} e^{2 \paren{p_2 \pi + 1} \sqrt{\left| x+ \Im(\Lambda_k )\right|}} e^{\frac{-2\theta_1}{N^3} \frac{C_{N,T}}{\theta_0} \left| x+ \Im(\Lambda_k )\right|^2 } \,dx \\
	\noalign{\smallskip}
\displaystyle \phantom{I_1} \leq e^{2 \paren{p_2 \pi + 1}  \paren{\sqrt{\left | \Im(\Lambda_k )\right |}+ \sqrt{X_{N,T}}}}  \int_{\calA_1} e^{\frac{-2 \theta_1}{N^3} \frac{C_{N,T}}{\theta_0} \left|  x  + \Im(\Lambda_k ) \right|^2 } \,dx \\
	\noalign{\smallskip}
\displaystyle \phantom{I_1} \leq e^{2 \paren{p_2 \pi + 1} \paren{\sqrt{\left | \Im(\Lambda_k )\right |}+ \sqrt{X_{N,T}}} }  \vabsolut{\calA_1} = 2 e^{2 \paren{p_2 \pi + 1} \paren{\sqrt{\left | \Im(\Lambda_k )\right |}+ \sqrt{X_{N,T}}} } X_{N,T} . 
	\end{array}
	\right.
	\end{equation}

Let us now estimate $I_2$. If we denote
	$$
L :=  \frac{\theta_1}{2^2}\sqrt{\frac{C_{N,T}}{\theta_0}} - 2 \paren{p_2 \pi + 1},
	$$
and we use again~\eqref{f9}, we get
	\begin{equation}\label{f12}
	\left\{
	\begin{array}{l}
\displaystyle I_2 \leq e^{2 \paren{p_2 \pi + 1} \sqrt{\left |  \Im(\Lambda_k ) \right |}}  \int_{\calA_2}  e^{  2 \paren{p_2 \pi + 1} \sqrt{\left | x+ \Im(\Lambda_k ) \right |} } \left| P_{N,T}(x+\Im(\Lambda_k )) \right|^2\,dx, \\
	\noalign{\smallskip}
\displaystyle \phantom{I_2} \leq e^{2 \paren{p_2 \pi + 1} \sqrt{\left |  \Im(\Lambda_k ) \right |}} \int_{\calA_2} e^{  -L \sqrt{\left | x+\Im(\Lambda_k ) \right | }} \,dx \le 2 e^{2 \paren{p_2 \pi + 1} \sqrt{\left |  \Im(\Lambda_k ) \right |}} \int_0^\infty e^{  -L \sqrt{ x }} \,dx \\
	\noalign{\smallskip}
\displaystyle \phantom{I_2} = 4 e^{2  \paren{p_2 \pi + 1} \sqrt{\left |  \Im(\Lambda_k ) \right |}} \frac 1{L^2}.
	\end{array}
	\right.
 	\end{equation}

As before and in order to bound $L$, we use again~\eqref{f6} and~\eqref{f4}. Thus,
	\begin{equation*}
	\left\{
	\begin{array}{l}
\displaystyle L \ge \sqrt{\frac{\theta_1^2}{2^5 \theta_0} \paren{T + \frac{2^7\theta_0 (d - 1)^2 }{\theta_1^2}}} - 2  \paren{p_2 \pi + 1} = \sqrt{ \frac{\theta_1^2 T}{2^5 \theta_0} + 2^2  \paren{p_2 \pi + 1}^2} - 2  \paren{p_2 \pi + 1} \\
	\noalign{\smallskip}
\displaystyle \phantom{L} =  \frac{\frac{\theta_1^2 T}{2^5 \theta_0}}{\sqrt{ \frac{\theta_1^2 T}{2^5 \theta_0} + 2^2  \paren{p_2 \pi + 1}^2} + 2  \paren{p_2 \pi + 1}} \ge \frac{\frac{\theta_1^2 T}{2^5 \theta_0}}{2 \sqrt{ \frac{\theta_1^2 T}{2^5 \theta_0} + 2^2  \paren{p_2 \pi + 1}^2} } > 0 ,
	\end{array}
	\right.
	\end{equation*}
and 
	\begin{equation}\label{f10}
\frac{1}{L^2} \le \gamma \left( \frac 1T + \gamma  \paren{p_2 \pi + 1}^2 \frac{1}{T^2} \right),
	\end{equation}
with $\gamma$ given in~\eqref{f4}.

Coming back to~\eqref{f9} and taking into account the inequality
	$$
x \le e^{\sqrt x}, \quad \forall x \ge 0,
	$$
assumption~\ref{item3},~\eqref{f11} and~\eqref{f12}, with $X_{NT}$ and $L $ satisfying~\eqref{f8} and~\eqref{f10}, we deduce
	\begin{equation*}
	\begin{split}
\| G \|_{L^2 (\R )} &\le \frac{e^{C \paren{1 + \sqrt{\paren{ 1 + T + p_2} \Re(\Lambda_k)}} }}{ \left| f_k(\Lambda_k) \right| } e^{C \paren{1 + p_2 } \paren{\sqrt{\Im \paren{\Lambda_k }} + \sqrt{X_{NT}}}} \left( \frac 1T + \gamma  \paren{p_2 \pi + 1}^2 \frac{1}{T^2} \right) \\
	\noalign{\smallskip}
& \le \frac{e^{C \paren{1 + \sqrt{\paren{ 1 + T + p_2 } \Re(\Lambda_k)}} }}{ \left| f_k(\Lambda_k) \right| }  e^{C \paren{1 + p} \paren{ 1 + \Im \paren{\Lambda_k } + \frac{1}{\sqrt T } + \frac{p_2 + 1}{ T } }
} \\
	\noalign{\smallskip}
& \le \frac{e^{C \paren{1 + p_2 + \sqrt{\paren{ 1 + T + p_2} \left|\Lambda_k \right| }} }}{ \left| f_k(\Lambda_k) \right| } e^{C \paren{1 + p_2} \paren{ \sqrt{\left|\Lambda_k \right|} + \frac{p_2 + 1}{ T } } }. 
	\end{split}
	\end{equation*}

Finally, the previous inequality and~\eqref{f19} provide estimate~\eqref{estimatt}, for $G_k(x)$. This ends the proof.
\end{proof}


\section{A lower bound for the norm of arbitrary biorthogonal families: Proof of Theorem \ref{estimbas famillbio}} \label{s4} 

This section will be devoted to prove the result the second main result, Theorem~\ref{estimbas famillbio}, of this paper.  To this end, we will follow some ideas developed by~G\"uichal in~\cite{Gui} (see also~\cite{CMV2}).

Let us consider a sequence $\Lambda = \left \{ \Lambda_k \right \}_{k \geq 1}\subset\C$ satisfying property~\eqref{item6}, for $\nu >0$, and such that $\Lambda \in \mathcal{L} (\beta,\rho,q,p_0, p_1, p_2,  \alpha) $,  for $\beta \in [0, \infty)$, $ \rho , p_0, p_1, p_2, \alpha \in (0, \infty) $ and $q \in \N$. On the other hand, let us also consider $ \{q_k\} _{k\geq1} \subset L^2(0,T; \C) $, a biorthogonal family to $\{e^{-\Lambda_k\,t}\}_{k\geq1}$ in $L^2(0,T; \C)$. 

Associated to the sequence $\Lambda$ we introduce the spaces:
	\begin{equation*}
	\begin{array}{l}
E(\Lambda,T) := \overline{ \hbox{span} \ens{ e_n : n\geq1\}} }^{L^2(0,T;\C )} , \\
	\noalign{\smallskip}
E_k(\Lambda,T)  := \overline{ \hbox{span} \ens{ e_n : n\geq1, n \neq k\}} }^{L^2(0,T; \C )} , \quad \forall k \ge 1. 
	\end{array}
	\end{equation*}
where $e_k$ is the function given by~\eqref{f0}. With this notation, one has:
%


\begin{lemma}\label{l2}
Assume that $\Lambda = \left \{ \Lambda_k \right \}_{k \geq 1}\subset\C$ is a complex sequence satisfying~\eqref{f18} for a positive constant $\delta$. Then, the closed space $E(\Lambda,T)$ is a proper subspace of $L^2(0,T;\C ) $. Moreover, the family of exponentials $\ens{e_k }_{k \ge 1}$ is minimal in $L^2(0,T;\C ) $, that is to say, for every $k \ge 1$, one has 
	$$
e_k \notin E_k (\Lambda,T). 
	$$ 
\end{lemma}

The previous lemma is a well-know result for sequences that satisfy~\eqref{f18} (see for instance~\cite{S}, \cite{Re}, \cite{AKBGBdTJMPA}, \cite{AKBGBdTJFA} and Remark~\ref{r2.4}). 

As a consequence of Lemma~\ref{l2}, we can consider $d_{T,k} > 0 $, the distance between $e^{-\Lambda_k t}$ and  $E_k (\Lambda,T)$, i.e.,
	\begin{equation*}
d_{T,k}^2 =\inf_ {p\in E_k (\Lambda,T) } \norma{e_k - p }^2_{L^2(0,T; \C )} = \int^T_0 \left \vert e^{-\Lambda_k t}-p_k(t) \right \vert^2 \,dt, \quad \forall k \ge 1,
	\end{equation*}
where $p_k \in E_k (\Lambda,T) $ is the orthogonal  projection of the function $e_k (t) = e^{-\Lambda_k t}$ on $E_k (\Lambda,T)$. Observe that the function $p_k$ is characterized by: $p_k \in E_k (\Lambda,T)$ and
	$$
\paren{ e_k - p_k ,  e^{-\Lambda_n t}}_{L^2 (0, T; \C )} = 0, \quad \forall n \ge 1: n \not= k. 
	$$
Thus, if we consider the function $s_k$ given by 
	\begin{equation*}
s_k(t):= \frac{e_k (t) -p_k(t)}{d_{T,k}^2} = \frac{e^{-\Lambda_k t}-p_k(t)}{d_{T,k}^2}, \quad t\in(0,T), \quad \forall k \ge 1,
	\end{equation*}
we deduce that the sequence $\ens{s_k }_{k\geq1 } \subset E(\Lambda,T) $ is biorthogonal to $\ens{e^{-\Lambda_k t} }_{k \geq 1}$ in $L^2(0,T; \C )$. This biorthogonal family is optimal in the following sense: if we consider another biorthogonal family $\ens{\tilde{ q }_k }_{k \geq 1}$ to $\ens{e_k }_{k \geq 1}$ in $L^2(0,T; \C )$,  then $\tilde{ q }_k - s_k \in E(\Lambda,T)^{\perp}$. Since $s_k \in E(\Lambda,T)$, we deduce
	$$
\norma{ \tilde q _k}^2_{L^2(0,T; \C)} =  \norma{ s_k }^2_{L^2(0,T; \C )}+ \norma{ \tilde{q }_k -  s_k }^2_{L^2(0,T; \C)} \geq \| s_k \|^2_{L^2(0,T; \C)} = \frac{1}{d_{T,k}^2} , \quad \forall k \ge 1.
	$$
The previous inequality proves the optimality of the sequence $\ens{s_k }_{k\geq1}$. In particular, 
	\begin{equation*}
\norma{ q_k }^2_{L^2(0,T; \C)} \geq  \frac{1}{d_{T,k}^2},  \quad \forall k\geq1.
	\end{equation*}

The goal now is to obtain an upper bound of $d_{T,k}$, for any $k \ge 1$. From the definition of $d_{T,k}$ we clearly  have
	$$
d_{T,k} \leq \norma{ e_k - p  }_{L^2(0,T; \C )}, \quad \forall p \in E_k (\Lambda,T), \quad \forall k \ge 1 .
	$$
Then,
	\begin{equation}\label{f48}
\norma{ q_k }_{L^2(0,T; \C)} \geq  \frac{1}{d_{T,k}} \ge \frac{1}{\norma{ e_k - p  }_{L^2(0,T; \C )} },  \quad \forall p \in E_k (\Lambda,T), \quad \forall k \ge 1 .
	\end{equation}

In order to obtain~\eqref{lowerbound}, we are going to apply the previous inequality to two appropriate functions $p \in E_k (\Lambda,T) $. Inequality~\eqref{lowerbound} will be a direct consequence of inequality~\eqref{f48}, written for these two functions.


\subsection{A lower bound for the norm of arbitrary biorthogonal families. First part}\label{ss41}
Let us prove that, for any $k \ge 3$, one has
	\begin{equation}\label{lowerbound1}
\norma{q_{k} }_{L^2(0,T; \C)} \geq \frac{6}{\pi^2}  \calB_k \, \mathcal P_k \, e^{\frac{1}{T\nu }}, 
	\end{equation}
where $\mathcal P_k$ and $\calB_k$ are respectively given in~\eqref{Pk} and~\eqref{f49}.

Following~\cite{Gui}, the idea is to construct a particular function $p$ in $E_k (\Lambda,T)$. To this end, let us fix a positive integer $M \geq q + k$, where $q$ is given in assumption~\ref{item5}. On the other hand, let us take
	\begin{equation}\label{f44}
f_ 1 (t) := \sum_{n=1} ^{M+1} A_n e_n ( t ) = \sum_{n=1} ^{M+1} A_n e^{-\Lambda_n t } , \quad t\in(0,T), 
	\end{equation}
with coefficients $ A_1,A_2,..., A_{M+1} \in \C$. Observe that $f_1 \in E_k (\Lambda,T )$ if and only if  $A_k=0$ and, when $A_k \neq 0$ then 
	$$
\frac{ 1 }{A_k }  f_1 (t) = e^{-\Lambda_kt }+ \sum_{n=1}^{k-1}\frac{A_n}{A_k}  e^{-\Lambda_nt }+ \sum_{n=k+1}^{M+1}\frac{A_n}{A_k}  e^{-\Lambda_nt } = e_k (t ) - p( t) , \quad t\in(0,T) .
	$$
Therefore, 
	\begin{equation}\label{f45}
d_{T,k} \leq  \left \| \frac{1}{A_k}f_1 \right \|_{L^2(0,T; \C )}, \quad \forall k \ge 1.
	\end{equation}

One has:
\begin{lemma}\label{l3}
Let us fix $k \ge 1$ and $M \ge q + k$. Let us consider the coefficients $A_1,A_2,..., A_{M+1} \in \C $ given by 
	\begin{equation}\label{coefguichal}
A_n := \frac{1}{\displaystyle \prod_{\substack{ i = 1 \\  i \neq n }}^{M+1} \left( \Lambda_i - \Lambda_n \right)}, \quad 1 \le n \le M + 1. 
	\end{equation}
Then, the function $f_1$ introduced in~\eqref{f44} satisfies
	$$
f_1 (0) = f'_1 (0) = \cdots = f_1 ^{(M-1)}(0) = 0 \hbox{ and } f_1 ^{(M)}(0) =1 .
	$$
\end{lemma}
The previous result is due to  G\"uichal. For a proof we refer to~\cite{Gui} or~\cite[Lemma 4.1]{CMV2}.

The next task will be to estimate $\left \| f_1 \right \|_{L^2(0,T; \C)}$, with $f_1$ and $A_n$, $1 \le n \le M + 1$, respectively  given in~\eqref{f44} and~\eqref{coefguichal}. To this aim, we recall the following technical results:
%


\begin{lemma} \label{l4}
Let $B := \ens{ a_n  }_{1 \leq n \leq r + 1} \subset \C$ be a set of distinct points, $r \geq 1$, and let us fix $g$ an analytic function in a convex domain $\Omega \subset \mathbb{C}$ such that $B \subset \Omega$. 
Then, there exists $\theta  \in [-1,1]$ and $\xi \in \hbox{Conv} \,(B)$, the convex hull of $B$, such that
	$$
\sum_{n = 1}^{r + 1} \frac{g(a_n)}{\displaystyle \prod_{\substack{ a_i \in B  \\ i \neq n }} (a_n-a_i)}= \frac{\theta}{r!} \frac{\partial ^{r} g}{\partial z^{r}}(\xi).
	$$
\end{lemma}

We also have:

\begin{lemma}\label{l5}
The following properties hold:
\begin{enumerate}

\item 
	$
\displaystyle \int_0^T t^{N} e^{- \lambda t} \, dt \leq \frac{2 T^{N + 1} }{N + 1 + \lambda T}
 	$, for any $N \ge 1$ and $\lambda > 0$. 
\item 
	$
\displaystyle \frac{1}{N !} \paren{ \frac{ x }{ 1 + x }}^N e^x \le \sum_{ n = N }^\infty \frac{ x^n }{ n !}
	$, for any $x \ge 0$ and $N \ge 1$. 
\end{enumerate}
 \end{lemma}

Lemma~\ref{l4} is a formula due to Jensen. On the other hand, the proof of Lemma~\ref{l5} can be found in~\cite[Lemma 4.2, Lemma 4.3]{CMV2}. 


Now, using assumption~\eqref{item6}, we can provide an estimate of $\vabsolut{A_k}^{-1}$. One has:
\begin{lemma}\label{l6}
Let us fix $k \ge 1$ and $M \ge q + k$. Then, under the assumptions of Theorem~\ref{estimbas famillbio}, we have
	\begin{equation}\label{f46}
\vabsolut{A_k}^{-1}  \le 
	\left\{
	\begin{array}{ll}
\displaystyle \nu^{M + 2 - q - k } \frac{ 1 }{ \paren{q - 1}!  (2 k + q - 1)!  } \paren{M + 1 - k}! \paren{M + 1 + k}!  \, \calP_k^{-1} , & \hbox{if } k \le q , \\
	\noalign{\smallskip}
\displaystyle \nu^{M -2 (q-1)} \frac{ (2 k - q)! }{ k \left[ \paren{q - 1}! \right]^2 (2 k + q - 1)!  } \paren{M + 1 - k}! \paren{M + 1 + k}!  \, \calP_k^{-1} , & \hbox{if } k > q ,
	\end{array}
	\right.
	\end{equation}
where $A_k$ and $\calP_k$ are respectively given in~\eqref{coefguichal} and~\eqref{Pk}. 
\end{lemma}

\begin{proof} 
The proof is a direct consequence of assumption~\eqref{item6}. Indeed, let us first assume that $k > q $. From the expression of $A_k$ (see~\ref{coefguichal}), we obtain, 
	$$
	\begin{array}{l}
\displaystyle \vabsolut{A_k}^{-1} = \prod_{\substack{ n = 1 \\  n \neq k } }^{M+1} \vabsolut{ \Lambda_k - \Lambda_n } = \prod_{n = 1}^{ k-q } \vabsolut{ \Lambda_k -\Lambda_n } \prod_{ \{ n \ge 1:  \ 1 \leq \left | k-n \right | <q \} } \vabsolut{ \Lambda_k - \Lambda_n } \prod_{ n = q + k }^{M+1} \vabsolut{ \Lambda_k - \Lambda_n } \\
	\noalign{\smallskip}
\displaystyle\phantom{\vabsolut{A_k}^{-1}} :=S_{1,k}\,  \calP_k^{-1} \, S_{2,k},
	\end{array}
	$$
where $\calP_k$ is given in~\eqref{Pk}. 
On the other hand, assumption~\eqref{item6} provides the following estimate
	$$
	\begin{array}{l}
\displaystyle S_{1,k} \,S_{2,k} \leq \prod_{n = 1}^{ k-q } \paren{ \nu \vabsolut{ k^2 - n^2 }} \prod_{ n = q + k }^{M+1} \paren{ \nu \vabsolut{ k^2 - n^2 }} \\
	\noalign{\smallskip}
\displaystyle \phantom{ S_{1,k} \,S_{2,k} } = \nu^{M -2 (q-1)} \frac{ (2 k - q)! }{ k \left[ \paren{q - 1}! \right]^2 (2 k + q - 1)!  } \paren{M + 1 - k}! \paren{M + 1 + k}! \, .
	\end{array}
	$$
Putting both inequalities together we deduce~\eqref{f46} in the case $k > q$. 

We can reason as before in the case $k \le q$. In this case, the first product $S_{1,k}$ in the expression of $\vabsolut{A_k}^{-1}$ does not appear. It is not difficult to deduce the following estimate:
	$$
\displaystyle S_{2,k} \leq \prod_{ n = q + k }^{M+1} \paren{ \nu \vabsolut{ k^2 - n^2 }} = \nu^{M + 2 - q - k } \frac{ 1 }{ \paren{q - 1}!  (2 k + q - 1)!  } \paren{M + 1 - k}! \paren{M + 1 + k}! \, .
	$$
The previous inequality implies~\eqref{f46}  for $k \le q$. This ends the proof. 
\end{proof}


Let us continue with the proof of inequality~\eqref{lowerbound1} when $k \ge 3 $. Observe that we can apply Lemma~\ref{l4} to $f_ 1$ with coefficients $A_n$ given in~\eqref{coefguichal}, $r = M$, $B = \ens{ \Lambda_n  }_{1 \leq n \leq M + 1}$ and $g ( z ) = e^{- t z}$ ($ t \in [0, T]$ is fixed). We obtain,
	$$
f_ 1 (t) = \sum_{n=1}^{M+1} \frac{(-1)^M}{\displaystyle \prod_{\substack{ i = 1 \\  i \neq n } }^{M+1} \left( \Lambda_n - \Lambda_i \right)} e^{-\Lambda_n t } = (-1)^M \frac{\theta }{M!} \frac{\partial ^{M} g}{\partial z^{M} }(\xi) = \frac{\theta t^M }{M!} e^{-t\xi},
	$$
where $ \theta = \theta (t)$ satisfies $\left | \theta \right | \leq 1$ and 
	$$
\xi = \sum_{n=1}^{M+1} \alpha_n \Lambda_n ,\quad \hbox{with} \quad \alpha_n \geq 0 \quad \hbox{and} \quad \sum_{n=1}^{M+1}\alpha_n=1.
	$$
Recall that $\left| \Lambda_1 \right| \leq \left| \Lambda_{2} \right| \leq \cdots \leq \left| \Lambda_{M + 1}  \right|$ (see~\ref{item4}). On the other hand, we can write (see~\eqref{f18}):
	$$
\Re( \xi ) = \sum_{n=1}^{M+1} \alpha_n \Re \paren{ \Lambda_n } \ge \delta \sum_{n=1}^{M+1} \alpha_n \vabsolut{ \Lambda_n } \ge \delta \vabsolut{ \Lambda_1 } , 
	$$
where $\delta > 0 $ is a constant only depending on $\beta$ ($ \delta = 1$ when $\beta = 0$). 
Thus, 
	$$
\left| f_1 (t) \right| = \left| \frac{\theta t^{M} e^{-t \xi}}{M!} \right| \leq \frac{t^M}{M!} e^{-  t \Re ( \xi ) } \leq \frac{t^M}{M!} e^{ - \delta \vabsolut{ \Lambda_1 } t} , \quad \forall t \in [0, T].
	$$

Coming back to~\eqref{f45} with $A_k$ given in~\eqref{coefguichal}, we deduce that 
	$$
\displaystyle 
d_{T,k} \leq \frac{1}{\left | A_k \right |} \norma{f_1 }_{L^2(0,T; \C )} \leq \frac{ 1 }{ M! } \vabsolut{A_k}^{-1 } \paren{\int_0^T t^{2M} e^{- 2 \delta \vabsolut{ \Lambda_1 } t} \, dt}^{1/2} , \quad \forall k \ge 1 .
	$$

Let us introduce the quantity
	\begin{equation}\label{f47}
\calD_k = 
	\left\{
	\begin{array}{ll}
\displaystyle \nu^{ k + q -2} \paren{q - 1}!  \,  (2 k + q - 1)!  \sqrt{\delta \vabsolut{\Lambda_1 } + \frac{1}{2T} } \ \calP_k , & \hbox{if} \quad k \le q , \\
	\noalign{\smallskip}
\displaystyle \nu^{ 2 (q-1)} \left[ \paren{q - 1}! \right]^2  \, \frac{ k (2 k + q - 1)!  }{ (2 k - q)! } \sqrt{\delta \vabsolut{\Lambda_1 } + \frac{1}{2T} } \ \calP_k , & \hbox{if} \quad k \ge q. 
	\end{array}
	\right.
	\end{equation}

Let us first work with $k \ge \max \{3, q \} $. If we use Lemma~\ref{l6} and item~1 of Lemma~\ref{l5}, we deduce
	$$
	\begin{array}{l}
\displaystyle d_{T,k} \leq \frac{ \nu^{M -2 (q-1)} }{ M! }  \frac{ (2 k - q)! }{ k \left[ \paren{q - 1}! \right]^2 (2 k + q - 1)!  } (M+k+1)! (M-k+1)!  \frac{T^{M}\sqrt{2T} }{\sqrt{ 2 M + 1 + 2 \delta  T \vabsolut{\Lambda_1 } }} \calP_k^{-1} \\
	\noalign{\medskip}
\displaystyle \phantom{ d_{T,k} } \le \nu^{ -2 (q-1)} \sqrt{\frac{ 2 T}{ 1 + 2 \delta  T \vabsolut{\Lambda_1 } }} \, \frac{ (2 k - q)! }{ k \left[ \paren{q - 1}! \right]^2 (2 k + q - 1)!  } \calP_k^{-1} \frac{ (M-k+1)! }{ M! } (M+k+1)! \paren{ \nu T}^M \\
	\noalign{\medskip}
\displaystyle \phantom{ d_{T,k} } = \calD_k^{-1} \frac{(M+k+1)!}{M (M-1) \cdots (M-k + 4) (M-k + 3) (M-k + 2)}  \paren{ \nu T}^M \\
	\noalign{\medskip}
\displaystyle \phantom{ d_{T,k} } \le  \calD_k^{-1} \frac{1}{(k+q) (k+q-1) \cdots (q+ 4)} \frac{(M+k+1)! }{(M +1 -k - q)^2} \paren{ \nu T}^M \\
	\noalign{\medskip}
\displaystyle \phantom{ d_{T,k} } = \calD_k^{-1} \frac{(q + 3)!}{(k + q)!} \frac{(M+k+1)! }{(M +1 -k - q)^2} \paren{ \nu T}^M   , \quad \forall k \ge \max \{3, q \}  ,

	\end{array}
	$$
where $\calD_k$ is given in~\eqref{f47}. In the previous inequalities we have used that $k \ge 3$ and $M \ge k + q$. 

Now, if $k$ is such that $3 \le k \le q$, we can argue as before and deduce the same inequality. Summarizing, for any $k \ge 3$, one has  
	$$
d_{T,k} \leq \calD_k^{-1} \frac{(q + 3)!}{(k + q)!} \frac{(M+k+1)! }{(M +1 -k - q)^2} \paren{ \nu T}^M , \quad \forall k \ge 3,
	$$
where $\calD_k$ is given in~\eqref{f47}.

Let us finalize the proof of inequality~\eqref{lowerbound1} when $k \ge 3$. The previous estimate of $d_{T,k}$ and item~2 of Lemma~\ref{l5} allow us to write: 
	\begin{equation*}
	\begin{split}
\frac{1}{d_{T,k}} & = \frac{6}{\pi^2} \sum_{ M = k+q }^\infty \frac{1}{(M + 1 - k -q )^2} \frac{1}{d_{T,k}} \\
	\noalign{\smallskip}
& \ge \frac{6}{\pi^2} \calD_k  \frac{(k + q)!}{(q + 3)!}  \sum_{ M = k+q }^\infty \frac{ 1 }{ (M+k+1)! } \frac{1}{ \paren{ \nu T}^M } =  \frac{6}{\pi^2}  \paren{ \nu T}^{ k + 1 } \calD_k \frac{(k + q)!}{(q + 3)!}  \sum_{ n = 2 k + q + 1 }^\infty \frac{ 1 }{ n ! } \frac{1}{ \paren{ \nu T}^n} \\
	\noalign{\smallskip}
& \ge \frac{6}{\pi^2}  \paren{ \nu T}^{ k + 1 } \calD_k \frac{(k + q)!}{(q + 3)!} \frac{1 }{ \paren{2k + q + 1} ! } \frac{ 1 }{\paren{ 1 + \nu T }^{2 k + q + 1}} e^{\frac {1}{ \nu T }},
	\end{split}
	\end{equation*}
where $\calD_k$ is given in~\eqref{f47}. Coming back to~\eqref{f48}, the previous inequality proves~\eqref{lowerbound1}.

\subsection{A lower bound for the norm of arbitrary biorthogonal families. Second part}\label{ss42}
In order to finish the proof of Theorem~\ref{estimbas famillbio}, let us now show that, for any $k \ge 1$, one has
	\begin{equation}\label{lowerbound2}
\norma{q_{k} }_{L^2(0,T; \C)} \geq \calE_k \, \mathcal P_k, 
	\end{equation}
where $\calE_k$ and $\mathcal P_k$ are respectively given in~\eqref{f49'} and~\eqref{Pk}. The approach is close to the previous one.

Let us introduce the function
	\begin{equation}\label{f59}
f_2 (t ) = \sum_{ \{ n \ge 1: \left | k-n \right | <q \} } \widetilde A_n e_n (t) = \sum_{ \{ n \ge 1: \left | k-n \right | <q \} } \widetilde A_n e^{- \Lambda_n t}, \quad t \in (0, T),
	\end{equation}
with coefficients $\widetilde A_n \in \C$ given by
	\begin{equation}\label{f60}
\widetilde A_n := \frac{1}{\displaystyle \prod_{\substack{\{ i \ge 1:   \left | k-i \right | <q \} \\  i \neq n }}  \left( \Lambda_i - \Lambda_n \right)}, \quad n \ge 1: \left | k-n \right | < q . 
	\end{equation}

Observe that $ \vabsolut{\widetilde A_k} = \calP_k \not= 0 $ ($\calP_k$ is given in~\eqref{Pk}). As in the previous subsection, we can write
	$$
\frac{ 1 }{\widetilde A_k }  f_2 (t) = e^{-\Lambda_kt }+ \sum_{ \{ n \ge 1: 0 < \left | k-n \right | <q \} } \frac{\widetilde A_n}{\widetilde A_k}  e^{-\Lambda_nt } = e_k (t ) - \widetilde p( t) , \quad t\in(0,T) .
	$$
Therefore,
	\begin{equation}\label{f61}
d_{T,k} \leq \norma{e_k - \widetilde p }_{L^2(0,T; \C )} = \calP_k^{-1} \left \| f_2 \right \|_{L^2(0,T; \C )}, \quad \forall k \ge 1.
	\end{equation}

Given $k \ge 1$, we consider the set
	$$
B = \{ \Lambda_n : \vabsolut{k - n} < q \}. 
	$$
and the number $r + 1 = \# B $. It is not difficult to see that
	$$
r = \left\{
	\begin{array}{ll}
k + q -2, & \hbox{if } 1 \le k < q, \\
	\noalign{\smallskip}
2(q - 1), & \hbox{if } 1 \le k \ge q,
	\end{array}
	\right.
	$$
and, therefore $r \ge 1$ ($q \ge 2$). Now, if we apply Lemma~\ref{l4} to $f_2$ with coefficients $\widetilde A_n $ given by~\eqref{f60}, the set $B$, $r$  and $g ( z ) = e^{- t z}$ ($ t \in [0, T]$ is fixed), we deduce
	$$
f_2 ( t ) = \frac{\widetilde \theta  }{ r ! } t^{r } e^{-t \widetilde \xi},
	$$
where $\widetilde \theta = \widetilde \theta (t)$ is such that $\vabsolut{\widetilde \theta} \le 1$ and $\widetilde \xi \in \hbox{Conv}\,( B ) $, i.e.,
	$$
\widetilde \xi = \sum_{ \{ n \ge 1: \left | k-n \right | <q \} } \widetilde \alpha_n \Lambda_n \quad \hbox{with} \quad \widetilde \alpha_n \ge 0 \quad \hbox{and} \quad \sum_{ \{ n \ge 1: \left | k-n \right | <q \} } \widetilde \alpha_n = 1. 
	$$
The previous expression of $\widetilde \xi$ and assumption~\ref{item4} also allow us to deduce
	$$
\Re  (\widetilde \xi ) = \sum_{ \{ n \ge 1: \left | k-n \right | <q \} } \widetilde \alpha_n \Re \paren {\Lambda_n}  \ge \delta \sum_{ \{ n \ge 1: \left | k-n \right | <q \} } \widetilde \alpha_n \vabsolut{\Lambda_n} \ge \delta \min_{\Lambda_ n \in B }\vabsolut{\Lambda_n  } = \delta \vabsolut{\Lambda_{\max \{ 1, k + 1 - q \} }  } , 
	$$
with $\delta > 0$ as in~\eqref{f18} ($ \delta = 1 $ when $\beta = 0$). Summarizing, we have proved
	$$
\vabsolut{f_2 (t)} \le \frac{ 1 }{ r ! } \vabsolut{t}^{r } e^{- \delta \vabsolut{\Lambda_{k + 1 - q} } t }, \quad \forall t \in [0, T]. 
	$$

Let us finalize the proof of~\eqref{lowerbound2}. To this end, we work with the previous inequality, inequality~\eqref{f61}, item 1 of Lemma~\ref{l5} and the expression of $r$. Thus, if $1 \le k < q$, we obtain
	$$
	\left\{
	\begin{array}{l}
\displaystyle d_{T,k} \leq \frac{ 1 }{ (k + q -2 )! }  \paren{ \int_0^T \vabsolut{t}^{2(k + q -2 )} e^{- 2 \delta \vabsolut{ \Lambda_1 } t } \, dt  }^{1/2} \calP_k^{-1} \\
	\noalign{\smallskip}
\displaystyle \phantom{d_{T,k}} \leq \frac{ 1 }{ (k + q -2 )! } \, \frac{T^{k + q -2 } \sqrt{2 T}}{ \sqrt{ 2(k + q -2 ) + 1 + 2 \delta \vabsolut{\Lambda_1 } T } } \, \calP_k^{-1} .
	\end{array}
	\right.
	$$

Now, if $k \ge q$, $r = 2(q-1)$ and a similar computation provides
	$$
\displaystyle d_{T,k} \leq \frac{ 1 }{ (2q - 2)! } \, \frac{T^{2 (q -1) } \sqrt{2 T}}{ \sqrt{ 4 (q - 1) + 1 + 2 \delta \vabsolut{\Lambda_{k + 1 - q} } T } } \, \calP_k^{-1} .
	$$
Of course, inequality~\eqref{lowerbound2} is a direct consequence of these inequalities and inequality~\eqref{f48}. This finally ends the proof of Theorem~\ref{estimbas famillbio}. \fin

\section{Application to the boundary controllability problem for some parabolic systems}\label{s5}
This section will be devoted to apply Theorems~\ref{principal theorem} and~\ref{estimbas famillbio} to two particular parabolic systems in order to provide some new results on the control cost of the corresponding boundary controllability problem associated to these systems. To be precise, we will revisit the controllability problems analyzed in~\cite{O} and~\cite{GB-SN} and we will prove new estimates of the control cost with respect to the final time of controllability $T > 0$. Some results in this section have been previously announced in~\cite{GB-O}.


\subsection{A $2 \times 2$ linear coupled parabolic system}\label{ss51}
Let us consider the one-dimensional Dirichlet-Laplace operator $\widetilde L_1 := - \partial_{xx} $ with domain $D(\widetilde L_1) = H^{2} \left ( 0,\pi \right )\cap H_{0}^{1}(0,\pi)$. It is well-known that $\paren{\widetilde L_1,D(\widetilde L_1)}$ is a self-adjoint operator and admits a sequence of eigenvalues $ \Lambda_1 = \ens{ \lambda_k^{(1)}}_{ k \ge 1} = \ens{ k^2 }_{k\geq1}$ and normalized eigenfunctions given by 
	$$ 
\varphi_k^{(1) } (x) := \sqrt {\frac 2 \pi } \sin \paren{k  x }, \quad \forall k\geq1, \quad x \in (0, \pi).
	$$

On the other hand, let $Q$ be a given function in  $ L^2(0, \pi)$ and consider the operator $\widetilde L_2 := - \partial_{xx} + Q $ with domain $D(\widetilde L_2)=D(\widetilde L_1)$. Again, $(\widetilde L_2,D(\widetilde L_2))$ is a self-adjoint operator and admits a sequence of increasing eigenvalues $ \Lambda_2 = \ens{ \lambda_k^{(2)}}_{ k \ge 1}$  and a sequence of normalized eigenfunctions $\ens{ \varphi_k^{(2) } }_{ k \ge 1}$  which is an orthonormal basis of $L^2(0,1)$. 

In this section we will revisit the boundary controllability problem of the system
%
%
	\begin{equation}\label{cbord}
	\left\{
	\begin{array}{ll} 
\partial_t y+  L_2 y = 0 & \hbox{in } Q_T := (0,T)\times (0,\pi),\\ 
	\noalign{\smallskip}
y(\cdot ,0) = B v , \quad  y(\cdot ,\pi)=0 & \hbox{on } (0,T),\\ 
	\noalign{\smallskip}
y(0,\cdot ) = y_{0}    &\hbox{in }  (0,\pi),
	\end{array}
	\right.
	\end{equation}
where $y_0 \in H^{-1} \paren{ 0,\pi,\R^2 }$ is the initial datum, $v \in L^2 (0,T )$ is a scalar control, the operator $(L_2 ,D(L_2))$ is given by: 
	\begin{equation}\label{LBD}
L_2:=\left(
	\begin{array}{cc}
- \partial_{xx}  &0\\0& - \partial_{xx} + Q
	\end{array}
	\right) ,
\quad D(L_2) = H^2(0,\pi;\R^2) \cap H_0^1(0,\pi; \R^2) ,
	\end{equation} 
and $B\in \R^2$ is the control vector. It is interesting to observe that we want to control system~\eqref{cbord}, two variables, with a unique control function $v \in L^2 (0,T )$.

For every $y_0 \in H^{-1} \paren{ 0,\pi ; \R^2 }$, system~\eqref{cbord} admits a unique solution defined by transposition, $y$, which satisfies
	$$
y\in L^2 \paren{ Q_T ; \R^2 } \cap C^0 \paren{ [0,T] ;  H^{-1} \paren{ 0,\pi ; \R^2 } }.
	$$

It is well-known that, when $Q \in L^2(0, \pi )$ satisfies
	\begin{equation}\label{f50}
\int_0^\pi Q(x) \, dx = 0,
	\end{equation}
then one has 
	\begin{equation*}
\lambda_k^{(2)} = \lambda_k^{(1)}  + \varepsilon_k =  k^2 + \varepsilon_k , \quad \forall k \ge 1 ,
	\end{equation*}
with $\ens{ \varepsilon_k }_{k \ge 1} \in \ell^2$. In particular, $\lim \varepsilon_k = 0$ (see for instance~\cite{K}). Observe that in this case, the eigenvalues of the operator $L$ does not fulfill the gap condition~\eqref{f29} and in this case,  the null controllability of system~\eqref{cbord} has a minimal time $T_0 $ of null controllability which is defined as:
	\begin{equation}\label{f51}
T_0 = \limsup \frac{- \log \vert  \varepsilon_k \vert }{ k^2} \in [0, \infty] .
	\end{equation}
To be precise, one has:
%
%
\begin{theorem}\label{t7}
Let us consider $Q \in L^2(0, \pi)$, a function satisfying $Q \not \equiv 0$ and~\eqref{f50}. Given $T >0$ and  $B = \paren{ b_1 , b_2}^t $, one has
\begin{enumerate}
\item System~\eqref{cbord} is approximately controllable at time $T > 0$ if and only if 
	\begin{equation}\label{f67}
b_1 b_2 \not= 0 \quad \hbox{and} \quad \lambda_k^{(1)} \not= \lambda_n^{(2)} \quad \forall k,n \ge 1. 
	\end{equation}
\item Assume that~\eqref{f67} holds and consider $T_0$ given in~\eqref{f51}. Then
\begin{enumerate} 
\item If $T > T_0$,  system~\eqref{cbord} is null controllable at time $T$.
\item If $T < T_0$,  system~\eqref{cbord} is not null controllable at time $T$. \fin
\end{enumerate}
\end{enumerate}
\end{theorem}

The previous result has been proved in~\cite{O}. In this reference, the author also shows that $T_0$ depends on $Q \in L^2(0, \pi )$ and satisfies this property: given $\tau \in [0, \infty] $, there exists $Q \in L^2 (0, \pi)$ satisfying~\eqref{f50}  such that $T_0 = \tau$. Thus, the minimal time $T_0$ associated to system~\eqref{cbord} could reach any value in the interval $[0, \infty]$. Therefore, there exist coefficients $Q \in L^2 (0,\pi )$ such that the corresponding minimal time of system~\eqref{cbord} satisfies $T_0 > 0$. 


\begin{remark}
The study of the controllability of system~\eqref{cbord} is easier when $Q \in L^2(0, \pi)$ does not satisfy condition~\eqref{f50}, i.e., when
	\begin{equation*}
\int_0^\pi Q(x) \, dx \not= 0. 
	\end{equation*}
In fact, we have the following property: system~\eqref{cbord} is null controllable at time $ T >0$ if and only if the system is approximately controllable at this time, i.e., if and only if~\eqref{f67} holds. In this case, we have that $T_0 = 0$ and the null controllability of the system is valid for any $T > 0$ (see~\cite{O}). On the other hand, it is not difficult to check that we can apply Theorem~\ref{Olive} to  the sequence $\Lambda$. As a consequence, the associated control cost $\mathcal{K} (T) $ for system~\eqref{cbord} can be estimated as follows:
	$$
\exp \left( \frac{C_0} T  \right) \le \mathcal{K}(T) \le \exp \left( \frac{C_1} T  \right), \quad \forall T \in (0, \tau_0) ,
	$$
for appropriate positive constants $C_0$, $C_1$ and $\tau_0$ independent of $T$. \fin
\end{remark}

From now on, let us suppose that~\eqref{f50} and~\eqref{f67} hold. Then, when $T > T_0$, we deduce that  system~\eqref{cbord}  is null controllable at time $T$. So, for any $y_0 \in H^{-1}(0, \pi; \R^2)$, the set 
	$$
\calC_{T}\paren{y_0} :=\left \{ v \in L^2(0,T) : y(T, \cdot )=0 \hbox{ in } (0, \pi), \ y \hbox{ solution of } \eqref{cbord} \right \},
	$$
is non empty and therefore, we can define the control cost of system~\eqref{cbord} in time $T$, $\calK(T)$, when $T > T_0$ (see~\eqref{controocost}).

The positive part of the null controllability result for system~\eqref{cbord} at time $T > 0$ stated in Theorem~\ref{t7} is proved in~\cite{O} by using the moment method. Let us briefly describe this method for system~\eqref{cbord}. 

From the previous assumptions, we deduce that $(L_2, D(L_2))$ is a self-adjoint operator. Its spectrum is given by 
	\begin{equation}\label{f62}
\sigma(L_2) := \Lambda = \Lambda_1 \cup \Lambda_2 = \ens{ \lambda_k^{(1)} , \lambda_k^{(2)} }_{k\geq 1} = \ens{ k^2,  k^2 + \varepsilon_k  }_{k\geq 1} ,
	\end{equation}
and the eigenspaces of $L_2$ corresponding to $\lambda_k^{(1)} $ and $ \lambda_k^{(2)}$ are respectively generated by
	\begin{equation}\label{spectrum}
\phi_k^{(1) }  = 
	\left(
	\begin{array}{c} 
\varphi_k^{(1) } \\ \noalign{\smallskip} 0
	\end{array}
	\right) 
\quad \hbox{and} \quad \phi_k^{(2) }  =
	\left(
	\begin{array}{c} 0 \\ \noalign{\smallskip} \varphi_k^{(2) } 
	\end{array}
	\right),\quad \forall k\geq 1.
	\end{equation}
Moreover, the  sequence $\ens{ \phi_k^{(1) } ,\phi_k^{(2) } }_{ k \ge 1} $  is an orthonormal basis of $L^2 \paren{ 0, \pi; \R^2 }$ and an orthogonal basis of $H_0^1 \paren{ 0, \pi; \R^2 }$ and $H^{-1} \paren{ 0,\pi ; \R^2 }$. 

Using the spectral properties of the operator $L_2$ (see~\eqref{LBD}) we can rewrite the null controllability problem for system~\eqref{cbord} at time $T$ as a moment problem. To be precise, one has: 
%

%
%

\begin{proposition} \label{p9}
Under the previous assumptions, given $y_0 \in H^{-1} \paren{ 0,\pi ; \R^2 }$, the control $v \in L^2(0, T)$ is such that the corresponding solution of~\eqref{cbord} satisfies 
	$$
y (T, \cdot) = 0 \quad \hbox{in} \quad (0, \pi ) , 
	$$
if and only if $v \in L^2(0, T)$ satisfies
	\begin{equation}\label{f64}
	\left\{
	\begin{array}{l}
\displaystyle b_1 \varphi_{k, x}^{( 1 ) } (0) \int_0^T v (T - t) e^{- \lambda_k^{(1)} t } \, dt = - e^{- \lambda_k^{(1)} T } \langle y_0 , \phi_{k}^{(1) }  \rangle_{H^{-1}, H_0^1} , \\
	\noalign{\smallskip}
\displaystyle b_2 \varphi_{k, x}^{( 2 ) } (0) \int_0^T v (T - t) e^{- \lambda_k^{(2)} t } \, dt = - e^{- \lambda_k^{(2)} T } \langle y_0 , \phi_{k}^{(2) }  \rangle_{H^{-1}, H_0^1}.  
	\end{array}
	\right.
	\end{equation}
for any $k \ge 1$, where $\lambda_k^{(i)}$ and $\phi_k^{(i)}$ are respectively given in~\eqref{f62} and~\eqref{spectrum}. \hfill $\Box$
\end{proposition}

For a proof of the previous property, see~\cite{O}.

In fact, when~\eqref{f67} holds and $T > T_0$, $T_0$ given in~\eqref{f51}, the corresponding null controllability problem at time $T$ for system~\eqref{cbord}  (or equivalently, the moment problem stated in Proposition~\ref{p9}) can be explicitly solved as follows (see~\cite{O} for the details). The sequence $\Lambda$ given in~\eqref{f62} satisfies~\eqref{f18}. Therefore, Lemma~\ref{l2} can be applied to deduce the existence of a sequence $\ens{q_k^{(1)}, q_k^{(2)}}_{k \ge 1} \subset L^2(0,T)$ biorthogonal to $\ens{e_k^{(1) },  e_k^{(2) } }_{k \ge 1} \subset L^2(0,T)$, where
	\begin{equation}\label{f0'}
e_k^{(i) } (t) = e^{- \lambda_k^{(i) } t}, \quad \forall t \in (0, T), \quad i= 1,2. 
	\end{equation}
Thus, a formal solution of the moment problem~\eqref{f64} is:
	\begin{equation}\label{f65}
v (t) = \sum_{k \ge 1}  \paren{ e^{- \lambda_k^{(1)} T} m_k^{(1)} q_k^{(1)} (T -t)  + e^{- \lambda_k^{(2)} T} m_k^{(2)} q^{( 2 )}_{k } (T -t) }, \quad \forall t \in (0, T),
	\end{equation}
where
	\begin{equation}\label{f66}
m_k^{(i)} = \frac{- 1}{b_i \varphi_{k, x}^{(i) } (0) } \langle y_0 , \phi_{k}^{(i) }  \rangle_{H^{-1}, H_0^1}  , \quad \forall k \ge 1, \quad i=1,2. 
	\end{equation}
Furthermore, when $T > T_0$, with $T_0$ given in~\eqref{f51}, the series~\eqref{f65} converges absolutely in $L^2(0,T)$ and provides a null control $v \in L^2(0,T)$ which in fact is a solution of the moment problem~\eqref{f64}.

Let us see that we can conveniently choose the sequence $\ens{q_k^{(1)}, q_k^{(2)}}_{k \ge 1} \subset L^2(0,T)$ in order to select a null control for system~\eqref{cbord} associated to $y_0 \in H^{-1} \paren{ 0, \pi; \R^2 }$ with minimal norm in $L^2 (0,T)$. For that purpose, we define (see Section~\ref{s4})
	\begin{equation*}
	\begin{array}{l}
E(\Lambda,T) := \overline{ \hbox{span} \ens{ e_n^{(i) }  : n \geq1, \ i=1,2 \}} }^{L^2(0,T)} , \\
	\noalign{\smallskip}
E^{(1)} _k(\Lambda,T)  := \overline{ \hbox{span} \ens{ e_n^{(1) }  , e_l^{ (2) }  : n\geq1,\  n \neq k , \ l \ge 1 }}^{L^2(0,T )} , \quad \forall k \ge 1\\
	\noalign{\smallskip}
E^{(2)} _k(\Lambda,T)  := \overline{ \hbox{span} \ens{ e_n^{(1) }  , e_l^{ (2) }  : n\geq1, \  l \ge 1, \ l \not= k }}^{L^2(0,T )} , \quad \forall k \ge 1. 
	\end{array}
	\end{equation*}

We have:

%
%

\begin{proposition} \label{p10}
Under the previous assumptions, let us suppose that~\eqref{f67} holds. Let us also consider $T > T_0$ and the sequence of functions
	\begin{equation*}
s_k^{(i)}(t) := \frac{e^{-\lambda^{(i)}_k t} - p^{(i)}_k(t)}{d_{T,k,i}^2}, \quad t\in(0,T), \quad \forall k \ge 1, \quad i = 1,2,
	\end{equation*}
where $d_{T,k,i}$ and $p^{(i)}_k \in E^{(i)}_k (\Lambda,T) $ are defined by
	\begin{equation*}
\displaystyle d_{T,k,i}^2 =\inf_ {p\in E^{(i)}_k (\Lambda,T) } \norma{ e_k^{(i) } - p }^2_{L^2(0,T )} = \int^T_0 \left \vert e^{- \lambda^{(i)}_k t} - p^{(i) }_k(t) \right \vert^2 \,dt, \quad \forall k \ge 1, \quad i = 1,2. 
	\end{equation*}
Then, the sequence $\ens{ s_k^{(1)}, s_k^{(2)}}_{k \ge 1} $ lies in $ E(\Lambda,T) $ and is biorthogonal to $\ens{e_k^{(1) },  e_k^{(2) } }_{k \ge 1} $ in $ L^2(0,T)$ (the function $ e_k^{(i)} $ is given in~\eqref{f0'}). Moreover, given $y_0 \in H^{-1}(0, \pi ; \R^2)$, the control $u \in L^2(0,T)$ given by 
	\begin{equation}\label{f68}
u (t) = \sum_{k \ge 1}  \paren{ e^{- \lambda_k^{(1)} T} m_k^{(1)} s_k^{(1)} (T -t)  + e^{- \lambda_k^{(2)} T} m_k^{(2)} s_{k }^{(2)} (T -t) }, \quad \forall t \in (0, T),
	\end{equation}
where $ m_k^{(i)} $ is given in~\eqref{f66}, satisfies $u \in \calC_{T} \paren{y_0}$, $ \widehat u \in E(\Lambda,T) $ ($\widehat u$ is the function $\widehat u (t) = u( T -t)$, $t \in (0, T)$) and
	$$
\| u \|_{L^2(0,T)} = \inf_{v \in \calC_{T} \paren{y_0} } \|v\|_{L^2(0,T)}. 
	$$
\end{proposition}

\begin{proof}
As said before, under assumption~\eqref{f67}, the sequence $\Lambda$ satisfies~\eqref{f18}. Then, the family 
	$$
\ens{e_k^{(1) },  e_k^{(2) } }_{k \ge 1} \subset L^2(0,T)
	$$ 
is minimal in $L^2(0,T) $. In particular, we deduce that the functions $s_k^{(i)}$ are well defined, live in $E ( \Lambda , T) $, for any $k \ge 1$ and $i=1,2$, and are biorthogonal to $\ens{e_k^{(1) },  e_k^{(2) } }_{k \ge 1} $. These properties together with $T > T_0$ imply that the function $u$ defined in~\eqref{f68} satisfies $\widehat u \in E (\Lambda, T)$ and solves the null controllability problem at time $T$ for system~\eqref{cbord} and $y_0 \in H^{-1 }(0, \pi; \R^2)$, i.e., $u \in \calC_T (y_0)$. 

Let us now consider another null control $v \in \calC_T(y_0)$. Using the equivalence stated in Proposition~\ref{p9} we infer that $v $ satisfies the moment problem~\eqref{f64}. Therefore,
	$$
\int_0^T \left[ v (T-t) - u(T -t) \right] e^{- \lambda_k^{(i)} t } \, dt = 0, \quad \forall k \ge 1 \hbox{ and } i=1,2,
	$$
that is to say, $\widehat v - \widehat u \in E(\Lambda,T)^{\perp}$ ($\widehat v$ and $ \widehat u $ are defined as $\widehat v (t) = v ( T -t)$ and $\widehat u (t) = u( T -t)$, $t \in (0, T)$). Using that $\widehat u \in E(\Lambda,T)$, we deduce
	$$
	\left\{
	\begin{array}{l}
\displaystyle \norma{ v}_{L^2(0,T)} = \norma{\widehat v}_{L^2(0,T)} = \norma{\paren{\widehat v - \widehat u} + \widehat u}_{L^2(0,T)} = \norma{\widehat v - \widehat u}_{L^2(0,T)} + \norma{ \widehat u}_{L^2(0,T)} \\
	\noalign{\smallskip}
\displaystyle \phantom{\norma{ v}_{L^2(0,T)}} \ge \norma{ \widehat u}_{L^2(0,T)} = \norma{ u}_{L^2(0,T)}. 
	\end{array}
	\right.
	$$
The previous inequalities prove the result. This finalizes the proof. 
\end{proof}

Our objective is to apply Theorems~\ref{principal theorem} and~\ref{estimbas famillbio} to system~\eqref{cbord} in a particular case. To this end, let us state a technical result of inverse spectral theory whose proof can be found in~\cite{PT} (see also~\cite{O}):
%


\begin{lemma}\label{l7}
Let us consider $\ens{ \varepsilon_k }_{k \ge 1}$, a sequence in $ \ell^2$. Then, there exists a function $Q \in L^2 (0, \pi)$ satisfying~\eqref{f50} such that 
	$$
\sigma(\widetilde L_2) = \Lambda_2 = \ens{ k^2 + \varepsilon_k}_{k \ge 1} ,
	$$
where $\widetilde L_2 := - \partial_{xx} + Q $ with domain $D(\widetilde L_2) = H^{2} \left ( 0,\pi \right )\cap H_{0}^{1}(0, \pi )$. \fin
\end{lemma}

From now on, we will take 
	$$
\varepsilon_k = e^{-  k^{2\gamma }} , \quad k \ge 1,
	$$
with $\gamma \in (0, \infty) $, and $B = (b_1, b_2)^t$ with $b_1 b_2 \not= 0$. Clearly $\ens{ \varepsilon_k }_{k \ge 1} \in \ell^2$ and we can apply Lemma~\ref{l7}. We will work with the function $Q$ associated to the previous sequence provided by Lemma~\ref{l7} and the corresponding sequences of eigenvalues and orthogonal basis $\Lambda_1$, $\Lambda_2$ and $\ens{ \varphi_k^{(1) } }_{ k \ge 1}$ and $\ens{ \varphi_k^{(2) } }_{ k \ge 1}$ associated to the operators $\widetilde L_1$ and $\widetilde L_2$. With this choice, we consider the parabolic control system~\eqref{cbord} with $L_2$ given in~\eqref{LBD}. 

Observe that the sequence $\Lambda$ of eigenvalues of the operator $L_2$ can be rearranged as an increasing sequence $\Lambda = \ens{ \Lambda_k }_{k \ge 1}$ ($ \gamma \in (0, \infty)$) doing:
	\begin{equation}\label{expression}
\Lambda_{2k-1} = k^2, \quad \Lambda_{2k} = k^2 + e^{- k^{2 \gamma }},\quad \forall  k \geq 1.
	\end{equation}
It is clear that the functions 
	\begin{equation}\label{spectrum2}
\phi_{2k-1}  = 
	\left(
	\begin{array}{c} 
\varphi_k^{(1) } \\ \noalign{\smallskip} 0
	\end{array}
	\right) 
\quad \hbox{and} \quad \phi_{2k}  =
	\left(
	\begin{array}{c} 0 \\ \noalign{\smallskip} \varphi_k^{(2) } 
	\end{array}
	\right),\quad \forall k\geq 1.
	\end{equation}
are an orthonormal basis of eigenfunctions of the operator $L_2$ in $L^2(0, 1; \R^2)$ and an orthogonal basis of $H_0^1 ( 0, \pi ; \R^2)$ and $H^{-1} ( 0, \pi ; \R^2)$.

The controllability properties of system~\eqref{cbord} at time $T>0$ can be deduced from Theorem~\ref{t7}. In this case, system~\eqref{cbord} is approximately controllable for any final time $ T > 0$. The expression of the minimal time is (see~\eqref{f51})
	$$
T_0 = \lim \frac{- \log \paren{ e^{-  k^{2\gamma }} }}{ k^2 } =
	\left\{
	\begin{array}{ll}
0 & \hbox{if } \gamma \in (0, 1), \\
1 & \hbox{if } \gamma = 1, \\
\infty & \hbox{if } \gamma >  1. \\
	\end{array}
	\right.
	$$
We deduce then
\begin{enumerate}
\item If $ \gamma \in (0,1)$, system~\eqref{cbord} is null controllable at any final time $ T > 0$. 
\item If $\gamma = 1$, system~\eqref{cbord} is null controllable at any final time $ T > 1$ and is not controllable at time $T$ when $T < 1$. 
\item When $\gamma > 1$, system~\eqref{cbord} is never null controllable at any final time $ T > 0$. 
\end{enumerate}

Observe that, when $\gamma \in (0, 1 )$ and $Q \in L^2 (0, \pi )$ is the function provided by Lemma~\ref{l7} associated to $\varepsilon_k = e^{-  k^{2\gamma }} $, system~\eqref{cbord} is null controllable at time $T $, for any $T > 0$. We can introduce the control cost $\calK ( T)$ associated to this system (see~\eqref{controocost}). Our objective is to analyze the dependence of $\calK (T)$ with respect to $T$ and $\gamma \in (0,1)$.

First, let us see that the sequence $\Lambda = \ens{\Lambda_k}_{k \ge 1}$ (see~\eqref{expression}) of eigenvalues of the operator $L_2$ is in the class $\mathcal{L} (\beta ,\rho, q, p_0, p_1, p_2,  \alpha)$ (see Definition~\ref{d1}) for appropriate parameters $\beta \ge 0$, $\rho, p_0, p_1,p_2 \in (0, \infty)$ and $q \in \N$. One has:
%


\begin{proposition}\label{p7}
Let us fix $\gamma \in (0, 1 )$ and consider the sequence $\Lambda = \ens{\Lambda_k}_{k \ge 1} $ with $\Lambda_k$ given in~\eqref{expression}, $k \ge 1$. Then, the sequence $\Lambda$ satisfies $\Lambda \in  \mathcal{L} (\beta ,\rho, q, p_0, p_1, p_2,  \alpha) $ and~\eqref{item6}, with $\beta = 0$, $q = 2$, 
	$$
\rho = \frac{ 1 }{16}, \quad p_0 = 1, \quad p_1 = p_2 = 2, \quad \alpha = 2 + \frac 1{ \sqrt{e}} \quad \hbox{and} \quad \nu = \frac 12 \paren{ 1 + \frac 1e}. 
	$$
\end{proposition}

\begin{proof}
The proof of this result is a direct consequence of Proposition~\ref{p5} and Remark~\ref{r5}. Indeed, the sequence $\Lambda$ can be written as  $\Lambda = \Lambda_1 \cup \Lambda_2$ with
	$$
\Lambda_1 = \ens{ \lambda_k^{(1)} }_{k\geq 1} = \ens{ k^2 }_{k\geq 1}  \quad \hbox{and} \quad \Lambda_2 = \ens{ \lambda_k^{(2)} }_{k\geq 1} = \ens{  k^2 + e^{-  k^{2\gamma }} }_{k\geq 1}.  
	$$

It is easy to see that $\Lambda_1 \in  \mathcal{L} \paren{ \beta , \rho_1 , q , \pi_0 , \pi_1 , \pi_2 , \alpha_1 }$ and satisfies~\eqref{item6} with $\beta = 0$, $\rho_1 = 1$, $q = 1$, $\pi_0 = \pi_1 = \pi_2 = 1 $,  $\alpha_1 = 1$ and $\nu = \rho_1 = 1 $. On the other hand, 
	$$
\varepsilon_0 = \sup_{k \ge 1}  \vabsolut{\varepsilon_k } = \sup_{k \ge 1}  e^{-  k^{2\gamma }} = e^{- 1} . 
	$$
In addition, the sequence $\Lambda $ can be explicitly defined by
	\begin{equation*}
\Lambda_k = \left\{
	\begin{array}{ll}
\displaystyle \lambda_{ \ell }^{(1)}, & \hbox{if } k = 2 \ell - 1, \\
	\noalign{\medskip}
\displaystyle \lambda_{ \ell }^{(2)}, & \hbox{if } k = 2 \ell.
	\end{array}
	\right.
	\end{equation*}
(see~\eqref{f40'}) for any $ k \ge 1 $. So, from Proposition~\ref{p5} and Remark~\ref{r5}, we deduce that the sequence $ \Lambda $ lies in $\mathcal{L} (\beta ,\rho, q, \pi_0, p_1, p_2,  \alpha)$ and satisfies~\eqref{item6} with parameters given in the statement of the result. This finalizes the proof. 
\end{proof}

With the previous choice, the sequence $\Lambda$, in particular, satisfies property~\ref{item5} for $q = 2$. In this case, let us see how the term $\calP_k$ (see~\eqref{Pk}) can be estimated. One has:


\begin{proposition}\label{p8}
Let us fix $\gamma \in (0,1) $ and consider the sequence $\Lambda $ given by~\eqref{expression}. Then, $\calP_1 = e$ and 
	$$
	\left\{
	\begin{array}{c}
\displaystyle \frac{1}{ \paren{2n - 1}} e^{ n^{2 \gamma } } \le \calP_{2n - 1} \le \frac{1}{ \paren{2n - 1} - e^{-1}} e^{  n^{2 \gamma } } , \quad \forall n \ge 2, \\
	\noalign{\smallskip}
\displaystyle  \frac{1}{ \paren{2n + 1}} e^{ n^{2 \gamma } } \le \calP_{2n } \le \frac{1}{  \paren{2n + 1} - e^{-1}} e^{ n^{2 \gamma } }  , \quad \forall n \ge 1 , 
	\end{array}
	\right.
	$$
where $\calP_k$ and the sequence $\Lambda$ are respectively given in~\eqref{Pk} and~\eqref{expression}.  
\end{proposition}

\begin{proof}
Let us prove the result when $k = 2n$, with $n \ge 1$. The case $k = 2n - 1 $, with $ n \ge 1$, is similar. From~\eqref{Pk} and~\eqref{expression}, we deduce
	$$
\calP_{2n }^{-1} = \paren{ \Lambda_{2n } - \Lambda_{2n - 1 }} \paren{ \Lambda_{2n + 1} - \Lambda_{2n }} = e^{ - n^{2 \gamma } }  \left[ \paren{2n + 1} - e^{ -  n^{2 \gamma } }\right], \quad \forall n \ge 1. 
	$$
The previous formula provides the proof of the result. 
\end{proof}


\begin{remark}\label{r7}
Let us take a sequence $\ens{\varepsilon_k}_{k \ge 1}$ in $\ell^2$ such that $0 < {\varepsilon_k } < 1$ for any $k \ge 1$. From Lemma~\ref{l7}, there exists a function $Q \in L^2(0, 1)$ such that $\sigma(\widetilde L_2) = \ens{ \pi^2k^2 + \varepsilon_k}_{k \ge 1} $. As before, we can consider the operator $L_2$ associated to system~\eqref{cbord} (see~\eqref{LBD}) and the corresponding sequence of eigenvalues $\Lambda$ given by~\eqref{f62}. In this case, we can repeat the computations in Proposition~\ref{p8} and deduce $\calP_1 = \varepsilon_1^{-1}$ and 
	$$
	\left\{
	\begin{array}{l}
\displaystyle \frac{1}{ \paren{2n - 1}  \varepsilon_n} \le \calP_{2n - 1} \le \frac{1}{ \paren{2n - 2}  \varepsilon_n } , \quad \forall n \ge 2, \\
	\noalign{\smallskip}
\displaystyle  \frac{1}{ \paren{2n + 1}  \varepsilon_n}  \le \calP_{2n } \le \frac{1}{ 2n \, \varepsilon_n}  , \quad \forall n \ge 1 , 
	\end{array}
	\right.
	$$
($\calP_k$ is given in~\eqref{Pk} with $q = 2$). The previous estimates prove that we can construct functions $Q \in L^2(0, 1) $ such that the sequence $\{ \calP_k \}_{k \ge 1 }$ associated to $\sigma(\widetilde L_2) = \ens{ \pi^2k^2 + \varepsilon_k}_{k \ge 1} $ can have any arbitrary explosive behavior.  \fin
\end{remark}


The main results of this section concern the control cost $\calK (T) $ associated to system~\eqref{cbord}. First, let us state a bound from above of the control cost:

\begin{theorem}\label{t8}
Let us fix $\gamma \in (0,1) $ and take the function $Q \in L^2(0, \pi)$ provided by Lemma~\ref{l7} associated to $\varepsilon_k = e^{- k^{2\gamma } }$. If we denote $\mathcal{K} (T) $ the control cost of system~\eqref{cbord} in $L^2(0,T)$ at time $T > 0$, then, there exists a positive constant $\calC$, independent of $\gamma$, such that
	\begin{equation}\label{f52}
\mathcal{K}(T) \le \exp \left[  \calC \paren{ 1 + \frac{1}{T}}  + \frac{\calC }{\paren{ 1 - \gamma } T} + \frac{ 1- \gamma }{T^{\frac{\gamma }{1-\gamma }}} \right] , \quad \forall T >0.
	\end{equation}
\end{theorem}


\begin{proof}
Under the assumptions of the theorem, we can apply Proposition~\ref{p7} and deduce that the sequence $\Lambda = \ens{\Lambda_k}_{k \ge 1} $  ($\Lambda_k$ given in~\eqref{expression}) of eigenvalues of the operator $L_2$ (see~\eqref{LBD}) satisfies $\Lambda \in  \mathcal{L} (\beta ,\rho, q, p_0, p_1, p_2,  \alpha) $ and~\eqref{item6}, with $\beta$, $\rho$, $q$, $p_0$, $p_1$, $p_2$ and $  \alpha$ given in the statement of this result. 

Let us now take $T > 0$. Remember that the minimal time associated to system~\eqref{cbord} is $T_0 = 0$. Therefore, without loss of generality, we can assume that $T \in (0,1)$. Thus, Theorem~\ref{principal theorem} can be applied to $\Lambda$ and we deduce the existence of  a family of functions $\{q_k\} _{k\geq1} \subset L^2(0,T)$, biorthogonal to $\{e_k \}_{k\geq1}$ in $L^2(0,T)$ (for the expression of $e_k$, see~\eqref{f0}) which satisfies~\eqref{bounds}. In particular, there exists a positive constant $\calC$, independent of $\gamma$, such that
	$$
\displaystyle \|q_k \|_{L^2(0,T)} \leq \exp \left[ \calC \paren{ 1 + \sqrt{\vabsolut{\Lambda_k } } + \frac{1}{T}} \right] \mathcal P_k , \quad \forall k \geq 1. 
	$$
If we combine the previous inequality with Proposition~\ref{p8} and~\eqref{expression}, we get
	\begin{equation}\label{f54}
	\left\{
	\begin{array}{l}
\displaystyle \|q_{2k -1} \|_{L^2(0,T)} \leq \exp \left[ \calC \paren{ 1 + k + \frac{1}{T}} + k^{2 \gamma } \right]  , \quad \forall k \geq 1, \\
	\noalign{\smallskip}
\displaystyle \|q_{2k } \|_{L^2(0,T)} \leq \exp \left[ \calC \paren{1 + k + \frac{1}{T}} +  k^{2 \gamma }\right]  , \quad \forall k \geq 1,
	\end{array}
	\right.
	\end{equation}
for a new positive constant $\calC$, independent of $\gamma$. 

Let us prove the result. To this end, we consider $ y_0 \in H^{-1} (0, \pi ; \R^2) $ with 
	$$ 
\| y_0 \|_{ H^{-1} (0, \pi ; \R^2) } \le 1.  
	$$
Using the moment method, in~\cite{O}, the author proves that, taking 
	\begin{equation}\label{f55}
v (t) = \sum_{k \ge 1}  \paren{ e^{- k^2 T} m_k^{(1)} q_{2k - 1} (T -t)  + e^{- \paren{ k^2 + e^{-  k^{2\gamma }} } T} m_k^{(2)} q_{2k } (T -t) }, \quad \forall t \in (0, T),
	\end{equation}
where $m_k^{(i)}$ is given in~\eqref{f66}, one has $v \in L^2 (0, T)$ and the corresponding solution of system~\eqref{cbord}, $y \in C^0 ([0,T]; H^{-1}(0,\pi ; \R^2))$, satisfies $ y( T, \cdot ) = 0$ in $(0, \pi)$ ($\phi_{k}^{(i)} $ is given in~\eqref{spectrum}). In~\cite{O} the author also shows that there exists a positive constant $\calC$ (independent of $k$) such that
	$$
\vabsolut{ m_k^{(i)}  } \le \calC \| y_0 \|_{ H^{-1} (0, \pi ; \R^2) } \le \calC ,  \quad \forall k \ge 1, \quad i=1,2. 
	$$

Coming back to~\eqref{f55} and taking into account~\eqref{f54} and the previous estimate, we deduce
	$$
\| v \|_{L^2 (0, T) } \le e^{ \calC \paren{ 1 + \frac{1}{T}} } \sum_{k \ge 1} e^{- k^2 T + \calC k +  k^{2 \gamma } } ,
	$$
for a new positive constant $\calC$, independent of $\gamma$. Let us now take $\varepsilon \in ( 0, 1/2)$, which will be fixed later. Observe that Young inequality implies
	$$
\calC k \le \varepsilon k^2 T + \frac{\calC^2}{4 \varepsilon T}, \quad \forall k \ge 1,
	$$
and therefore, we can write
	$$
- k^2 T + \calC k + k^{2 \gamma } \le - k^2 T + \varepsilon k^2 T + \frac{\calC^2}{4 \varepsilon T} + k^{2 \gamma } = h_\varepsilon \paren{ k^2 } + \frac{\calC^2}{4 \varepsilon T} - \varepsilon k^2 T, \quad \forall k \ge 1,
	$$
where the function $h_\varepsilon$ is given by 
	$$
h_\varepsilon (x) = - \paren{1 - 2 \varepsilon } T x + x^\gamma,\quad x\in (0, \infty).
	$$
To summarize, the previous calculations provide the following estimate for $\| v \|_{L^2 (0, T) } $:
	\begin{equation}\label{f56}
\| v \|_{L^2 (0, T) } \le  e^{ \calC \paren{ 1 + \frac{1}{T}} } e^{ \frac{\calC^2}{4 \varepsilon T} } \sum_{k \ge 1} e^{ h_\varepsilon ( k^2 )} e^{ - \varepsilon k^2 T},
	\end{equation}
for any $\varepsilon > 0 $. 

It is easy to see that $h_\varepsilon$ possesses an absolute maximum in $(0, \infty)$ at point 
	$$
x^*= \left ( \frac{ \gamma}{1- 2 \varepsilon} \frac1{T} \right )^{\frac{1}{1-\gamma }}.
	$$
Thus, if we take $\varepsilon = (1 - \gamma)/ 2$, we can write
	$$
h_\varepsilon (x) \leq h_\varepsilon (x^*) = \paren{1 - \gamma} \paren{ \frac{ \gamma}{1- 2 \varepsilon}  }^{\frac{\gamma }{1-\gamma }} \frac{1}{T^{\frac{\gamma }{1-\gamma }}} = \frac{ 1- \gamma }{T^{\frac{\gamma }{1-\gamma }}} , \quad \forall x \in (0,\infty). 
	$$

Going back to the formula~\eqref{f56}, we deduce
	$$
\| v \|_{L^2 (0, T) } \le \exp \left[  \calC \paren{ 1 + \frac{1}{T}}  + \frac{\calC^2}{2 \paren{ 1 - \gamma } T} + \frac{ 1- \gamma }{T^{\frac{\gamma }{1-\gamma }}} \right] \sum_{k \ge 1}  e^{ - \frac 12 \paren{ 1 - \gamma }  k^2 T}. 
	$$
Finally, a comparison with Gauss integral gives
	$$
\sum_{k \ge 1}  e^{ - \frac 12 \paren{ 1 - \gamma } k^2 T} \leq \int_0^\infty e^{ - \frac 12 \paren{ 1 - \gamma } T x^2 } \, dx = \frac{\sqrt{2 \pi}}{2}  \frac{1}{\sqrt{ (1 - \gamma)T}} \leq e^{ \frac{1}{2 \paren{ 1 - \gamma } T} },
	$$
and then,
	$$
\| v \|_{L^2 (0, T) } \le \exp \left[  \calC \paren{ 1 + \frac{1}{T}}  + \frac{\calC^2 + 1 }{2 \paren{ 1 - \gamma } T} + \frac{ 1- \gamma }{T^{\frac{\gamma }{1-\gamma }}} \right] . 
	$$
It is clear that, from the previous inequality, we can deduce~\eqref{f52} for a new positive constant $\calC$, independent of $\gamma$. This completes the proof. 
\end{proof}


Our second result provides an estimate from below of the control cost $\calK (T)$ for system~\eqref{cbord} in $L^2 (0,T)$ at the final time $T > 0$. As before, we are going to fix $\gamma \in (0,1) $ and take the function $Q \in L^2(0, \pi)$ provided by Lemma~\ref{l7} associated to $\varepsilon_k = e^{- k^{2\gamma } }$. One has:


\begin{theorem}\label{t9}
Under the assumptions of Theorem~\ref{t8}, there exists two positive constants $\tau_0$ and $\calC$, independent of $\gamma$, such that
	\begin{equation}\label{f53}
\mathcal{K}(T) \ge \calC \exp \left( \frac{\calC} T + \frac{\calC \paren{1 - \gamma}}{T^{\frac\gamma{1 - \gamma}}} \right), \quad \forall T \in (0, \tau_0) .
	\end{equation}
\end{theorem}

Before starting the proof of Theorem~\ref{t9}, we will show a technical result that we will use in its proof:


\begin{lemma}\label{l8}
Let us consider $T > 0$ and $\gamma \in (0, 1)$ and define the function
	$$
\widetilde h (x ) = - T x + x^\gamma, \quad \forall x \in (0, \infty). 
	$$
Let us assume that
	\begin{equation}\label{f58}
T < \gamma \paren{ \frac{ \sqrt 2 - 1}{ \sqrt 2} }^{2 (1 - \gamma)}. 
	\end{equation}
Then, there exists $k_0 \ge 1 $ such that
	$$
\widetilde h \paren{k_0^2 } \ge  \frac{\paren{ 1 + \log 2} }{2e} \frac{\paren{1 - \gamma } }{ T^{\frac{ \gamma }{ 1 - \gamma }} }. 
	$$
\end{lemma}


\begin{proof}
Under the assumptions of the lemma, it is easy to see that the function $ \widetilde h $ is increasing in $( 0 , \widetilde x)$ and decreasing in $(\widetilde x , \infty)$, where
	$$
\widetilde x = \paren{\frac \gamma T}^{\frac{ 1 }{ 1 - \gamma }}.
	$$
Thus, if $ k_0 \ge 1$ is such that
	\begin{equation}\label{f57}
\frac 12 \widetilde x \le k_0^2 \le \widetilde x,
	\end{equation}
then
	$$
\widetilde h \paren{ k_0^2 } \ge \widetilde h (\widetilde x /2) = \paren{ \frac{1}{2^\gamma} - \frac \gamma 2} \paren{\frac \gamma T}^{\frac{ \gamma }{ 1 - \gamma }} > \frac 12 \paren{ 1 + \log 2} \paren{1 - \gamma } \paren{\frac \gamma T}^{\frac{ \gamma }{ 1 - \gamma }} \ge \frac{\paren{ 1 + \log 2} }{2e} \frac{\paren{1 - \gamma } }{ T^{\frac{ \gamma }{ 1 - \gamma }} },
	$$
and we would have the proof of the result. 

In order the finish the proof, let us check that there exists $k_0 \ge 1$ such that~\eqref{f57} holds. Indeed,~\eqref{f57} is equivalent to 
	$$
\frac{1}{ \sqrt2} \sqrt{ \widetilde x} \le k_0 \le \sqrt{\widetilde x}. 
	$$
Observe that this property occurs if 
	$$
\sqrt{\widetilde x} - \frac{1}{ \sqrt2} \sqrt{ \widetilde x} > 1,
	$$
i.e., if $T$ satisfies~\eqref{f58}. This ends the proof. 
\end{proof}

%
%
\begin{proof}[Proof of Theorem~\ref{t9}]
As before, under the assumptions of the theorem, we know that the sequence $ \Lambda$ of eigenvalues of the operator $L_2$ (see~\eqref{LBD}) satisfies $\Lambda \in  \mathcal{L} (\beta ,\rho, q, p_0, p_1, p_2,  \alpha) $ and~\eqref{item6}, with $\beta$, $\rho$, $q$, $p_0$, $p_1$, $p_2$ and $  \alpha$ given in the statement of Proposition~\ref{p7}.

Let us fix $T > 0$. The minimal time $T_0$ for system~\eqref{cbord} associated to the function $Q$ is $T_0 = 0$. In addition, we can apply Proposition~\ref{p10} and Theorem~\ref{estimbas famillbio}. We deduce that the optimal family $\{s_k\} _{k\geq1} \subset E( \Lambda , T)$ biorthogonal to $\{ e_k\}_{k\geq1}$ in $L^2(0,T)$ satisfies~\eqref{lowerbound} ($e_k$ is given in~\eqref{f0}). 

We will divide the proof of the result into two parts:

\smallskip

\textbf{1. }Assume that $\gamma \in (0 , 1/2]$. In this case, it is easy to check that, for any $\tau_0 \in (0, 1] $, one has
	$$
 \frac{1} T \ge \frac{{1 - \gamma}}{T^{\frac\gamma{1 - \gamma}}} , \quad \forall T < \tau_0. 
	$$
Therefore, inequality~\eqref{f53} is equivalent to prove the existence of a positive constant $\calC_0$, independent of $\gamma$, such that
	\begin{equation}\label{f53bis}
\mathcal{K}(T) \ge\calC_0  \exp \left( \frac{\calC_0} T \right), \quad \forall T \in (0, \tau_0) .
	\end{equation}
Our objective is to find $\calC_0 > 0$ and $\tau_0 \in (0, 1]$, independent of $\gamma$, such that one has inequality~\eqref{f53bis}. 

From inequality~\eqref{lowerbound} written for the function $s_3$, we deduce ($ \nu = \frac 12 \paren{ 1 + \frac 1e}$):
	$$
\norma{s_3}_{L^2(0,T)} \ge \frac{6}{\pi^2}  \calB_{3}  \, \calP_ 3 \, e^{\frac{1}{T\nu }}
	$$
 where (see~\eqref{f49} for $q = 2$)
 	$$
\calB_{3} = \displaystyle \calC \frac{ \paren{\nu T}^{4} }{ \paren{1 + \nu T }^{9} }  \sqrt{ \vabsolut{\Lambda_1 } + \frac{1}{2T} },
	$$
and $ \calC $ is a positive constant ($\beta = 0$ and then $ \delta = 1$). From the previous expression, it is not difficult to see that there exist $\calC_0 > 0$ and $\tau_0 \in (0, 1]$, independent of $\gamma$, such that
	$$
\calB_{3} \, \calP_ 3 \ge \calC_0 e^{- \frac{1 }{2 T \nu}}, \quad \forall T \in (0, \tau_0). 
 	$$

Coming back to the expression of $\norma{s_3}_{L^2(0,T)}$, we finally deduce:
	\begin{equation}\label{f69}
\norma{s_3}_{L^2(0,T)} \ge \calC e^{\frac{1 }{2 T \nu}}, \quad \forall T \in (0, \tau_0). 
	\end{equation}

Let us take $y_0 = \phi_3/ \norma{\phi_3}_{H^{-1}}$ (see~\eqref{spectrum2}). Then, applying Proposition~\ref{p10} to $y_0$, it is possible to construct the null control with minimal $L^2$-norm for system~\eqref{cbord} associated to $y_0$ (see~\eqref{f68}):
	$$
u (t) = e^{- 4 T} \frac{1}{b_1 \varphi_{2, x}^{(1) } (0) } \frac{1}{\norma{\phi_3 }_{H^{-1}} } s_{3 } (T -t) , \quad \forall t \in (0, T). 
	$$
From~\eqref{f69}, we also have
	$$
\calK (T ) \ge \inf_{v \in \calC_{T} \paren{y_0} } \|v\|_{L^2(0,T)} = \| u \|_{L^2(0,T)} = \calC \norma{s_3}_{L^2(0,T)} \ge \calC e^{\frac{1 }{2 T \nu}}, \quad \forall T \in (0, \tau_0). 
	$$
This proves inequality~\eqref{f53bis} and inequality~\eqref{f53} when $\gamma \in (0, 1/2]$. 

\smallskip

\textbf{2. }Let us now assume that $\gamma \in (1/2 , 1) $. In this case, inequality~\eqref{f53} is equivalent to
	\begin{equation}\label{f53bbis}
\mathcal{K}(T) \ge\calC_0  \exp \left( \frac{\calC_0}{T^{\frac\gamma{1 - \gamma}}} \right), \quad \forall T \in (0, \tau_0) .
	\end{equation}
and therefore, our goal is to prove that there exist two positive constants $\calC_0$ and $\tau_0$, independent of $\gamma$, in such a way that the previous inequality holds. As before, we are going to work with an appropriate element  $s_{k_0}$ of the optimal biorthogonal family $\{s_k\} _{k\geq1} \subset E( \Lambda , T)$ provided by Proposition~\ref{p10}. 

Let us define $\tau_0$ as
	$$
\tau_0 = \frac 12 \paren{ \frac{ \sqrt 2 - 1}{ \sqrt 2} }. 
	$$
Observe that if $T \in (0, \tau_0)$, then inequality~\eqref{f58} is valid for any $\gamma \in (1/2 , 1)$. From Lemma~\ref{l8}, we can infer the existence of $k_0 \ge 1$ such that
	\begin{equation}\label{f70}
\widetilde h \paren{ k_0^2 } = - k_0^2 T + \paren{ k_0 }^{2 \gamma } \ge  \frac{\paren{ 1 + \log 2} }{2e} \frac{\paren{1 - \gamma } }{ T^{\frac{ \gamma }{ 1 - \gamma }} } =  \frac{\calC \paren{1 - \gamma } }{ T^{\frac{ \gamma }{ 1 - \gamma }} }  . 
	\end{equation}

Consider $y_0 = \phi_{2 k_0 -1} / \norma{\phi_{2 k_0 -1} }_{H^{-1}}$, i.e. (see~\eqref{spectrum2}),
	$$
y_0 (x) = k_0 \sqrt \frac{  2   }{ \pi} 
	\left(
	\begin{array}{c} 
\sin \paren{ k_0 x } \\ \noalign{\smallskip} 0
	\end{array}
	\right). 
	$$
On the other hand, let us also consider the null control for system~\eqref{cbord} associated to $y_0$ provided by Proposition~\ref{p10}:
	$$
u (t) = e^{- k_0^2 T} \frac{1}{b_1 \varphi_{k_0, x}^{(1) } (0) } \langle y_0 , \phi_{2k_0 - 1}  \rangle_{H^{-1}, H_0^1} s_{ 2 k_0 - 1 } (T -t) = \frac{1}{b_1 } \sqrt \frac{  2   }{ \pi} e^{- k_0^2 T}s_{ 2 k_0 - 1 } (T -t)  , \  \forall t \in (0, T). 
	$$

Using inequality~\eqref{lowerbound}, written for the function $s_{2 k_0 - 1 }$, and taking into account Proposition~\ref{p8} ($q = 2$ and $\delta = 1$) and~\eqref{f70}, we deduce 
	\begin{equation*}
	\left\{
	\begin{array}{l}
\displaystyle \norma{ u }_{L^2(0,T)}  = \calC e^{- k_0^2 T} \norma{s_{ 2 k_0 - 1 } }_{L^2(0,T)} \ge \frac{ \calC }{T^2 }  \paren{ \frac{5 }{ 2 T} + \Lambda_{2k_ 0 - 2} }^{1/2}  e^{- k_0^2 T} \, \calP_{2k_0 - 1} \\
	\noalign{\smallskip}
\displaystyle \phantom{ \norma{ u }_{L^2(0,T)} } \ge \frac{ \calC }{T^2 } \paren{ \frac{5 }{ 2 T} + \Lambda_{2k_ 0 - 2} }^{1/2} \frac{1}{2k_0 - 1} \, e^{- k_0^2 T + \paren{  k_0 }^{2 \gamma }}  \\
	\noalign{\smallskip}
\displaystyle \phantom{ \norma{ u }_{L^2(0,T)} } = \frac{ \calC }{T^2 } \paren{ \frac{5 }{ 2 T} + \Lambda_{2k_ 0 - 2} }^{1/2} \frac{  e^{\widetilde h \paren{  k_0^2 } } }{2k_0 - 1}  \ge \frac{ \calC }{T^2 } \exp \paren{ \frac{\calC \paren{1 - \gamma } }{ T^{\frac{ \gamma }{ 1 - \gamma }} } }. 
	\end{array}
	\right.
	\end{equation*}
where $\calC $ is a constant independent of $\gamma $ and $k_0 $. 

As before, 
	$$
\calK (T ) \ge \inf_{v \in \calC_{T} \paren{y_0} } \|v\|_{L^2(0,T)} = \| u \|_{L^2(0,T)} \ge  \frac{ \calC }{T^2 } \exp \paren{ \frac{\calC \paren{1 - \gamma } }{ T^{\frac{ \gamma }{ 1 - \gamma }} } } , \quad \forall T \in (0, \tau_0). 
	$$
Thus, we can conclude that inequality~\eqref{f53bbis} holds. This ends the proof of Theorem~\ref{t9}. 
\end{proof}


\begin{remark}
Observe that inequalities~\eqref{f52} and~\eqref{f53} are valid when $\gamma \in (0,1)$. In fact, these inequalities are equivalent to: 
\begin{enumerate}
\item If $\gamma \in (0, 1/2]$, then, there exist three  positive constants $\tau_0 $, $\calC_0$ and $\calC_1$ (independent of $\gamma$) such that
	$$
\exp \left[ \calC_0 \left( 1 + \frac{1 } T  \right) \right] \le \mathcal{K}(T) \le \exp \left[ \calC_1 \left( 1 + \frac{1 } T  \right) \right], \quad \forall T \in (0, \tau_0) .
	$$
Observe that the previous estimates for the control cost of system~\eqref{cbord} are similar to those obtained for the control cost of the heat equation (see for instance~\cite{Gui} and~\cite{FCZ}).
\item If $\gamma \in (1/2, 1)$, again, there exist three  positive constants $\tau_0 $, $\calC_0$ and $\calC_1$ (independent of $\gamma$) such that
	$$
\exp \left[ \calC_0 \left( 1 + \frac{1}{T^{\frac\gamma{1 - \gamma}}}  \right) \right] \le \mathcal{K}(T) \le \exp \left[ \calC_0 \left( 1 + \frac{1}{T^{\frac\gamma{1 - \gamma}}}  \right) \right]  , \quad \forall T \in (0, \tau_0) .
	$$
The previous expressions prove that the control cost blows up when $\gamma \to 1^-$. This is natural because the minimal time for system~\eqref{cbord} when $\gamma = 1$ is $T_0 = 1$ and the system is not null controllable at time $T$ when $T < 1$. \fin
\end{enumerate}
\end{remark}


\subsection{The linear phase-field system}\label{ss52}
Let us now apply Theorem~\ref{principal theorem} and Theorem~\ref{estimbas famillbio} to the linear version of~\eqref{PFSy} around the constant trajectory $(0,c)$ with $c = 1$ or $c = -1$. To be precise, we will work with the linear system~\eqref{generic} with $L = L_3$ (see~\eqref{f78}) and $\rho, \tau , \xi \in (0, \infty)$. As said above, the  controllability properties of this system has been analyzed in~\cite{GB-SN} under the condition $\xi \not= \frac{1}{j^2}\dfrac{\rho}{\tau}$, for any  $ j \in \N$. The approximate controllability of this system  is given by the next result:
%


%
\begin{theorem}[Approximate controllability] \label{CAproximada_}
Fix $T>0$. Then, system~\eqref{generic} with $L = L_3$ (see~\eqref{f78}) is approximately controllable in $H^{-1} (0, \pi; \R^2)$ at time $T > 0$ if and only if $\lambda_k^{(3,1)} \not= \lambda_n^{(3,2)} $ for any $k,n \ge 1$ (see~\eqref{f78'}), that is to say, if and only if 
	\begin{equation}\label{H2}
{\xi^2\tau^2(\ell^2 - k^2)^2 - 2\xi\rho\tau(\ell^2 + k^2)- 2 \rho-1 \neq 0}, \quad \forall k,\ell\geq1, \quad \ell >k.
	\end{equation}
\end{theorem}

The proof of this result can be found in~\cite{GB-SN}.

Now, our objective is to give a null controllability result at time $T > 0$ for this system when~\eqref{H2} holds (which, in fact, is a necessary condition for the null controllability at time $T$ of system~\eqref{generic} with $L = L_3$) and obtain a bound for the corresponding control cost $\calK (T)$. This problem has analyzed in~\cite{GB-SN} under additional assumptions on the parameters $\xi$, $\rho$ and $\tau$. To be precise, in~\cite{GB-SN} the authors prove:
%


\begin{theorem} \label{CNull_}
Let us us fix $T>0$ and consider $\xi$, $\rho$ and $\tau$, positive real numbers satisfying~\eqref{H2} and
	\begin{equation}\label{H1}
\xi \neq \dfrac{1}{j^2}\dfrac{\rho}{\tau},\quad \forall j\geq1.
	\end{equation}
Then, system~\eqref{generic} with $L = L_3$ (see~\eqref{f78}) is exactly controllable to zero in $H^{-1} (0, \pi; \R^2)$ at time $T>0$. Moreover, there exist  two positive constants $C$ and $M$, only depending on $\xi$, $\rho$ and $\tau$, such that 
	$$
\mathcal{K} (T) \le C e^{M / T}, \quad \forall T>0,
	$$
where $\mathcal{K} (T) $ is the control cost for system~\eqref{generic} with $L = L_3$:
	\begin{equation*}
\mathcal{K}(T) = \sup_{\|y_0\|_{H^{-1} (0,\pi; \R^2 )}=1} \left( \inf_{v\in\mathcal{Z}_{T}(y_0)} \| v \|_{ L^2(0,T)} \right), \quad \forall T>0.
	\end{equation*}
and
	$$
\mathcal{Z}_T(y_0) := \left\{ v \in L^2(0, T) :  y (\cdot , T) = 0 \hbox{ in } (0,\pi), \hbox{ with } y \hbox{ solution of }\eqref{generic} \hbox{ with } L = L_3 \right\}.
	$$
	\end{theorem}

Conditions~\eqref{H2} and \eqref{H1} implies that the sequence $ \Lambda^{(3)} = \ens{\lambda_k^{(3,1)}, \lambda_k^{(3,2)} }_{k\geq1}$ (see~\eqref{f78'}) satisfies the conditions in Theorem~\ref{Olive} (see Remark~\ref{r0}). In fact, condition~\eqref{H1} provides the gap condition~\eqref{f29} for the sequence $\Lambda_3$. Therefore, Theorem~\ref{CNull_} is a consequence of Theorem~\ref{Olive}.

As said before, our objective is to analyze the null controllability of system~\eqref{generic} with $L = L_3$ without imposing condition~\eqref{H1} to the sequence $\Lambda_3$ of eigenvalues of the operator $L_3$.  Let us first see that this sequence is in the class $ \mathcal{L}(\beta,\rho,q,p_0, p_1, p_2,  \alpha)$ with $\beta = 0$ and appropriate parameters $ \rho , p_0, p_1, p_2, \alpha \in (0, \infty) $ (see Definition~\ref{d1}):


\begin{proposition} \label{p12}
Let us consider $\xi$, $\rho$ and $\tau$, positive real numbers satisfying~\eqref{H2}. Then, the sequence $ \Lambda^{(3)} = \ens{\lambda_k^{(3,1)}, \lambda_k^{(3,2)} }_{k\geq1}$, with $\lambda_k^{(3,i)} $ given in~\eqref{f78'}, can be rearranged as a positive increasing sequence $\Lambda^{(3)}  = \ens{\Lambda_k}_{k \ge 1} $ satisfying $\Lambda^{(3)} \in  \mathcal{L} (0 ,\rho, q, p_0, p_1, p_2,  \alpha) $ and~\eqref{item6}, with   

	$$
p_0 = p_1 = p_2 = \frac {2}{\sqrt \xi} \quad \hbox{and} \quad \alpha = \frac 1{2 \sqrt \xi} \left( \sqrt{\dfrac{\rho}{\tau}} + \sqrt{ \dfrac{3\rho+4}{\tau} } \right) + 2 ,
	$$
and $q \ge 2$, $\rho$ and $\nu$ positive constants only depending on $\xi$, $\rho$ and $\tau$.

\end{proposition}


\begin{proof}
The proof of this result is a direct consequence of the results in~\cite{GB-SN}. Indeed, from Proposition~3.2 of~\cite{GB-SN} one has,
	\begin{equation*}
0 < \lambda_{k}^{(3,1)} < \lambda_k^{(3, 2)},\quad \forall k\geq1.		
	\end{equation*}
Secondly, as a consequence of assumption~\eqref{H2} and Theorem~\ref{CAproximada_}, we deduce that the elements of the sequence $\Lambda^{(3)} $ are pairwise different. Thus, this sequence can be rearranged into a positive increasing sequence $\Lambda^{(3)} = \{ \Lambda_k\}_{k\geq1}$ that satisfies~\ref{item1} and, of course,~\ref{item2}, \ref{item3}, with $\beta = 0$, and~\ref{item4}. 

On the other hand, taking into account the proof of Proposition~3.3 in~\cite{GB-SN}, we also have that $\Lambda^{(3)} $ satisfies condition~\ref{item7} in Definition~\ref{d1} with parameters $p_0$, $p_1$, $p_2$ and $\alpha$ as in the statement of the proposition. 

Finally, we deduce properties~\ref{item7} and~\eqref{item6} from Proposition~\ref{p1} with $q$, $\rho$ and $\nu$ given in~\eqref{f30}. This ends the proof of the proposition. 
\end{proof}


In the next result we will provide further properties of the sequence $\Lambda^{(3)} $ that will be used later. Again, we will use some properties that has been proved in~\cite{GB-SN}. One has:


\begin{proposition} \label{p14}
Let us consider $\xi$, $\rho$ and $\tau$, positive real numbers. Then, 
	\begin{equation}\label{f79}
\lambda_k^{(3,2)} - \lambda^{(3, 1)}_{k + i} = \xi \paren{\sqrt{ \frac{\rho}{ \xi \tau } } - i} \paren{2 k + i} + \paren{ \frac{\epsilon_{k + i}} {k + i} + \frac{\epsilon_{k }} {k }  }, \quad \forall k, i \ge 1,
	\end{equation}
where $\lambda_k^{(3,i)} $ is given in~\eqref{f78'} and $\ens{ \epsilon_k }_{k \ge 1}$ is the increasing sequence given by
	\begin{equation}\label{eps}
\epsilon_{k } = \paren{ \frac{\rho + 1}{2 \tau} }^2 \frac{1}{ \sqrt{ \frac{ \xi \rho}{ \tau } + \frac{1}{k^2}  \paren{ \frac{\rho + 1}{2 \tau} }^2}  + \sqrt{ \frac{ \xi \rho}{ \tau } } } , \quad \forall k \ge 1. 
	\end{equation}
\end{proposition}


\begin{proof}
The proof of the result can be found in~\cite{GB-SN} but it is included here for the sake of completeness. From the expression of $r_k$ (see~\eqref{f78''}), we get
	$$
r_k = r_k - \sqrt{ \frac{ \xi \rho}{ \tau} } \, k + \sqrt{ \frac{\xi \rho}{ \tau} } \, k = \frac{ r_k^2 - \frac{\xi \rho}{ \tau} k^2}{ r_k + \sqrt{ \frac{ \xi \rho}{ \tau} } \, k } + \sqrt{ \frac{\xi \rho}{ \tau} } \, k = \sqrt{ \frac{\xi \rho}{ \tau} } \, k + \frac{\epsilon_k }k, \quad \forall k \ge 1,
	$$
(the expression of $\epsilon_k$ is given in~\eqref{eps}). If we take into account the previous expression and~\eqref{f78}, we also deduce
	\begin{equation*}
\lambda_k^{(3,1)} = \xi k^2 + \dfrac{\rho + 1}{2\tau} - \sqrt{\frac{\xi \rho}{\tau}} \, k - \frac{\epsilon_k}{k} ,\quad \lambda_k^{(3,2)} = \xi k^2 + \dfrac{\rho + 1}{2\tau} + \sqrt{\frac{\xi \rho}{\tau}} \, k + \frac{\epsilon_k}{k} , \quad \forall k \ge 1,			
	\end{equation*}
and~\eqref{f79}. This proves the result. 
\end{proof}

Let us now analyze the control cost for the linear phase-field system, i.e., the control cost for system~\eqref{generic} with $L = L_3$. One has:

\begin{theorem}\label{t6}
Let us consider $\xi$, $\rho$ and $\tau$, positive real numbers satisfying~\eqref{H2}. Then, system~\eqref{generic} with $L = L_3$ (see~\eqref{f78}) is exactly controllable to zero at any time $T>0$. Moreover, there exist positive constants $C_0$, $C_1$, $M_0$ and $M_1$ (only depending on $\xi$, $\rho$ and $\tau$) such that 
	\begin{equation}\label{f82}
C_0 e^{M_0/T} \le \mathcal{K} (T) \le C_1 e^{M_1/T}, \quad \forall T \in (0, 1],
	\end{equation}
where $\mathcal{K} (T) $ is the control cost for system~\eqref{generic} with $L = L_3$ defined in the statement of Theorem~\ref{CNull_}.
\end{theorem}

\begin{proof}
The result is proved in~\cite{GB-SN} when the coefficients $\xi$, $\rho$ and $\tau$ satisfy conditions~\eqref{H2} and~\eqref{H1}. Thus, let us prove the result when these coefficients do not satisfy~\eqref{H1}, that is to say, when one has 
	$$ 
\xi = \frac{1}{ j_0^2 }\frac{\rho}{\tau} , 
	$$
for some integer $j_0 \ge 1$. In this case,~\eqref{f79} becomes
	\begin{equation}\label{f79'}
\lambda_k^{(3,2)} - \lambda^{(3, 1)}_{k + i} = \xi \paren{j_0 - i} \paren{ 2 k + i } + \paren{ \frac{\epsilon_{k + i}} {k + i} + \frac{\epsilon_{k }} {k }  }, \quad \forall k, i \ge 1,
	\end{equation}
where $\lambda_k^{ ( 3, i) }$ is given in~\eqref{f78'} and $\ens{ \epsilon_k }_{k \ge 1}$ is the increasing sequence given by~\eqref{eps}. In particular, we can estimate the terms $\epsilon_k$ of the sequence as follows: 
	$$
\paren{ \frac{\rho + 1}{2 \tau} }^2 \frac{1}{ \sqrt{ \frac{ \xi \rho}{ \tau } + \paren{ \frac{\rho + 1}{2 \tau} }^2}  + \sqrt{ \frac{ \xi \rho}{ \tau } } } = \epsilon_1 \le \epsilon_k < \lim_{k \to \infty} \epsilon_k = \paren{ \frac{\rho + 1}{2 \tau} }^2 \frac{{ \sqrt \tau }}{2 \sqrt{ \xi \rho } } := L, \quad \forall k \ge 1. 
	$$
We will use the previous inequalities in what follows. 

If we choose $i$ such that $1 \le i \le j_0 - 1$, from~\eqref{f79'}, we infer
	\begin{equation}\label{ff1}
	\left\{
	\begin{array}{l}
\displaystyle \lambda_k^{(3,2)} - \lambda^{(3, 1)}_{k + i} > \xi \paren{j_0 - i} \paren{ 2 k + i }  , \quad \forall k \ge 1, \\
	\noalign{\smallskip}
\displaystyle \lambda_k^{(3,2)} - \lambda^{(3, 1)}_{k + i} < \xi \paren{j_0 - i} \paren{ 2 k + i } +  \frac{2L } {k } \le \xi \paren{j_0 - 1} \paren{ 2 k + j_0 -1  } +  2L , \quad \forall k \ge 1.
	\end{array}
	\right.
	\end{equation}
Now, if we take $i = j_0$, using again~\eqref{f79'}, we deduce
	\begin{equation}\label{ff2}
	\left\{
	\begin{array}{l}
\displaystyle \lambda_k^{(3,2)} - \lambda^{(3, 1)}_{k + j_0} =  \frac{\epsilon_{k + j_0}} {k + j_0} + \frac{\epsilon_{k }} {k }   > \frac{2 \epsilon_1 }{k + j_0}  , \quad \forall k \ge 1, \\
	\noalign{\smallskip}
\displaystyle \lambda_k^{(3,2)} - \lambda^{(3, 1)}_{k + j_0 } = \frac{\epsilon_{k + j_0}} {k + j_0} + \frac{\epsilon_{k }} {k }  <   \frac{2L } {k } , \quad \forall k \ge 1.
	\end{array}
	\right.
	\end{equation}
Finally, if $i \ge j_0 + 1$, equality~\eqref{f79'}  provides the formula
	\begin{equation*}
\lambda^{(3, 1)}_{k + i} - \lambda_k^{(3,2)} = \xi \paren{i - j_0 } \paren{ 2 k + i } - \paren{ \frac{\epsilon_{k + i}} {k + i} + \frac{\epsilon_{k }} {k }  }, \quad \forall k \ge 1, \quad \forall i \ge j_0 + 1. 
	\end{equation*}
If we take $k_0 \ge 1$ (only depending on $\xi$, $\rho$ and $\tau$) such that
	$$
\frac{2 L}{k_0 } \le \frac \xi2 \left( 2 k_0 + j_0 + 1 \right) ,
	$$
in particular, for any $k \ge k_0 $ and $i \ge j_0 + 1$, one has
	$$
\frac{\epsilon_{k + i}} {k + i} + \frac{\epsilon_{k }} {k } < \frac{2 L}{k} \le \frac{2 L}{k_0} \le \frac \xi2 \left( 2 k_0 + j_0 + 1 \right) \le \frac\xi2 \paren{i - j_0 } \paren{ 2 k + i }  ,
	$$
and 
	\begin{equation}\label{ff3}
	\left\{
	\begin{array}{l}
\displaystyle \lambda^{(3, 1)}_{k + i} - \lambda_k^{(3,2)} > \frac\xi 2 \paren{i - j_0 } \paren{ 2 k + i }  \ge \frac \xi 2 \paren{ 2 k + j_0 + 1 }  , \quad \forall k \ge k_0 , \quad \forall i \ge j_0 + 1, \\
	\noalign{\smallskip}
\displaystyle \lambda^{(3, 1)}_{k + i} - \lambda_k^{(3,2)} < \xi \paren{i - j_0 } \paren{ 2 k + i }  , \quad \forall k \ge k_0 , \quad \forall i \ge j_0 + 1 .
	\end{array}
	\right.
	\end{equation}

The first consequence that we can obtain from~\eqref{ff1}--\eqref{ff3} is the following one: for any $k \ge k_0$, we can write
	$$
\lambda_{k + j_0 }^{(3,1)} < \lambda_{k}^{(3,2)} < \lambda_{k+1+  j_0 }^{(3,1)} < \lambda_{k+1}^{(3,2)}<\cdots , \quad \forall k\geq k_0,
	$$
($\lambda_k^{(3,i)} $ is given in~\eqref{f78'}). Thus, we can give an explicit expression of the elements of the increasing sequence $\Lambda^{(3)} = \ens{\Lambda_k}_{k \ge 1}$ (see Proposition~\ref{p12}): if $1 \le k \le 2 k_0 + j_0 -2$, we define $\Lambda_k$ such that
	$$
\ens{ \Lambda_k }_{1 \le k \le 2 k_0 + j_0 -2} \equiv \ens{ \lambda_k^{(3,1)} }_{1 \le k \le k_0+ j_0 -1} \cup \ens{ \lambda_k^{(3,2)} }_{1 \le k \le k_0 -1}  ,
	$$
and $\Lambda_k < \Lambda_{k+1}$, for any $ k: 1 \le k \le 2 k_0 + j_0 - 3$. From the  $(2k_0 + j_0 - 1)$-th term, we define
	\begin{equation*}
\Lambda_{2k_0 + j_0 + 2 s -1} = \lambda^{(3,1)}_{k_0 + j_0 + s} \quad \hbox{and} \quad
\Lambda_{2k_0 + j_0 + 2 s} = \lambda^{(3,2)}_{k_0 + s }, \quad \forall s \ge 0.
	\end{equation*}
Equivalently, in the case $k \ge 2k_0 + j_0 - 1$, we have
	\begin{equation}\label{Lambda}
	\left\{
	\begin{array}{ll}
\Lambda_{k} = \lambda^{(3,1)}_{\frac 12 \paren{k + j_0 + 1}}, & \hbox{if } k \ge 2k_0 + j_0 - 1 \hbox{ and }  k + j_0 \hbox{ is odd} , \\
	\noalign{\smallskip}
\Lambda_{k} = \lambda^{(3,2)}_{\frac 12 \paren{k - j_0}}, & \hbox{if } k \ge 2k_0 + j_0 - 1 \hbox{ and }  k + j_0  \hbox{ is even} . 
	\end{array}
	\right.	
	\end{equation}

Our next objective will be to obtain appropriate estimates of the products $\calP_k$ (see~\eqref{Pk}) for the sequence $\Lambda^{(3)} $. Remember that $\Lambda^{(3)} \in  \mathcal{L} (0 ,\rho, q, p_0, p_1, p_2,  \alpha) $ and satisfies~\eqref{item6}, with $p_0$, $ p_1$ and $ p_2$ given in Proposition~\ref{p12}, and $q \ge 2 $, $\rho$ and $\nu$ positive constants only depending on $\xi$, $\rho$ and $\tau$. We will reason for arbitrary $k \ge 2 k_0 + j_0 + q - 2$ because if $k$ is such that $1 \le k <2 k_0 + j_0 + q - 2$, taking into account that $\vabsolut{ \Lambda_k - \Lambda_n } > 0$ for any $k \not= n$ (assumption~\eqref{H2}), we deduce the existence of two positive constants $c_0$ and $c_1$ (only depending on $\xi$, $\rho$ and $\tau$) such that
	\begin{equation}\label{f80}
0 < c_0 \le \calP_k \le c_1, \quad \forall k : 1 \le k <2 k_0 + j_0 + q - 2. 
	\end{equation}

Let us then take $k \ge 2 k_0 + j_0 + q - 2$. In this case, if $n \ge 1$ is such that $1\leq \left | k-n \right | <q$, in particular $n \ge 2 k_0 + j_0 -1$. This means that we can use inequalities~\eqref{ff1}--\eqref{ff3} for appropriate indexes.

We will reasoning when $k + j_0$ is odd. A similar argument will provide the proof when $k + j_0 $ is even. Indeed, if $k + j_0$ is odd, from~\eqref{Lambda}, one has $\Lambda_k = \lambda^{(3,1)}_{\widetilde k}$ and  $\Lambda_{k + 1} = \lambda^{(3,2)}_{\widetilde k - j_0 }$ with $\widetilde k = \frac 12 \paren{k + j_0 + 1}$. Thus, we can apply~\eqref{ff2} for $\widetilde k - j_0$ and write
	$$
 \frac{2 \epsilon_1}{ \frac 12 \paren{k + j_0 + 1} } \le \Lambda_{k + 1} - \Lambda_k \le  \frac{2 L }{ \frac 12 \paren{k - j_0 + 1} }. 
 	$$
On the other hand, let us take $n \not= k+1$ with $1\leq \left | k-n \right | <q$. Using properties~\eqref{ff1} and~\eqref{ff3} and the expression of $\lambda_k^{(3, i)}$ (see~\eqref{f78'} and~\eqref{f78''}) and $\Lambda_n$ (see~\eqref{Lambda}), it is not difficult to check the existence of positive constants $c_0$ and $c_1$ (as before, only depending on $\xi$, $\rho$ and $\tau$) such that
	$$
c_0 k \le \vabsolut{ \Lambda_k - \Lambda_n } \le c_1 k, \quad \forall n \not= k+1 \quad \hbox{with} \quad 1\leq \left | k-n \right | <q. 
	$$
As a consequence of the previous inequalities, again, we deduce the existence of positive constants $c_0$ and $c_1$ (only depending on $\xi$, $\rho$ and $\tau$) such that
	$$
c_0 k^{2q -  4 } \le {\displaystyle \prod_{ \{ n \ge 1:  \ 1\leq \left | k-n \right | <q \} } \vabsolut{ \Lambda_k - \Lambda_n } } \le c_1 k^{2q - 4} ,
	$$
or, equivalently (see~\eqref{Pk}),
	\begin{equation}\label{f81}
\frac {c_0}{k^{2q -  4 }} \le \calP_k \le \frac {c_1}{k^{2q -  4 }}  , \quad \forall  k \ge 2 k_0 + j_0 + q - 2,
	\end{equation}
($c_0$ and $c_1$ are new positive constants only depending on $\xi$, $\rho$ and $\tau$). We will use this inequality later.

Let us now revisit some properties on null controllability of system~\eqref{generic} with $L=L_3$ proved in~\cite{GB-SN}: Given $T > 0 $ and $y_0 \in H^{-1}(0,\pi; \R^2 )$, there exists a control $v \in L^2(0,T)$ such that the solution of~\eqref{generic} with $L=L_3$ satisfies $y (\cdot, T) = 0$ in $(0, \pi)$ if and only if $v \in L^2(0, T) $ solves the moment problem
	\begin{equation}\label{f84}
\int_{0}^{T}e^{-\Lambda_k t} v (T -  t ) \, dt = e^{-\Lambda_k T} m_{k}, \quad \forall k\geq 1.
	\end{equation}
In the previous equality $m_k$ only depends on $ y_0$ and satisfies 
	\begin{equation}\label{f83}
| m_k | \le C k \norma{y_0}_{H^{-1}}, \quad \forall k \ge 1,
	\end{equation}
with $C >0$ only depending on $\xi$, $\rho$ and $\tau$. The sequence $\Lambda^{(3)} = \ens{\Lambda_k}_{k \ge 1} = \ens{\lambda_k^{(3,1)}, \lambda_k^{(3,2)}}_{k \ge 1} $  ($\lambda_k^{(3,i)}$ is given in~\eqref{f78'}) provides the eigenvalues of the operator $L_3$ (for the expression of $L_3$, see~\eqref{f78}). 

On the other hand, the real positive sequence $\Lambda^{(3)}$ belongs to $ \mathcal{L} (0 ,\rho, q, p_0, p_1, p_2,  \alpha) $ and satisfies~\eqref{item6} ($p_0$, $ p_1$ and $ p_2$ are given in Proposition~\ref{p12}, $q \ge 2 $, and $\rho$ and $\nu$ are positive constants only depending on $\xi$, $\rho$ and $\tau$). Then, we can apply Theorems~\ref{principal theorem} and~\ref{estimbas famillbio} to the sequence $\Lambda^{(3)}$. We deduce the existence of a biorthogonal family $\ens{ q_k }_{ k \ge 1}$ to the exponentials $\ens{e_k}_{ k \ge 1 }$ (see~\eqref{f0}) associated to the sequence $\Lambda^{(3)}$ satisfying~\eqref{bounds} and~\eqref{lowerbound}.

Let us first prove that, under the assumptions of Theorem~\ref{t6}, system~\eqref{generic} with $L = L_3$ is null controllable at any time $T> 0$ and satisfies the second inequality in~\eqref{f82}. To this end, we will solve the previous moment problem for any $y_0 \in H^{-1}(0,\pi; \R^2 ) $. An explicit solution of this problem is 
	$$
v (t) = \sum_{ k \ge 1} e^{-\Lambda_k T} m_{k} q_k (T - t), \quad \forall t \in (0, T). 
	$$
Since $q_k$, $\calP_k$ and $m_k$ respectively satisfy~\eqref{bounds},~\eqref{f80} or~\eqref{f81}, and~\eqref{f83}, we can prove that the previous series is absolutely convergent in $L^2(0,T)$ and provide an estimate of the $L^2$-norm of $v$. Indeed,
	$$
	\left\{
	\begin{array}{l}
\displaystyle  e^{-\Lambda_k T}  | m_k |  \| q_k \|_{L^2(0,T)} \le C k \,  e^{-\Lambda_k T}  e^{C  \sqrt{\Re(\Lambda _{k})} } e^{C /T} \calP_k \le C e^{-\Lambda_k T}  e^{C  \sqrt{ \Lambda _{k} } } e^{C /T} \norma{y_0}_{H^{-1}} \\
	\noalign{\smallskip}
\displaystyle \phantom{ e^{-\Lambda_k T}   | m_k |  \| q_k \|_{L^2(0,T)} }  \le C  e^{-\Lambda_k T}  e^{C /T}  e^{\frac{C^2}{ 2 T}  + \frac T2 \Lambda_k  } =  C   e^{C /T} e^{- \frac T2 \Lambda_k } , \quad \forall k\geq 1 , 
	\end{array}
	\right.
	$$
for a new positive constant $C $, only depending on $\xi$, $\rho$ and $\tau$. If we use~\eqref{condi**} ($p_2$ and $\alpha$ are given in the statement of Proposition~\ref{p12}), we deduce that $v \in L^2 (0,T)$ and 
	$$
	\left\{
	\begin{array}{l}
\displaystyle \norma{v}_{L^2(0,T)}  \le C e^{C /T} \sum_{k = 1}^\infty e^{- \frac T2 \Lambda_k } \le C e^{C /T} \sum_{k = 1}^\alpha e^{- \frac T2 \Lambda_k }  \sum_{k > \alpha } e^{- \frac T8 \xi (k - \alpha)^2 }  \\
	\noalign{\smallskip}
\displaystyle \phantom{\norma{v}_{L^2(0,T)} } \le C e^{C /T} \int_\R e^{- \frac T8 \xi (x - \alpha)^2 } \, dx = C \sqrt{\frac{ 8 \pi }{ \xi T }} \, e^{C /T}. 
	\end{array}
	\right.
	$$
From this inequality we deduce the estimate from above of $\calK (T)$ in~\eqref{f82}.

Let us now prove the first inequality in~\eqref{f82}. To this end, we will reason as in Subsection~\ref{ss51} and, to be precise, as in Proposition~\ref{p10} and the first point of the proof of Theorem~\ref{t9}. We first construct the sequence $\ens{s_k}_{k \ge 1 }$ biorthogonal to the exponentials $\ens{e_k}_{ k \ge 1 }$ associated to the sequence $\Lambda^{(3)}$. Given $y_0 \in H^{-1} (0, \pi; \R^2 )$, we know that the null control with minimal $L^2$-norm for system~\eqref{generic} with $L = L_3$ (see~\eqref{f78}) associated to $y_0 \in H^{-1} (0, \pi; \R^2)$ is
	$$
u (t) = \sum_{k \ge 1}  e^{- \Lambda_k T} m_k s_k (T -t) , \quad \forall t \in (0, T),
	$$
where $m_k$ depends on $y_0$ and appears in the corresponding moment problem~\eqref{f84}. 

Let us take $\ell = \max \ens{3, q}$ and $y_0 = \Psi_\ell $, with $\Psi_\ell $ the eigenvector of $L_3 $ associated to $\Lambda_\ell $ with $\norma{ \Psi_\ell}_{H^{-1}} = 1$ (for the expression of $\Psi_\ell$ see~Proposition~3.1 in~\cite{GB-SN}). In this case, the corresponding null control with minimal $L^2$-norm is
	$$
u (t) =  e^{- \Lambda_\ell T} m_\ell  s_\ell (T -t) , \quad \forall t \in (0, T),
	$$
and $\calK (T) \ge \norma{ u }_{L^2 (0,T)} = e^{- \Lambda_\ell T} | m_\ell | \norma{s_\ell }$ ($m_\ell$ only depends on $\rho$, $\xi$ and $\tau$). If we use inequalies~\eqref{lowerbound}, for the function $s_\ell$, and~\eqref{f80} or~\eqref{f81} for $k = \ell$, we deduce the existence of a positive constant $C$, only depending on $\rho$, $\xi$ and $\tau$, such that
	$$
\calK (T) \ge C \calB_{\ell } \, e^{\frac{1}{T\nu }} = C \frac{ \paren{\nu T}^{\ell + 1} }{ \paren{1 + \nu T }^{2 \ell + q + 1} }  \sqrt{ \vabsolut{\Lambda_1 } + \frac{1}{2T} } \, e^{\frac{1}{T\nu }} , \quad \forall T > 0. 
	$$
Finally, there exist $ C > 0$, only depending on $\rho$, $\xi$ and $\tau$, such that
	$$
\frac{ \paren{\nu T}^{\ell + 1} }{ \paren{1 + \nu T }^{2 \ell + q + 1} }  \sqrt{ \vabsolut{\Lambda_1 } + \frac{1}{2T} } \ge C \, e^{\frac{-1}{2 T\nu }}, \quad \forall T \in (0, 1]. 
 	$$
Therefore, 
	$$
\calK (T) \ge C \, e^{\frac{1}{2 T\nu }}, \quad \forall T \in (0, 1],
	$$
for a new constant $C > 0$ only depending on $\rho$, $\xi$ and $\tau$. This proves~\eqref{f82} and finalizes the proof of Theorem~\ref{t6}
\end{proof}

Theorem~\ref{t6} in particular provides a local boundary exact controllability result to the trajectory $(0,c)$ ($c = \pm 1$) for the nonlinear system~\eqref{PFSy} under assumption~\eqref{H2}. One has:

%
	\begin{theorem}\label{NullC}
Let us consider $\xi,\tau$ and $\rho$ three  positive numbers satisfying~\eqref{H2}, 
and let us fix $T>0$ and $c=-1$ or $c=1$. Then, there exists $\varepsilon>0$ such that, for any $({\theta}_0,
{\phi}_0)\in H^{-1}(0,\pi)\times (c + H_0^1(0,\pi))$ fulfilling
	\begin{equation*}
\| {\theta}_0\|_{H^{-1}} + \|\tilde{\phi}_0 - c\|_{H_0^1} \le \varepsilon,
	\end{equation*}
there exists $v\in L^2(0,T)$ for which system \eqref{PFSy} has a unique solution 
	$$
( \theta, {\phi})\in \left[L^2(Q_T)\cap C^0([0,T]; H^{-1}(0,\pi; \R^2)) \right] \times C^0(\overline{Q}_T)
	$$
which satisfies
	\begin{equation*}
 \theta (\cdot, T) = 0 \quad \hbox{and} \quad \phi (\cdot, T) = c \quad \hbox{in } (0, \pi).
	\end{equation*}
	\end{theorem}
%
%

In order to obtain the proof of the previous local controllability result for system~\eqref{PFSy}, it is enough to follow the reasoning of the reference~\cite{GB-SN} that combines inequality~\eqref{f82} with the general methodology developed in~\cite{LTT}. For further details, see~\cite{GB-SN}.

\begin{remark}
Theorem~\ref{NullC} is valid under the only assumption~\eqref{H2}. In this sense, Theorem~\ref{NullC} generalizes the local controllability result for system~\eqref{PFSy} stated in~\cite{GB-SN} where the authors prove the same result under assumptions~\eqref{H2} and~\eqref{H1}. \fin
\end{remark}


\section*{Acknowledgement}
Part of this work has been carried out when the second author was developing the research period of the Doc-Course "Partial Differential Equations: Analysis, Numerics and Control" in the Institute of Mathematics of the University of Seville (IMUS). The authors want to thank the IMUS for providing a great working framework. We also want to thank Prof.~Assia Benabdallah for his fruitful discussions that have helped to improve this work.


\appendix
\section{Proof of  Propositions~\ref{p6},~\ref{p3},~\ref{p4} and~\ref{p5}}


\subsection{Proof of Proposition~\ref{p6}} \label{a1}
Let us take $\Lambda = \ens{ \Lambda_k }_{ k \ge 1} \subset (0, \infty)$, a sequence under the assumptions of the proposition. It is clear that the sequence $\Lambda$ satisfies~\ref{item1}--\ref{item4} for $\beta = 0$.

Let us first see that property~\eqref{f27} implies property~\ref{item7}. Indeed, given $ r > 0$, one has $\calN (r) = k$ if and only if $\Lambda_k \le r $ and $\Lambda_{k + 1} > r $. Since the sequence $\Lambda $ satisfies 
	\begin{equation}\label{f26}
\displaystyle \gamma_{0} k + \sqrt{\Lambda_{1 }} - \gamma_0 \le \sqrt{\Lambda_{k }} \le \gamma_{1} k + \sqrt{\Lambda_{1 }} - \gamma_1, \quad \forall k \ge 1,
	\end{equation}
we can write
	\begin{equation*}
\displaystyle \gamma_{0} k + \sqrt{\Lambda_{1 }} - \gamma_0 \le \sqrt{ r } \quad \hbox{and} \quad \sqrt{ r } <  \gamma_{1} ( k + 1) + \sqrt{\Lambda_{1 }} - \gamma_1.
	\end{equation*}
The previous inequalities prove condition~\ref{item7} with $p_0$, $p_1$, $p_2$ and $\alpha$ as in the statement of the proposition.

Let us now see that we can deduce~\ref{item5} from property~\eqref{f27}. First, one has
	$$
\sqrt{\Lambda_{k }} - \sqrt{\Lambda_n} \ge \gamma_{0} \left( k -n \right) , \quad \forall k,n : k \ge n.
	$$
As a direct consequence, one also has
	\begin{equation*}
\displaystyle {\Lambda_{k }} - {\Lambda_n} \ge \gamma_{0} \left( k -n \right) \left(  \sqrt{\Lambda_{k }} + \sqrt{\Lambda_n} \right), \quad \forall k,n : k \ge n, 
	\end{equation*}
that together with~\eqref{f26} provides
	\begin{equation*}
\displaystyle{\Lambda_{k }} - {\Lambda_n} \ge \gamma_{0}^2 \left( k^2 -n^2 \right) + 2 \gamma_0 ( k - n) \paren{ \sqrt{\Lambda_{1 }} - \gamma_0 } , \\
	\end{equation*}
for any $k,n: k \ge n$. If $ \sqrt{\Lambda_{1 }} \ge  \gamma_0 $, clearly one gets~\ref{item5} with $\rho$ as in the statement. Otherwise, $ \sqrt{\Lambda_{1 }} <  \gamma_0 $ and, from the previous inequality, we deduce
	$$
\displaystyle{\Lambda_{k }} - {\Lambda_n} \ge \gamma_{0}^2 \left( k^2 -n^2 \right) - 2 \gamma_0  \frac{ k^2 -n^2 }{k + n}  \paren{  \gamma_0 -  \sqrt{\Lambda_{1 }} } \ge \left( \gamma_{0}^2 - \frac 23 \gamma_0  \paren{  \gamma_0 -  \sqrt{\Lambda_{1 }} } \right) \left( k^2 -n^2 \right) ,
	$$
for any $k,n: k \ge n$. In this case we also deduce~\ref{item5} with $\rho$ given in the statement.

Finally, let us prove~\eqref{item6}. Reasoning as before, we can write
	\begin{equation*}
{\Lambda_{k }} - {\Lambda_n} = \left(  \sqrt{\Lambda_{k }} - \sqrt{\Lambda_n} \right) \left(  \sqrt{\Lambda_{k }} + \sqrt{\Lambda_n} \right) \le \gamma_{1} \left( k -n \right) \left(  \sqrt{\Lambda_{k }} + \sqrt{\Lambda_n} \right), \quad \forall k,n : k \ge n, 
	\end{equation*}
that together with~\eqref{f26} gives
	\begin{equation*}
\displaystyle {\Lambda_{k }} - {\Lambda_n} \le \gamma_{1}^2 \left( k^2 -n^2 \right) + 2 \gamma_1 ( k - n) \paren{ \sqrt{\Lambda_{1 }} - \gamma_1 }  , 
	\end{equation*}
for any $k,n: k \ge n$. In the case in which $ \sqrt{\Lambda_{1 }} \le  \gamma_1 $, we deduce~\eqref{item6} with $\rho = \gamma_1^2$. Otherwise,

	\begin{equation*}
\displaystyle {\Lambda_{k }} - {\Lambda_n} \le \gamma_{1}^2 \left( k^2 -n^2 \right) + 2 \gamma_1 ( k - n) \paren{ \sqrt{\Lambda_{1 }} - \gamma_1 } \le \paren{\gamma_1^2 + \frac 23\gamma_1 \paren{ \sqrt{\Lambda_{1 }} - \gamma_1 } } \left( k^2 -n^2 \right)  , 
	\end{equation*}
for any $k,n: k \ge n$. We also obtain~\eqref{item6} in this case with $\rho$ given in the statement. This finalizes the proof of Proposition~\ref{p6}. \fin
%


\subsection{Proof of Proposition~\ref{p3}}\label{a2}

Let us consider two sequences $\ens{ \lambda_k^{(1)}}_{ k \ge 1}$ and $\ens{ \lambda_k^{(2)}}_{ k \ge 1}$ satisfying~\eqref{f28} and~\eqref{f31}. It is clear that, from~\eqref{f31} and the third condition in~\eqref{f28}, the sequence $ \ens{ \lambda_k^{(1)}}_{ k \ge 1} \cup \ens{ \lambda_k^{(2)}}_{ k \ge 1}$ can be rearranged as an increasing sequence $\Lambda  = \ens{\Lambda_k}_{k \ge 1} $. 

First, let us see that~\eqref{item6} holds and $\Lambda \in  \mathcal{L}(\beta,\rho,q,p_0, p_1, p_2,  \alpha) $ for appropriate positive constants $\rho$, $q$, $p_0$, $p_1$, $p_2$, $\alpha$ and $\nu$. It is clear that $\Lambda$ satisfies \ref{item1}--\ref{item4}. On the other hand, using that $\lambda_k^{(1)} \not= \lambda_n^{(2)}$ for any $k,n \ge 1$, we also have
	$$
\calN (r) =  \# \ens{ k : \lambda_k^{(1)}  \leq r} + \# \ens{ k : \lambda_k^{(2)}  \leq r} := \calN_1 (r) + \calN_2 (r), \quad \forall r>0,
	$$
where $\calN (r)$ is given in~\eqref{counting}. Using the first property in~\eqref{f28} we infer
	$$
\frac{1}{\pi_i^2} k^2 -c_1 k \le \lambda_k^{(i)} \le \frac{1}{\pi_i^2} k^2 + c_1 k , \quad \forall k \ge 1, \quad i=1,2. 
	$$
Therefore, we can follow the arguments in Remark~\ref{r4} (see~\eqref{f43}) and deduce
	$$
-1 - \frac{1}{2} \pi_i^2 c_1 +  \pi_i { \sqrt{ r }} < \calN_i (r)  \le \pi_i { \sqrt{ r }} +  c_1 \pi_i^2  , \quad i=1,2 . 
	$$
Coming back to the expression of $\calN (r) $, we finally obtain 
	$$
-2 - \frac 12 c_1\paren{ \pi_1^2 + \pi_2^2  } + \paren{ \pi_1 + \pi_2  } { \sqrt{ r }} < \calN (r)  \le \paren{ \pi_1 + \pi_2  }  { \sqrt{ r }} +  c_1\paren{ \pi_1^2 + \pi_2^2  }  , \quad \forall r>0. 
	$$
Thus, condition~\ref{item7} holds with $p_0$, $p_1$, $p_2$ and $\alpha $ as in the statement of Proposition~\ref{p3}. Finally, applying Proposition~\ref{p1}, we also have that the sequence $\Lambda$ satisfies \eqref{item6} and $\Lambda \in  \mathcal{L} (\beta,\rho,q,p_0, p_1, p_2,  \alpha) $ with the parameters $ \rho$, $q$, $p_0$, $p_1$, $p_2$, $\alpha$ and $\nu$ given in the statement of Proposition~\ref{p3}.

Let us now check that the gap condition~\eqref{f29} holds. Taking into account property~\eqref{f31}, we just have to check the following property
	$$
\vabsolut{ \lambda_k^{(1)} - \lambda_n^{(2)}} \ge c_2 >0 , \quad \forall k,n \ge 1,
	$$
and  this will be deduced from the third condition in~\eqref{f28}. Indeed, this condition implies 
	$$
\vabsolut{ \lambda_k^{(1)} - \lambda_n^{(2)}} \ge  \frac rk \paren{ \sqrt{ \lambda_k^{(1)} } + \sqrt{\lambda_n^{(2)} }} \ge \frac rk \sqrt{\lambda_k^{(1)} }, \quad \forall k \ge 1. 
	$$
If $k \le 2 c_1 p_1^2 $, from the previous inequality we deduce the existence of a constant $c>0$ such that
	$$
\vabsolut{ \lambda_k^{(1)} - \lambda_n^{(2)}} \ge c , \quad \forall n \in \N.  
	$$
If $k > 2c_1 p_1^2 $, then we can apply the first assumption in~\eqref{f28} and deduce
	$$
\vabsolut{ \lambda_k^{(1)} - \lambda_n^{(2)}} \ge \frac rk \sqrt{\lambda_k^{(1)} } \ge \frac rk \sqrt{ \frac{1}{\pi_1^2 } k^2 - c_1 k } = \frac rk \sqrt{ \frac{k}{\pi_1^2 }\paren{  k- c_1 \pi_1^2} } \ge \frac rk \sqrt{ \frac{k^2}{2 \pi_1^2 } } = \frac{1}{\sqrt{2}} \frac{r}{\pi_1} . 
	$$
This proves~\eqref{f29} and ends the proof of the result.  \fin
%


\subsection{Proof of Proposition~\ref{p4}}\label{a3}
Let us consider the sequence $\Lambda = \ens{k^2}_{k \ge 1} \cup  \ens{d k^2}_{k \ge 1}$ with $d > 0$. Thanks to assumption $\sqrt{d} \not\in \Q$, it is clear that $k^2 \not= d n^2$ for any $k, n \ge 1$. So, the sequence $\Lambda = \ens{k^2}_{k \ge 1} \cup  \ens{d k^2}_{k \ge 1} $ can be rearranged as an increasing sequence $\Lambda = \ens{\Lambda_k}_{k \ge 1}$ that satisfies~\ref{item1}--\ref{item4} with $\beta = 0$. On the other hand, 
	$$
\calN (r) =  \# \ens{ k : k^2  \leq r} + \# \ens{ k : d k^{2}  \leq r} := \lfloor \sqrt r \rfloor + \left\lfloor \frac{\sqrt r}{\sqrt d } \right\rfloor, \quad \forall r>0,
	$$
i.e., 
	$$
- 2 + \left( 1 +  \frac{1}{\sqrt d }  \right) \sqrt r \le \calN (r) \le \left( 1 +  \frac{1}{\sqrt d }  \right) \sqrt r, \quad \forall r>0. 
	$$
Thus, condition~\ref{item7} holds with $p_ 1 = p_2 = p$ and $\alpha$ given in~\eqref{f33}. 

As a direct consequence of Proposition~\ref{p1} we can deduce that condition~\eqref{item6} holds and
	\begin{equation*}
\Lambda \in  \mathcal{L} (\beta,\rho,q,p_0, p_1, p_2,  \alpha),
	\end{equation*}
$q$, $\rho$ and $\nu$ given in~\eqref{f30}. Nevertheless, we will provide better values of parameters $q$, $\rho$ and $\nu$ using the expression of the sequence $\Lambda$. Indeed, if we take $r = \Lambda_k$ with $k \ge 1$, one has $k = \calN (\Lambda_k)$ and
	$$
k = \calN (\Lambda_k) = \left\lfloor \sqrt{ \Lambda_k} \right\rfloor + \left\lfloor \frac{\sqrt{ \Lambda_k}}{\sqrt d } \right\rfloor, \quad \forall k \ge 1.
	$$
Observe that if $\Lambda_k = n_k^2 $ for some $n_k \ge 1$, from the previous inequality we deduce
	$$
k = n_k + \left\lfloor \frac{\sqrt{ \Lambda_k}}{\sqrt d } \right\rfloor. 
	$$
Using that $x- 1 \le \lfloor x \rfloor \le x$, for any $x > 0$, the previous inequality provides,
	$$
 \sqrt{ \Lambda_k} + \frac{\sqrt{ \Lambda_k}}{\sqrt d } - 1 = n_k + \frac{\sqrt{ \Lambda_k}}{\sqrt d } - 1 \le n_k + \left\lfloor \frac{\sqrt{ \Lambda_k}}{\sqrt d } \right\rfloor  = k \le n_k + \frac{\sqrt{ \Lambda_k}}{\sqrt d } = \sqrt{ \Lambda_k} + \frac{\sqrt{ \Lambda_k}}{\sqrt d },
	$$
and 
	\begin{equation}\label{f34}
k \le \left( 1 + \frac{1}{\sqrt d} \right) \sqrt{ \Lambda_k} \le k + 1, \quad \forall k \ge 1.
	\end{equation}
The same property can be proved in the case in which $\Lambda_k = d n_k^2$ for some $n_k \ge 1$.

Let us now prove conditions~\ref{item5}, with $q = 2$, and~\eqref{item6}. If $k - n \ge 2 $, from~\eqref{f34}, one has
	$$
\frac{1}{k^2 - n^2} \left( 1 + \frac{1}{\sqrt d} \right)^2 \left( \Lambda_k - \Lambda_n \right) \ge \frac{k^2 -(n + 1)^2}{k^2 - n^2} = \left( 1+ \frac{1}{k + n} \right) \left( 1 - \frac{1}{k - n} \right) \ge \frac 58.
	$$
Thus,~\ref{item5} holds with $\rho$ given in~\eqref{f33}. On the other hand, if $ k > n$, we deduce (see~\eqref{f34})
	$$
\frac{1}{k^2 - n^2} \left( 1 + \frac{1}{\sqrt d} \right)^2 \left( \Lambda_k - \Lambda_n \right) \le \frac{(k + 1)^2 - n^2}{k^2 - n^2} = \left( 1+ \frac{1}{k + n} \right) \left( 1 + \frac{1}{k - n} \right) \le \frac 83,
	$$
and property~\eqref{item6} with $\nu$ given in~\eqref{f33}. This ends the proof of Proposition~\ref{p4}.  \fin


\subsection{Proof of Proposition~\ref{p5}}\label{a4}
Let us consider two sequences $\Lambda_1 = \ens{ \lambda_k^{(1)}}_{ k \ge 1}$ and $\Lambda_2 = \ens{ \lambda_k^{(2)}}_{ k \ge 1}$ under the conditions of Proposition~\ref{p5}. In particular, the sequence $ \Lambda_1  \cup \Lambda_2 $ can be rearranged as a positive increasing sequence $\Lambda  = \ens{\Lambda_k}_{k \ge 1} $. Let us see that $\Lambda \in  \mathcal{L} (\beta,\rho,q,p_0, p_1, p_2,  \alpha) $, for $\beta =0$ and appropriate positive constants $\rho$, $q$, $p_0$, $p_1$, $p_2 $ and $\alpha$, and~\eqref{item6} holds for $\nu > 0$.
 
First, it is clear that $\Lambda$ satisfies \ref{item1}--\ref{item4} ($\beta = 0$). As above, using that $\lambda_k^{(1)} \not= \lambda_n^{(2)}$ for any $k,n \ge 1$, we also have
	$$
\calN (r) =  \# \ens{ k : \lambda_k^{(1)}  \leq r} + \# \ens{ k : \lambda_k^{(2)}  \leq r} := \calN_1 (r) + \calN_2 (r), \quad \forall r>0.
	$$

From Remark~\ref{r1} we deduce the following property: 
	\begin{equation}\label{f35}
\calN_1(r - \varepsilon_0) \le \calN_2(r) \le \calN_1(r + \varepsilon_0), \quad \forall r >0,
	\end{equation}
(in the previous inequality we have taken $\calN_1(r - \varepsilon_0) = 0$ when $r \le \varepsilon_0$). Indeed, given $r > 0$, if $k_2 = \calN_2(r) $, then $\lambda_{ k_2}^{(2)} \le r $ and $\lambda_{ k_2 + 1 }^{(2)} > r $. In particular,
	$$
\lambda_{ k_2}^{(1)} - \varepsilon_0 \le \lambda_{ k_2}^{(2)} \le r \quad \hbox{and} \quad r < \lambda_{ k_2 + 1 }^{(2)} \le \lambda_{ k_2 + 1 }^{(1)} + \varepsilon_0,
	$$
($\varepsilon_0 = \sup_{k \ge 1 } | \varepsilon_k |$) and $\lambda_{ k_2}^{(1)} \le r + \varepsilon_0$ and $r - \varepsilon_0 <  \lambda_{ k_2 + 1 }^{(1)} $. Applying item~2 of Remark~\ref{r1}, property~\eqref{f35} can be easily deduced. 

Recall that $\Lambda_1 = \ens{ \lambda_k^{(1)}}_{ k \ge 1} \in  \mathcal{L} ( 0 ,\rho_ 1, 1 , \pi_0, \pi_1, \pi_2,  \alpha_1 ) $. Thus, from~\ref{item7}, we deduce
	$$
\pi_1 \sqrt{r} - \alpha_1 \le \calN_1(r) \le \pi_2 \sqrt{r} + \alpha_1, \quad \forall r > 0.
	$$
Combining this inequality and the expression of $\calN (r) $ with~\eqref{f35}, we obtain
	$$
	\left\{
	\begin{array}{ll}
\displaystyle \pi_1 \sqrt{r} - \alpha_1 \le \calN (r) \le \pi_2 \sqrt{r} + \pi_2 \sqrt{r + \varepsilon_0 } + 2\alpha_1, &\hbox{if }  r \le \varepsilon_0,\\
	\noalign{\medskip}
\displaystyle \pi_1 \sqrt{r} + \pi_1 \sqrt{r - \varepsilon_0 } - 2\alpha_1 \le \calN (r) \le \pi_2 \sqrt{r} + \pi_2 \sqrt{r + \varepsilon_0 } + 2\alpha_1, &\hbox{if }  r > \varepsilon_0 .
	\end{array}
	\right.
	$$
Now, from the previous property and taking into account inequalities~\eqref{f71}, it is easy to deduce that $  \calN (r)$ satisfies
	$$
	\left\{
	\begin{array}{ll}
\displaystyle \pi_1 \sqrt{r} - \alpha_1 \le \calN (r) \le 2 \pi_2 \sqrt{r} + \pi_2 \sqrt{\varepsilon_0 } + 2\alpha_1, &\hbox{if }  r \le \varepsilon_0,\\
	\noalign{\medskip}
\displaystyle 2 \pi_1 \sqrt{r} - \pi_1 \sqrt{\varepsilon_0 } - 2\alpha_1 \le \calN (r) \le 2 \pi_2 \sqrt{r} + \pi_2 \sqrt{\varepsilon_0 } + 2\alpha_1, &\hbox{if }  r > \varepsilon_0 .
	\end{array}
	\right.
	$$
In particular,
	$$
\displaystyle 2 \pi_1 \sqrt{r} - \pi_1 \sqrt{\varepsilon_0 } - 2 \alpha_1 \le \calN (r) \le 2 \pi_2 \sqrt{r} + \pi_2 \sqrt{\varepsilon_0 } + 2\alpha_1,  \quad \forall r > 0 .
	$$
Therefore, condition~\ref{item7} holds with $p_1$, $p_2$ and $\alpha $ as in the statement of Proposition~\ref{p5}. 

Let us now see that the sequence $\Lambda$ satisfies~\ref{item5} and~\eqref{item6} with $q =2$ and appropriate positive parameters $ \rho$ and $\nu$. To this end, we will use that $\Lambda_1 \in  \mathcal{L} (0, \rho_1 , 1, \pi_0, \pi_1, \pi_2, \alpha_1)$ ($q=1$) and satisfies~\eqref{item6}, for $\nu_1 \in (0, \infty )$ or, more precisely, we will use
	\begin{equation}\label{f36}
\rho_1 \paren{k^2 -n^2} \le \lambda_{ k}^{(1)} - \lambda_{ n}^{(1)} \le \nu_1 \paren{k^2 -n^2} , \quad \forall k,n \in \N: k \ge n.
	\end{equation}

The sequence $\ens{\varepsilon_k }_{k \ge 1}$ is bounded. So, there exists $k_0 \ge 1 $, depending on $\rho_1$ and $ \varepsilon_0 $, such that 
	$$
\vabsolut{\varepsilon_k } \le \varepsilon_0 \le \frac{ \rho_1 }{4} \paren{2k -1} \le \frac{ \rho_1 }{4} \paren{ k^2 -n^2}, \quad \forall k,n \ge 1: k \ge k_0, \ n \le k-1.
	$$ 
With this value of $k_0$ and~\eqref{f36} written for $k $ and $n$, with $k \ge k_0 $ and $n \le k -1 $, we obtain
	\begin{equation*}
	\left\{
	\begin{array}{l}
\displaystyle \lambda_{ k }^{(1)} - \lambda_{ n }^{(1)}  \ge \rho_1 \paren{k^2 - n^2} \ge \frac{\rho_1}{2} \paren{k^2 - n^2} , \\
	\noalign{\medskip}
\displaystyle \lambda_{ k }^{(2)} - \lambda_{ n }^{(1)}  \ge \lambda_{ k }^{(1)} - \lambda_{ n }^{(1)} - \varepsilon_{0 } \ge \rho_1 \paren{k^2 - n^2} - \frac{\rho_1}{4} \paren{k^2 - n^2} \ge \frac{\rho_1}{2} \paren{k^2 - n^2}, \\
	\noalign{\medskip}
\displaystyle \lambda_{ k }^{(1)}  - \lambda_{ n }^{(2)} \ge \lambda_{ k }^{(1)}  - \lambda_{ n }^{(1)} - \varepsilon_0 \ge \rho_1 \paren{k^2 - n^2} - \frac{\rho_1}{4} \paren{k^2 - n^2} \ge \frac{\rho_1}{2} \paren{k^2 - n^2} , \\
	\noalign{\medskip}
\displaystyle \lambda_{ k }^{(2)}  - \lambda_{ n }^{(2)} \ge \lambda_{ k }^{(1)}  - \lambda_{ n }^{(1)} - 2 \varepsilon_0 \ge  \rho_1 \paren{k^2 - n^2} - \frac{\rho_1}{2} \paren{k^2 - n^2} \ge \frac{\rho_1}{2} \paren{k^2 - n^2} ,
	\end{array}
	\right.
	\end{equation*}
i.e.,
	\begin{equation}\label{f37}
\displaystyle \lambda_{ k }^{(i)}  - \lambda_{ n }^{(j)} \ge \frac{\rho_1}{2} \paren{k^2 - n^2} , \quad \forall k,n \ge 1: k \ge k_0, \ n \le k-1, \quad \forall i,j \in \{ 1 , 2 \}. 
	\end{equation}

As a consequence of~\eqref{f37}, we also obtain $ \lambda_{ k }^{(1)}  < \lambda_{ k + 1}^{(2)} $ and $\lambda_{ k }^{(2)} < \lambda_{ k + 1}^{(1)}$, for any $k \ge k_0$.  This provides the following explicit formula for the terms of the increasing sequence $\Lambda$ when $k \ge 2 k_0 -1$:
	\begin{equation}\label{f40}
\Lambda_k = \left\{
	\begin{array}{ll}
\displaystyle \min \ens{\lambda_{ \ell }^{(1)}, \lambda_{ \ell }^{(2)}}, & \hbox{if } k = 2 \ell - 1, \\
	\noalign{\medskip}
\displaystyle \max \ens{\lambda_{ \ell }^{(1)}, \lambda_{ \ell }^{(2)}}, & \hbox{if } k = 2 \ell.
	\end{array}
	\right.
	\end{equation}

We are going to use the previous expression of the terms $\Lambda_k$, with $k \ge k_0$, in order to prove condition~\ref{item5} with $q=2$. Remember that the sequence $\Lambda$ is real and increasing. Then,
	$$
\frac{\Lambda_k - \Lambda_n}{k^2 - n^2 } > 0, \quad \forall k,n \ge 1: k \ge n + 1.
	$$
Assume that, for every $n \in \{ 1, \dots, 2 k_0 - 2 \} $ fixed, one has
	\begin{equation}\label{f39}
\liminf_{k \to \infty}  \frac{\Lambda_k - \Lambda_n}{k^2 - n^2 } \ge \frac{\rho_1}{4} \quad \hbox{and} \quad \limsup_{k \to \infty}  \frac{\Lambda_k - \Lambda_n}{k^2 - n^2 } \le \frac{\nu_1}{4}.
	\end{equation}
Then, there exists a positive constant $\widetilde \rho$, only depending on $k_0$ and $\rho_1$ or, equivalently, on $\rho_1$ and $ \varepsilon_0 $, such that
	\begin{equation*}
\Lambda_{ k} - \Lambda_{ n} \ge \widetilde \rho \paren{k^2 -n^2}  , \quad \forall k,n \in \N: 1 \le n \le 2k_0 -2 \hbox{ and } n \le k .
	\end{equation*}
In this way, we have proved condition~\ref{item5} for $q=1$ and $k,n \in \N$ such that $1 \le n \le 2k_0 -2$ and $ n \le k$. We will prove~\eqref{f39} below.


Let us now see that the sequence $\Lambda$ satisfies~\ref{item5}, with $q=2$ and an appropriate value of the parameter $\rho$, when $k, n \ge 1$ with $k \ge n + 2 $ and $n \ge 2 k_0 - 1$. We divide the proof into four cases:
\begin{enumerate}
\item Assume that $k = 2 \ell -1 $ and $n = 2m -1$, with $\ell , m \ge k_0$ and $k - n \ge 2$. In particular, $\ell - m \ge 1$, $\Lambda_k = \lambda_\ell^{(i)}$ and $\Lambda_n = \lambda_m^{(j)}$, with $i, j \in \{ 1, 2\}$. Thus, from~\eqref{f37}
	\begin{equation*}
	\begin{split}
\displaystyle \Lambda_k - \Lambda_n &= \lambda_\ell^{(i)} - \lambda_m^{(j)} \ge \frac{\rho_1}{2} \paren{\ell^2 - m^2} = \frac{\rho_1}{8} \left[ \paren{k +1}^2 - \paren{n+1}^2 \right] = \frac{\rho_1}{8} \paren{k + n + 2} \paren{k - n} \\
	\noalign{\medskip}
\displaystyle &\ge \frac{\rho_1}{8} \paren{k^2 - n^2}.
	\end{split}
	\end{equation*}
\item Assume now that $k = 2 \ell -1 $ and $n = 2m $, with $\ell , m \ge k_0$ and $k - n \ge 2$. In particular, $\ell - m \ge 3/2$, $\Lambda_k = \lambda_\ell^{(i)}$ and $\Lambda_n = \lambda_m^{(j)}$, with $i, j \in \{ 1, 2\}$, and we can apply~\eqref{f37}:
	\begin{equation*}
\displaystyle \Lambda_k - \Lambda_n = \lambda_\ell^{(i)} - \lambda_m^{(j)} \ge \frac{\rho_1}{2} \paren{\ell^2 - m^2} = \frac{\rho_1}{8} \left[ \paren{k +1}^2 - n^2 \right] \ge \frac{\rho_1}{8} \paren{k^2 - n^2}.
	\end{equation*}
\item If $k = 2 \ell $ and $n = 2m $, with $\ell , m \ge k_0$ and $k - n \ge 2$, then $\ell - m \ge 1$, $\Lambda_k = \lambda_\ell^{(i)}$ and $\Lambda_n = \lambda_m^{(j)}$, with $i, j \in \{ 1, 2\}$. Applying again~\eqref{f37}, we get
	\begin{equation*}
\Lambda_k - \Lambda_n = \lambda_\ell^{(i)} - \lambda_m^{(j)} \ge \frac{\rho_1}{2} \paren{\ell^2 - m^2} =  \frac{\rho_1}{8} \paren{k^2 - n^2}.
	\end{equation*}

\item In the case $k = 2 \ell $ and $n = 2m - 1$, with $\ell , m \ge k_0$ and $k - n \ge 2$ we will use the inequality
	$$
k^2 -(n + 1)^2 \ge \frac 12 \paren{k^2 - n^2}
	$$
which is valid for any $k, n \ge 1$ such that $k \ge n + 2$. Also, $\ell - m \ge 1/2$, i.e., $\ell - m \ge 1$ and we can apply~\eqref{f37}. As before, $\Lambda_k = \lambda_\ell^{(i)}$ and $\Lambda_n = \lambda_m^{(j)}$, with $i, j \in \{ 1, 2\}$, and
	\begin{equation*}
\Lambda_k - \Lambda_n = \lambda_\ell^{(i)} - \lambda_m^{(j)} \ge \frac{\rho_1}{2} \paren{\ell^2 - m^2} =  \frac{\rho_1}{8} \paren{k^2 - \paren{n + 1}^2} \ge \frac{\rho_1}{16} \paren{k^2 - n^2} .
	\end{equation*}
We can conclude that property~\ref{item5} holds for the sequence $\Lambda$ with $q = 2$ and 
	$$
\rho = \min \ens{\widetilde \rho, \frac{\rho_1}{16}}. 
	$$
Remember that the constant $\widetilde \rho$ only depends on $\rho_1$  and $ \varepsilon_0 $. Therefore, $\rho$ only depends on $\rho_1$  and $ \varepsilon_0 $.
\end{enumerate}

The next task will be the proof of~\eqref{f39}. To this end, let us fix $n$ such that $ 1 \le n \le 2k_0 -2 $ and $k \ge 2k_0 -1 $. Then $k = 2 \ell$ or $k = 2 \ell -1$ with $\ell \ge k_0$. In both cases, $\Lambda_k = \lambda_{ \ell }^{(i)} $, with $ i \in \{ 1, 2\} $, and we can write (see~\eqref{f36}):
	$$
	\left\{
	\begin{array}{l}
\displaystyle \Lambda_k - \Lambda_n \ge \lambda_{ \ell }^{(1)} - \varepsilon_0 - \Lambda_{2 k_0 -2 } \ge \rho_1 \paren{\ell^2 - 1} + \lambda_1^{(1)} - \varepsilon_0 - \Lambda_{2 k_0 -2 } \\
	\noalign{\smallskip}
\displaystyle \phantom{ \Lambda_k - \Lambda_n} \ge \rho_1 \paren{\frac{k^2}{4} - 1}  - \varepsilon_0 - \Lambda_{2 k_0 -2 }, \\
	\noalign{\smallskip}
\displaystyle \Lambda_k - \Lambda_n \le \lambda_{ \ell }^{(1)} + \varepsilon_0 - \Lambda_{1 } \le \nu_1 \paren{\ell^2 - 1} + \lambda_1^{(1)} + \varepsilon_0 - \Lambda_{1 } \\
	\noalign{\smallskip}
\displaystyle \phantom{ \Lambda_k - \Lambda_n}  \le \nu_1 \paren{\frac{(k + 1)^2}{4} - 1} + \lambda_1^{(1)} + \varepsilon_0 - \Lambda_{1 }.
	\end{array}
	\right.
	$$
This proves~\eqref{f39}.

In order to finish the proof of Proposition~\ref{p5}, let us check that the sequence $\Lambda$ fulfills condition~\eqref{item6} for an appropriate $\nu >0$. The proof is very close to that of condition~\ref{item5}. First, one has
	$$
\vabsolut{\varepsilon_k } \le \varepsilon_0 \le \varepsilon_0 \paren{2k -1} \le \varepsilon_0 \paren{ k^2 -n^2}, \quad \forall k,n \ge 1:  n \le k-1.
	$$ 
From this inequality and the second inequality in~\eqref{f36} we deduce, for instance
	$$
 \lambda_{ k }^{(1)} - \lambda_{ n }^{(2)}  \le \lambda_{ k }^{(1)} - \lambda_{ n }^{(1)} + \varepsilon_{0 } \le \paren{ \nu_1 + \varepsilon_0} \paren{k^2 - n^2} , \quad \forall k,n \ge 1: n \le k-1.
	$$
Thus, as before, it is not difficult to show the following inequalities: 
	\begin{equation}\label{f38}
	\left\{
	\begin{array}{l}
\displaystyle \lambda_{ k }^{(1)} - \lambda_{ n }^{(1)}  \le \nu_1 \paren{k^2 - n^2} \le \paren{ \nu_1 + \varepsilon_0} \paren{k^2 - n^2}, \quad \forall k,n \ge1 :  n \le k, \\
	\noalign{\medskip}
\displaystyle \lambda_{ k }^{(1)}  - \lambda_{ n }^{(2)} \le \paren{ \nu_1 + \varepsilon_0} \paren{k^2 - n^2} , \quad \forall k,n \ge 1: n \le k-1,\\
	\noalign{\medskip}
\displaystyle \lambda_{ k }^{(2)}  - \lambda_{ n }^{(1)} \le \paren{ \nu_1 + \varepsilon_0} \paren{k^2 - n^2} , \quad \forall k,n \ge 1: n \le k-1 , \\
	\noalign{\medskip}
\displaystyle \lambda_{ k }^{(2)}  - \lambda_{ n }^{(2)} \le \paren{ \nu_1 + \varepsilon_0} \paren{k^2 - n^2} , \quad \forall k,n \ge 1: n \le k-1 .
	\end{array}
	\right.
	\end{equation}

Let us now prove condition~\eqref{item6} for the sequence $\Lambda$. As before, from the second property in~\eqref{f39} we deduce the existence of a positive constant $\widetilde \nu$, only depending on $k_0 $ and $\nu_1$, such that
	\begin{equation*}
\Lambda_{ k} - \Lambda_{ n} \le \widetilde \nu \paren{k^2 -n^2}  , \quad \forall k,n \in \N: 1 \le n \le 2k_0 -2 \hbox{ and } n \le k .
	\end{equation*}

Let us now see inequality~\eqref{item6} when $k,n \in \N$ are such that $2k_0 -1 \le n \le k$. Remember that, in this case, we have an explicit formula of the terms of the sequence $\Lambda$ (see~\eqref{f40}). Let us first consider the case $ n \ge 2k_0 -1$ and $k = n +1$. Thus, 
	$$
	\left\{
	\begin{array}{l}
\displaystyle \Lambda_{n + 1} - \Lambda_{n } = \Lambda_{2 \ell} - \Lambda_{2\ell - 1} = \vabsolut{\varepsilon_\ell } \le \varepsilon_0 \paren{ \paren{ n + 1}^2 - n^2 } , \\
	\noalign{\smallskip}
\displaystyle \Lambda_{n + 1} - \Lambda_{n } = \Lambda_{2 \ell + 1} - \Lambda_{2\ell} \le \lambda_{\ell + 1}^{(1)} - \lambda_\ell^{(1)} \le \nu_1 \paren{2 \ell + 1} \le \nu_1 \paren{4 \ell + 1} = \nu_1  \paren{ \paren{n+1}^2 - n^2 }. 
	\end{array}
	\right.
	$$

In the general case, i.e., when $k,n \in \N$ are such that $2k_0 -1 \le n \le k$ with $k \ge n+2$, we can repeat the arguments above and deduce inequality~\eqref{item6} . Indeed, as a consequence of~\eqref{f38}, we deduce
\begin{enumerate}
\item If $k = 2 \ell -1 $ and $n = 2m -1$, with $\ell , m \ge k_0$ and $k - n \ge 2$, then, $\ell - m \ge 1$, $\Lambda_k = \lambda_\ell^{(i)}$ and $\Lambda_n = \lambda_m^{(j)}$, with $i, j \in \{ 1, 2\}$. \eqref{f38} implies 
	\begin{equation*}
	\begin{split}
\displaystyle \Lambda_k - \Lambda_n &= \lambda_\ell^{(i)} - \lambda_m^{(j)} \le \paren{ \nu_1 + \varepsilon_0} \paren{\ell^2 - m^2} = \frac{ \nu_1 + \varepsilon_0}{4} \left[ \paren{k +1}^2 - \paren{n+1}^2 \right]  \\
	\noalign{\medskip}
\displaystyle &= \frac{ \nu_1 + \varepsilon_0}{4} \paren{k + n + 2} \paren{k - n} \le \frac{ \nu_1 + \varepsilon_0}{2} \paren{k^2 - n^2}.
	\end{split}
	\end{equation*}
\item Assume now that $k = 2 \ell -1 $ and $n = 2m $, with $\ell , m \ge k_0$ and $k - n \ge 2$. In this case, $\ell - m \ge 3/2$, $\Lambda_k = \lambda_\ell^{(i)}$ and $\Lambda_n = \lambda_m^{(j)}$, with $i, j \in \{ 1, 2\}$. On the other hand, it is not difficult to check that
	$$
\paren{k +1}^2 - n^2 \le 2 \paren{k^2 - n^2}, \quad \forall k,n \ge 1: k \ge n + 2. 
	$$
Thus, from~\eqref{f38} we get:
	\begin{equation*}
\displaystyle \Lambda_k - \Lambda_n = \lambda_\ell^{(i)} - \lambda_m^{(j)} \le \paren{ \nu_1 + \varepsilon_0}  \paren{\ell^2 - m^2} = \frac{ \nu_1 + \varepsilon_0}{4} \left[ \paren{k +1}^2 - n^2 \right] \le \frac{ \nu_1 + \varepsilon_0}{2} \paren{k^2 - n^2}.
	\end{equation*}
\item When $k = 2 \ell $ and $n = 2m $, with $\ell , m \ge k_0$ and $k - n \ge 2$, one has $\ell - m \ge 1$, $\Lambda_k = \lambda_\ell^{(i)}$ and $\Lambda_n = \lambda_m^{(j)}$, with $i, j \in \{ 1, 2\}$. Applying again~\eqref{f38}, we get
	\begin{equation*}
\Lambda_k - \Lambda_n = \lambda_\ell^{(i)} - \lambda_m^{(j)} \le \paren{ \nu_1 + \varepsilon_0}  \paren{\ell^2 - m^2} =  \frac{ \nu_1 + \varepsilon_0}{4} \paren{k^2 - n^2}.
	\end{equation*}

\item Finally, let us take $k = 2 \ell $ and $n = 2m - 1$, with $\ell , m \ge k_0$ and $k - n \ge 2$. Again, $\ell - m \ge 1/2$, i.e., $\ell - m \ge 1$ and we can apply~\eqref{f38}. As in the previous cases, $\Lambda_k = \lambda_\ell^{(i)}$ and $\Lambda_n = \lambda_m^{(j)}$, with $i, j \in \{ 1, 2\}$, and
	\begin{equation*}
\Lambda_k - \Lambda_n = \lambda_\ell^{(i)} - \lambda_m^{(j)} \le \paren{ \nu_1 + \varepsilon_0} \paren{\ell^2 - m^2} =  \frac{ \nu_1 + \varepsilon_0}{4} \paren{k^2 - \paren{n + 1}^2} \le \frac{ \nu_1 + \varepsilon_0}{4} \paren{k^2 - n^2}. 
	\end{equation*}

Summarizing, we have prove property~\ref{item6} for the sequence $\Lambda$ with  
	$$
\nu = \max \ens{\widetilde \nu, \frac{ \nu_1 + \varepsilon_0}{2}  }. 
	$$
Remember again that the constant $\widetilde \nu$ only depends on $k_0$ and $\nu_1$, that is to say, on $\rho_1$, $\varepsilon_0$ and $\nu_1$. Therefore, the paremeter $\nu$ only depends on $\rho_1$, $\nu_1$ and $\varepsilon_0$.
\end{enumerate}
With the proof of property~\eqref{item6} we end the proof of Proposition~\ref{p5}.  \fin
%


\end{document}